\noindent \centerline{\bf MITTAG-LEFFLER FUNCTIONS AND THEIR APPLICATIONS} 

\vskip.3cm \noindent \centerline{H.J. HAUBOLD}
\vskip0cm \noindent\centerline{Office for Outer Space Affairs, United Nations} 
\vskip0cm \centerline{Vienna International Centre, P.O. Box 500, A-1400 Vienna, Austria.}
\vskip0cm \centerline{and Centre for Mathematical Sciences, Pala Campus,}
\vskip0cm \noindent\centerline{Arunapuram P.O., Pala, Kerala-686574, India.}
\vskip0cm \centerline{hans.haubold@unvienna.org}

\vskip.3cm \noindent \centerline{A.M. MATHAI}
\vskip0cm \centerline{Department of Mathematics and Statistics, McGill University}
\vskip0cm \centerline{Montreal, Canada H3A 2K6}
\vskip0cm \centerline{and Centre for Mathematical Sciences, Pala Campus,}
\vskip0cm \noindent\centerline{Arunapuram P.O., Pala, Kerala-686574, India.}
\vskip0cm \centerline{mathai@math.mcgill.ca}

\vskip.2cm \noindent\centerline{R.K. SAXENA}
\vskip0cm \noindent\centerline{Department of Mathematics and Statistics, Jai Narain Vyas University}
\vskip0cm \noindent\centerline{Jodhpur-342005, India.}
\vskip0cm \centerline{ram.saxena@yahoo.com}

\vskip.5cm
\noindent
{\bf Abstract}

\vskip.5cm Motivated essentially by the success of the
applications of the Mittag-Leffler functions in many areas of
science and engineering, the authors present in a unified manner, a
detailed account or rather a brief survey of the Mittag-Leffler
function, generalized Mittag-Leffler functions, Mittag-Leffler type
functions, and  their interesting and useful properties.
Applications of Mittag-Leffler functions in certain areas of
physical and applied sciences are also demonstrated. During the last
two decades this function has come into prominence after about nine
decades of its discovery by a Swedish Mathematician G.M.
Mittag-Leffler, due its vast potential of its applications in
solving the problems of physical, biological, engineering and earth
sciences etc. In this survey paper, nearly all types of
Mittag-Leffler type functions existing in the literature are
presented. An attempt is made to present nearly an exhaustive list
of references concerning the Mittag-Leffler functions to make the
reader familiar with the present trend of research in Mittag-Leffler 
type functions and their applications.

\vskip.5cm \noindent {\bf 1.\hskip.3cm Introduction}

\vskip.5cm The special function

$$E_{\alpha}(z)=\sum_{k=0}^{\infty}{{z^k}\over{\Gamma(1+\alpha k)}},
\alpha\in C, \Re(\alpha)>0, z\in C\eqno(1.1)
$$and its general form

$$E_{\alpha,\beta}(z)=\sum_{k=0}^{\infty}{{z^k}\over{\Gamma(\beta+\alpha
k)}}, \alpha,\beta \in C, \Re(\alpha)>0, \Re(\beta)>0, z\in
C\eqno(1.2)
$$with $C$ being the set of complex numbers are called
Mittag-Leffler functions (Erd\'elyi, et al., 1955, Section 18.1).
The former was introduced by Mittag-Leffler (Mittag-Leffler,1903),
in connection with his method of summation of some divergent series.
In his papers (Mittag-Leffler, 1903, 1905), he investigated certain
properties of this function.  The function defined by (1.2) first
appeared in the work of Wiman (Wiman, 1905). The function (1.2) is
studied, among others, by Wiman (1905), Agarwal (1953), Humbert
(1953) and Humbert and Agarwal (1953) and others.  The main
properties of these functions are given in the book by Erd\'elyi et
al. (1955, Section 18.1) and a more comprehensive  and  a detailed
account of  Mittag-Leffler functions are presented in Dzherbashyan
(1966,  Chapter 2).  In particular, the functions (1.1) and (1.2)
are entire functions of order $\rho=1/{\alpha}$ and type $\sigma
=1$; see, for example, (Erd\'elyi, et al., 1955, p.118).

\vskip.2cm The Mittag-Leffler function arises naturally in the
solution of fractional order integral equations or fractional order
differential equations, and especially in the investigations of the
fractional generalization of the kinetic equation, random walks,
L\'evy flights, super-diffusive transport and in the study of
complex systems. The ordinary and generalized Mittag-Leffler
functions interpolate between a purely exponential law and power-law
like behavior of phenomena governed by ordinary kinetic equations
and their fractional counterparts, see Lang (1999a, 1999b), Hilfer
(2000), Saxena et al. (2002). \vskip.2cm The Mittag-Leffler function
is not given in the tables of Laplace transforms, where it naturally
occurs in the derivation of the inverse Laplace transform of the
functions of the type $p^{\alpha}(a+bp^{\beta})$, where $p$ is the
Laplace transform parameter and $a$ and $b$ are constants. This
function also occurs in the solution of certain boundary value
problems involving fractional integro-differential equations of
Volterra type (Samko et al., 1993). During the various developments
of fractional calculus in the last four decades this function has
gained importance and popularity on account of its vast applications
in the fields of science and engineering. Hille and Tamarkin (1930)
have presented a solution of the Abel-Volterra type equation in
terms of Mittag-Leffler function. During the last 15 years the
interest in Mittag-Leffler function and Mittag-Leffler type
functions is considerably increased among engineers and scientists
due to their vast potential of applications in several applied
problems, such as fluid flow, rheology, diffusive transport akin to
diffusion, electric networks, probability, statistical distribution
theory etc. For a detailed account of various properties,
generalizations, and application of this function, the reader may
refer to earlier important works of Blair (1974), Bagley and Torvik
(1984), Caputo and Mainardi (1971), Dzherbashyan (1966), Gorenflo
and Vessella (1991), Gorenflo and Rutman (1994), Kilbas and Saigo
(1995), Gorenflo et al. (1997),  Gorenflo and Mainardi (1994, 1996,
1997), Gorenflo, Luchko and Rogosin (1997),
 Gorenflo, Kilbas and Rogosin (1998), Luchko (1999), Luchko and
Srivastava (1995), Kilbas, Saigo and Saxena (2002, 2004), Saxena and
Saigo (2005), Kiryakova (2008a, 2008b), Saxena, Kalla and Kiryakova
(2003), Saxena, Mathai and Haubold (2002, 2004, 2004a, 2004b, 2006),
Saxena and Kalla (2008), Mathai, Saxena and Haubold (2006), Haubold
and Mathai (2000), Haubold, Mathai and Saxena (2007), Srivastava and
Saxena (2001), and others.

\vskip.2cm This paper is organized as follows: Section 2 deals with
special cases of $E_{\alpha}(z)$. Functional relations of
Mittag-Leffler functions are presented in Section 3. Section 4 gives
the basic properties. Section 5 is devoted to the derivation of
recurrence relations for Mittag-Leffler functions. In Section 6,
asymptotic expansions of the Mittag-Leffler functions are given.
Integral representations of Mittag-Leffler functions are given in
Section 7. Section 8 deals with the H-function and its special
cases. The Mellin-Barnes integrals for the Mittag-Leffler functions
are established in Section 9. Relations of Mittag-Leffler functions
with Riemann-Liouville fractional calculus operators are derived in
Section 10. Generalized Mittag-Leffler functions and some of their
properties are given in Section 11. Laplace transform, Fourier
transform, and fractional integrals and derivatives are discussed in
Section 12. Section 13 is devoted to the  application of
Mittag-Leffler function in fractional kinetic equations. In Section
14, time-fractional diffusion equation is solved. Solution of
space-fractional diffusion equation is discussed in Section 15. In
Section 16, solution of a fractional reaction-diffusion equation is
investigated in terms of the H-function. Section 17 is devoted to
the application of generalized Mittag-Leffler functions in nonlinear
waves. Recent generalizations of Mittag-Leffler functions are
discussed in Section 18.

\vskip.5cm \noindent {\bf 2.\hskip.3cm Some special cases}

\vskip.3cm We begin our study by giving the special cases of the
Mittag-Leffler function $E_{\alpha}(z)$.

$$E_0(z)={{1}\over{1-z}}, |z|<1\leqno(i)$$
$$E_1(z)={\rm e}^z\leqno(ii)$$
$$E_2(z)=\cosh (\sqrt{z}), z\in C\leqno(iii)$$
$$E_2(-z^2)=\cos z, z\in C\leqno(iv)$$
$$E_3(z)={1\over2}\left[{\rm e}^{z^{1\over3}}+2{\rm
e}^{-{1\over2}z^{1\over3}}\cos\left({{\sqrt{3}}\over{2}}z^{1\over3}\right)\right],
z\in C\leqno(v)$$
$$E_4(z)={1\over2}\left[\cos (z^{1\over4})+\cosh
(z^{1\over4})\right], z\in C\leqno(vi)$$
$$E_{1\over2}(\pm z^{1\over2})={\rm e}^z\left[1+{\rm
erf}(\pm z^{1\over2})\right]={\rm e}^z{\rm erfc}(\mp z^{1\over2}),
z\in C\leqno(vii)$$where erfc denotes the complimentary error
function and the error function is defines as

$${\rm erf}(z)={{2}\over{\sqrt{\pi}}}\int_0^z\exp (-t^2){\rm d}t,~~
{\rm erfc}(z)=1-{\rm erf}(z), z\in C.
$$For half-integer $n/2$ the function can be written explicitly as

$$E_{{{n}\over2}}(z)={_0F_{n-1}}\left(~~:{{1}\over{n}},{{2}\over{n}},...,{{n-1}\over{n}};{{z^2}\over{n^n}}\right)\leqno(viii)$$
$$+{{2^{(n+1)/2}z}\over{n!\sqrt{\pi}}}{_1F_{2n-1}}\left(1;{{n+2}\over{2n}},{{n+3}\over{2n}},...,{{3n}\over{2n}};{{z^2}\over{n^n}}\right)$$
$$E_{1,2}(z)={{{\rm e}^{z}-1}\over{z}},~~E_{2,2}(z)={{\sinh
(\sqrt{z})}\over{\sqrt{z}}}.\leqno(ix)$$

\vskip.3cm \noindent{\bf 3.\hskip.3cm Functional relations for the
Mittag-Leffler functions}

\vskip.3cm In this section, we discuss the Mittag-Leffler functions
of rational order $\alpha=m/n$, with $m,n\in N$ relatively prime.
The differential and other properties of these functions are
described in Erd\'elyi, et al.(1955) and Dzherbashyan (1966).

\vskip.3cm \noindent{\bf Theorem 3.1}.\hskip.3cm {\it The following
results hold:

$$\eqalignno{{{{\rm d}^m}\over{{\rm
d}z^m}}E_m(z^m)&=E_m(z^m)&(3.1)\cr {{{\rm d}^m}\over{{\rm
d}z^m}}E_{{{m}\over{n}}}(z^{{m}\over{n}})&=E_{{m}\over{n}}(z^{{m}\over{n}})+\sum_{r=1}^{n-1}{{z^{-{{rm}\over{n}}}}\over{\Gamma(1-rm/n)}},
~n=2,3,...&(3.2)\cr
E_{{m}\over{n}}(z)&={{1}\over{m}}\sum_{r=1}^{m-1}E_{{1}\over{n}}(z^{{1}\over{m}}\exp(i2\pi
r/m))&(3.3)\cr E_{{1}\over{n}}(z^{{1}\over{n}})&={\rm
e}^z\left[1+\sum_{r=1}^{n-1}{{\gamma(1-{{r}\over{n}},z)}\over{\Gamma(1-r/n)}}\right],~n=2,3,...&(3.4)\cr}
$$where $\gamma(a,z)$ denotes the incomplete gamma function, defined
by,

$$\gamma(a,z)=\int_0^z{\rm e}^{-t}t^{a-1}{\rm d}t.\eqno(3.5)
$$}In order to establish the above formulas, we observe that (3.1)
and (3.2) readily follow from the definition (1.2). For proving the
formula (3.3), we recall the identity

$$\sum_{r=0}^{m-1}\exp[i2\pi kr/m]=\cases{m\hbox{  if  }k=0\hbox{
(mod $m$)}\cr 0\hbox{  if  }k\ne 0\hbox{  (mod $m$)}.\cr}\eqno(3.6)
$$By virtue of the results (1.1) and (3.6), we find that

$$\eqalignno{\sum_{r=0}^{m-1}E_{\alpha}(z{\rm e}^{i2\pi
r/m})&=mE_{\alpha m}(z^m), m\in N&(3.7)\cr \noalign{\hbox{which can
be written as}}
E_{\alpha}(z)&={{1}\over{m}}\sum_{r=0}^{m-1}E_{{\alpha}\over{m}}(z^{{1}\over{m}}{\rm
e}^{i2\pi r/m}),~m\in N&(3.8)\cr}
$$and the result (3.3) now follows by taking $\alpha=m/n$. To prove
the relation (3.4), we set $m=1$ in (3.1) and multiply it by $\exp
(-z)$ to obtain

$${{{\rm d}}\over{{\rm d}z}}\left[{\rm
e}^{-z}E_{{1}\over{n}}(z^{{1}\over{n}})\right]={\rm
e}^{-z}\sum_{r=1}^{m-1}{{z^{-r/m}}\over{\Gamma(1-r/m)}}.\eqno(3.9)
$$On integrating both sides of the above equation with respect to $z$
and using the definition of incomplete gamma function (3.5), we
obtain the desired result (3.4). An interesting case of (3.8) is
given by
$$E_{2\alpha}(z^2)={1\over2}[E_{\alpha}(z)+E_{\alpha}(-z)].\eqno(3.10)
$$

\vskip.3cm \noindent {\bf 4.\hskip.3cm Basic properties}

\vskip.3cm This section is based on the paper of Berberan-Santos
(2005). From (1.1) and (1.2) it is not difficult to prove that

$$\eqalignno{E_{\alpha}(-x)&=E_{2\alpha}(x^2)-xE_{2\alpha,1+\alpha}(x^2)&(4.1)\cr
\noalign{\hbox{and}} E_{\alpha}(-ix)&=E_{2x}(-x^2)-ix
E_{2\alpha,1+\alpha}(-x^2), ~i=\sqrt{-1}.&(4.2)\cr}
$$It is shown in Berberan-Santos (2005, p.631) that the following
three equations can be used for the direct inversion of a function
$I(x)$ to obtain its inverse $H(k)$:

$$\eqalignno{H(k)&={{{\rm
e}^{ck}}\over{\pi}}\int_0^{\infty}[\Re[I(c+i\omega)]\cos(k\omega)-\Im[I(c+i\omega)]\sin(k\omega){\rm
d}\omega&(4.3)\cr &={{2{\rm e}^{ck}}\over{\pi}}\int_0^{\infty}
\Re[I(c+i\omega)\cos(k\omega){\rm d}\omega,~k>0&(4.4)\cr &=-{{2{\rm
e}^{ck}}\over{\pi}}\int_0^{\infty}\Im[I(c+i\omega)]\sin(k\omega){\rm
d}\omega,~~k>0.&(4.5)\cr}
$$With the help of the results (4.2) and (4.4), it yields the
following formula for the inverse Laplace transform $H(k)$ of the
function $E_{\alpha}(-x)$.

$$H_{\alpha}(k)={{2}\over{\pi}}\int_0^{\infty}E_{2\alpha}(-t^2)\cos(kt){\rm
d}t, ~k>0, 0\le \alpha\le 1.\eqno(4.6)
$$In particular, the following interesting results can be derived
from the above result.

$$\eqalignno{H_1(k)&={{2}\over{\pi}}\int_0^{\infty}\cosh (it)\cos
(kt){\rm d}t={{2}\over{\pi}}\int_0^{\infty}\cos(t)\cos(kt){\rm
d}t=\delta(k-1),i=\sqrt{-1}&(4.7)\cr
H_{{1}\over{2}}(k)&={{2}\over{\pi}}\int_0^{\infty}{\rm
e}^{-t^2}\cos(kt){\rm d}t={{1}\over{\sqrt{\pi}}}{\rm
e}^{-{{k^2}\over4}}&(4.8)\cr
H_{1\over4}(k)&={{2}\over{\pi}}\int_0^{\infty}{\rm e}^{t^2}{\rm
erfc}(t^2)\cos(kt){\rm d}t.&(4.9)\cr}
$$Another integral representation of $H_{\alpha}(k)$ in terms of the
L\'evy one-sided stable distribution $L_{\alpha}(k)$ was given by
Pllard (1948) in the form

$$H_{\alpha}(k)={{1}\over{\alpha}}k^{-(1+{{1}\over{\alpha}})}L_{\alpha}(k^{-{{1}\over{\alpha}}}).\eqno(4.10)
$$The inverse Laplace transform of $E_{\alpha}(-x^{\beta})$, denoted
by $H_{\alpha}^{\beta}(k)$ with $0< \alpha \le 1$, is obtained as

$$H_{\alpha}^{\beta}(k)=\int_0^{\infty}t^{{\alpha}\over{\beta}}L_{\alpha}(t)L_{\beta}(kt^{{\alpha}\over{\beta}}){\rm
d}t,\eqno(4.11)
$$where $L_{\alpha}(t)$ is the one-sided L\'evy probability density
function. From Berberan-Santos (2005, p.432) we have

$$H_{\alpha}(k)={{1}\over{\pi}}\int_0^{\infty}[E_{2\alpha}(-\omega^2)\cos(k\omega)+\omega
E_{2\alpha,1+\alpha}(-\omega^2)\sin (k\omega)]{\rm d}\omega, 0<
\alpha\le 1.\eqno(4.12)
$$Expanding the above equation in a power series, it gives

$$\eqalignno{H_{\alpha}(k)&={{1}\over{\pi}}\sum_{n=0}^{\infty}b_n(\alpha)k^{n},~0\le
\alpha<1&(4.13)\cr \noalign{\hbox{with}}
b_0(\alpha)&=\int_0^{\infty}E_{2\alpha}(-t^2){\rm d}t.&(4.14)\cr}
$$The Laplace transform of the equation (4.13) is the asymptotic
expansion of $E_{\alpha}(-x)$ as

$$E_{\alpha}(-x)={{1}\over{\pi}}\sum_{n=0}^{\infty}{{b_n(\alpha)}\over{x^{n+1}}},~0\le
\alpha<1.\eqno(4.15)
$$

\vskip.3cm \noindent {\bf 5.\hskip.3cm Recurrence relations}

\vskip.3cm By virtue of the definition (1.2), the following
relations are obtained in the form of

\vskip.3cm \noindent {\bf Theorem 5.1.}\hskip.3cm{\it We have

$$\eqalignno{E_{\alpha,\beta}(z)&=zE_{\alpha,\alpha+\beta}(z)+{{1}\over{\Gamma(\beta)}}&(5.1)\cr
E_{\alpha,\beta}(z)&=\beta E_{\alpha,\beta+1}(z)+\alpha z{{{\rm
d}}\over{{\rm d}z}}E_{\alpha,\beta+1}(z)&(5.2)\cr {{{\rm
d}^m}\over{{\rm
d}z^m}}\left[z^{\beta-1}E_{\alpha,\beta}(z^{\alpha})\right]&=z^{\beta-m-1}E_{\alpha,\beta-m}(z^{\alpha}),\Re(\beta
-m)>0,m\in N.&(5.3)\cr {{{\rm d}}\over{{\rm
d}z}}E_{\alpha,\beta}(z)&={{E_{\alpha,\beta-1}(z)-(\beta-1)E_{\alpha,\beta}(z)}\over{\alpha
z}}.&(5.4)\cr}
$$}The above formulae are useful in computing the derivative of the
Mittag-Leffler function $E_{\alpha,\beta}(z)$. The following theorem
has been established by Saxena (2002):

\vskip.3cm \noindent {\bf Theorem 5.2}.\hskip.3cm{\it If
$\Re(\alpha)>0, \Re(\beta)>0$ and $r\in N$ then there holds the
formula

$$z^rE_{\alpha,\beta+r\alpha}(z)=E_{\alpha,\beta}(z)-\sum_{n=0}^{r-1}{{z^n}\over{\Gamma(\beta+n\alpha)}}.\eqno(5.5)
$$}

\vskip.3cm \noindent {\bf Proof}.\hskip.3cm We have from the right
side of (5.5),

$$E_{\alpha,\beta}(z)-\sum_{n=0}^{r-1}{{z^n}\over{\Gamma(\beta+n\alpha)}}
=\sum_{n=r}^{\infty}{{z^n}\over{\Gamma(\beta+n\alpha)}}.
$$Put $n-r=k$ or $n=k+r$. Then

$$\eqalignno{\sum_{n=r}^{\infty}{{z^n}\over{\Gamma(\beta+n\alpha)}}
&=\sum_{k=0}^{\infty}{{z^{k+r}}\over{\Gamma(\beta+r\alpha+k\alpha)}}\cr
&=z^rE_{\alpha,\beta+r\alpha}(z).\cr}
$$For $r=2,3,4$ we obtain the following corollaries:

\vskip.3cm \noindent {\bf Corollary 5.1.}\hskip.3cm{\it If
$\Re(\alpha)>0, \Re(\beta)>0$ then there holds the formula
$$z^2E_{\alpha,\beta+2\alpha}(z)=E_{\alpha,\beta}(z)-{{1}\over{\Gamma(\beta)}}-{{z}\over{\Gamma(\alpha+\beta)}}.\eqno(5.6)
$$}{\bf Corollary 5.2}.\hskip.3cm{\it If
$\Re(\alpha)>0,\Re(\beta)>0$ then there holds the formula

$$z^3E_{\alpha,\beta+3\alpha}(z)=E_{\alpha,\beta}(z)-{{1}\over{\Gamma(\beta)}}-{{z}\over{\Gamma(\alpha+\beta)}}-{{z^2}\over{\Gamma(2\alpha+\beta)}}.\eqno(5.7)
$$}{\bf Corollary 5.3.}\hskip.3cm{\it If $\Re(\alpha)>0,
\Re(\beta)>0$ then there holds the formula

$$z^4E_{\alpha,\beta+4\alpha}(z)=E_{\alpha,\beta}(z)-{{1}\over{\Gamma(\beta)}}-{{z}\over{\Gamma(\alpha+\beta)}}-{{z^2}\over{\Gamma(2\alpha+\beta)}}
-{{z^3}\over{\Gamma(3\alpha+\beta)}}.\eqno(5.8)
$$}

\vskip.3cm \noindent{\bf Remark 5.1.}\hskip.3cm For a generalization
of the result  (5.5), see Saxena, Kalla and Kiryakova (2003).

\vskip.5cm \noindent{\bf 6.\hskip.3cm Asymptotic expansions}

\vskip.3cm The asymptotic behavior of Mittag-Leffler functions plays
a very important role in the interpretation of the solution of
various problems of physics connected with fractional reaction,
fractional relaxation, fractional diffusion and fractional
reaction-diffusion etc in complex systems. The asymptotic expansion
of $E_{\alpha}(z)$ is based on the integral representation of the
Mittag-Leffler function in the form

$$E_{\alpha}(z)={{1}\over{2\pi
i}}\int_{\Omega}{{t^{\alpha-1}\exp(t)}\over{t^{\alpha}-z}}{\rm d}t,
\Re(\alpha)>0, \alpha, z\in C,\eqno(6.1)
$$where the path of integration $\Omega$ is a loop starting and
ending at $-\infty$ and encircling the circular disk $|t|\le
|z|^{{1}\over{\alpha}}$ in the positive sense, $|\arg t|<\pi$ on
$\Omega$. The integrand has a branch point at $t=0$. The complex
$t$-plane is cut along the negative real axis and in the cut plane
the integrand is single-valued, the principal branch of $t^{\alpha}$
is taken in the cut plane. (6.1) can be proved by expanding the
integrand in powers of $t$ and integrating term by term by making
use of the well-known Hankel's integral for the reciprocal of the
gamma function, namely

$${{1}\over{\Gamma(\beta)}}={{1}\over{2\pi i}}\int_{H_a}{{{\rm
e}^{\zeta}}\over{\zeta^{\beta}}}{\rm d}\zeta.\eqno(6.2)
$$The integral representation (6.1) can be used to obtain the
asymptotic expansion of the Mittag-Leffler function at infinity
(Erd\'elyi, et at., 1955). Accordingly, the following cases are
mentioned below: $(i)$:~~If $0<\alpha<2$ and $\mu$ is a real number
such that

$${{\pi\alpha}\over{2}}<\mu<\min[\pi,\pi\alpha],\eqno(6.3)
$$then for $N^{*}\in N, N^{*}\ne 1$ there holds the following
asymptotic expansion:

$$E_{\alpha}(z)={{1}\over{\alpha}}z^{(1-\beta)/\alpha}\exp(z^{{1}\over{\alpha}})-\sum_{r=1}^{N^{*}}{{1}\over{\Gamma(1-\alpha r)}}
{{1}\over{z^r}}+O\left[{{1}\over{z^{N^{*}+1}}}\right],\eqno(6.4)
$$as $|z|\rightarrow\infty$, $|\arg z|\le \mu$; and

$$E_{\alpha}(z)=-\sum_{r=1}^{N^{*}}{{1}\over{\Gamma(1-\alpha r)}}{{1}\over{z^r}}+O\left[{{1}\over{z^{N^{*}+1}}}\right],\eqno(6.5)
$$as $|z|\rightarrow\infty, \mu\le |\arg z|\le \pi$.  $(ii):$~~When
$\alpha\ge 2$ then there holds the following asymptotic expansion:

$$E_{\alpha}(z)={{1}\over{\alpha}}\sum_nz^{{1}\over{n}}\exp[\exp({{2n\pi
i}\over{\alpha}})z^{{1}\over{\alpha}}]-\sum_{r=1}^{N^{*}}{{1}\over{\Gamma(1-\alpha
r)}}{{1}\over{z^r}}+O\left[{{1}\over{z^{N^{*}+1}}}\right]\eqno(6.6)
$$as $|z|\rightarrow\infty, |\arg z|\le {{\alpha \pi}\over{2}}$, and
where the first sum is taken over all integers $n$ such that

$$|\arg (z)+2\pi n|\le{{\alpha\pi}\over{2}}.\eqno(6.7)
$$The asymptotic expansion of $E_{\alpha,\beta}(z)$ is based on the
integral representation of the Mittag-Leffler function
$E_{\alpha,\beta}(z)$ in the form

$$E_{\alpha,\beta}(z)={{1}\over{2\pi
i}}\int_{\Omega}{{t^{\alpha-\beta}\exp(t)}\over{t^{\alpha}-z}}{\rm
d}t, \Re(\alpha)>0, \Re(\beta)>0, z,\alpha,\beta\in C,\eqno(6.8)
$$which is an extension of (6.1) with the same path. As in the
previous case, the Mittag-Leffler function has the following
asymptotic estimates: $(iii)$:~~If $0<\alpha<2$ and $\mu$ is a real
number such that

$${{\pi\alpha}\over{2}}<\mu<\min[\pi,\pi\alpha],\eqno(6.9)
$$then there holds the following asymptotic expansion:

$$E_{\alpha,\beta}(z)={{1}\over{\alpha}}z^{(1-\beta)/\alpha}\exp(z^{{1}\over{\alpha}})-\sum_{r=1}^{N^{*}}{{1}\over{\Gamma(\beta-\alpha
r)}}{{1}\over{z^r}}+O\left[{{1}\over{z^{N^{*}+1}}}\right]\eqno(6.10)
$$as $|z|\rightarrow\infty, |\arg z|\le \mu$; and

$$E_{\alpha,\beta}(z)=-\sum_{r=1}^{N^{*}}{{1}\over{\Gamma(\beta-\alpha
r)}}{{1}\over{z^r}}+O\left[{{1}\over{z^{N^{*}+1}}}\right],\eqno(6.11)
$$as $|z|\rightarrow\infty, \mu\le |\arg z|\le \pi$.  $(iv)$:~~When
$\alpha\ge 2$ then there holds the following asymptotic expansion:

$$E_{\alpha,\beta}(z)={{1}\over{\alpha}}\sum_nz^{{1}\over{n}}\exp[\exp({{2n\pi
i}\over{\alpha}})z^{{1}\over{\alpha}}]^{1-\beta}-\sum_{r=1}^{N^{*}}{{1}\over{\Gamma(\beta-\alpha
r)}}{{1}\over{z^r}}+O\left[{{1}\over{z^{N^{*}+1}}}\right],\eqno(6.12)
$$as $|z|\rightarrow\infty, |\arg z|\le{{\alpha \pi}\over{2}}$ and
where the first sum is taken over all integers $n$ such that

$$|\arg (z)+2\pi n|\le{{\alpha \pi}\over{2}}.\eqno(6.13)
$$

\vskip.3cm \noindent {\bf 7.\hskip.3cm Integral representations}

\vskip.3cm In this section several integrals associated with
Mittag-Leffler functions are presented, which can be easily
established by the application by means of beta and gamma function
formulas and other techniques, see Erd\'elyi, et al. (1955),
Gorenflo et al. (1997, 2002).

$$\eqalignno{ \int_0^{\infty}{\rm
e}^{-\zeta}E_{\alpha}(\zeta^{\alpha}z){\rm d}\zeta
&={{1}\over{1-z}},|z|<1, \alpha\in C, \Re(\alpha)>0&(7.1)\cr
\int_0^{\infty}{\rm
e}^{-x}x^{\beta-1}E_{\alpha,\beta}(x^{\alpha}z){\rm
d}x&={{1}\over{1-z}}, |z|<1,\alpha,\beta\in C,
\Re(\alpha)>0,\Re(\beta)>0&(7.2)\cr
\int_0^{x}(x-\zeta)^{\beta-1}E_{\alpha}(\zeta^{\alpha}){\rm
d}\zeta&=\Gamma(\beta)x^{\beta}E_{\alpha,\beta+1}(x^{\alpha}),
\Re(\beta)>0&(7.3)\cr \int_0^{\infty}{\rm
e}^{-s\zeta}E_{\alpha}(-\zeta^{\alpha}){\rm
d}\zeta&={{s^{\alpha-1}}\over{1+s^{\alpha}}}, \Re(s)>0&(7.4)\cr
\int_0^{\infty}{\rm
e}^{-s\zeta}\zeta^{m\alpha+\beta-1}E_{\alpha,\beta}^{(m)}(\pm
a\zeta^{\alpha}){\rm
d}\zeta&={{m!s^{\alpha-\beta}}\over{(s^{\alpha}\mp a)^{m+1}}},
\Re(s)>0,\Re(\alpha)>0,\Re(\beta)>0,&(7.5)\cr}
$$where $\alpha,\beta\in C$ and

$$\eqalignno{E_{\alpha,\beta}^{(m)}(z)&={{{\rm d}^m}\over{{\rm
d}z^m}}E_{\alpha,\beta}(z)\cr
E_{\alpha}(-x^{\alpha})&={{2}\over{\pi}}\sin(\alpha\pi/2)\int_0^{\infty}{{\zeta^{\alpha-1}\cos(x\zeta)}\over{1+2\zeta^{\alpha}\cos(\alpha\pi/2)
+\zeta^{2\alpha}}}{\rm d}\zeta,\alpha\in C,\Re(\alpha)>0&(7.6)\cr
E_{\alpha}(-x)&={{1}\over{\pi}}\sin(\alpha\pi)\int_0^{\infty}{{\zeta^{\alpha-1}}\over{1+2\zeta^{\alpha}\cos(\alpha\pi)
+\zeta^{2\alpha}}}{\rm e}^{-\zeta x^{{1}\over{\alpha}}}\zeta~{\rm
d}\zeta,\alpha\in C, \Re(\alpha)>0&(7.7)\cr
E_{\alpha}(-x)&=1-{{1}\over{2\alpha}}+{{x^{{1}\over{\alpha}}}\over{\pi}}\int_0^{\infty}\arctan
\left[{{\zeta^{\alpha}+\cos(\alpha\pi)}\over{\sin(\alpha\pi)}}\right]{\rm
e}^{-\zeta x^{{1}\over{\alpha}}}\zeta~{\rm d}\zeta,\alpha\in
C,\Re(\alpha)>0.&(7.8)\cr}
$$

\vskip.3cm \noindent {\bf Note 7.1.}\hskip.3cm Equation (7.7) can be
employed to compute the numerical coefficients of the leading term
of the asymptotic expansion of $E_{\alpha}(-x)$. Equation (7.8)
yields

$$b_0(\alpha)={{\alpha}\over{\pi}}\Gamma(\alpha)\sin(2\alpha\pi)\int_0^{\infty}{{\zeta^{\alpha-1}}\over{1+2\zeta^{2\alpha}\cos(2\alpha\pi)
+\zeta^{4\alpha}}}{\rm d}\zeta,\alpha<{1\over2}.\eqno(7.9)
$$From Berberan-Santos (2005) and Gorenflo et al.(1997) the
following results hold:

$$\eqalignno{E_{\alpha}(-x)&={{2x}\over{\pi}}\int_0^{\infty}{{E_{2\alpha}(-t^2)}\over{x^2+t^2}}{\rm d}t,
0\le \alpha\le 1, \alpha\in C.&(7.10)\cr \noalign{\hbox{In
particular, the following cases are of importance}}
E_1(-x)&={{2x}\over{\pi}}\int_0^{\infty}{{\cosh(it)}\over{x^2+t^2}}{\rm
d}t=\exp(-x)&(7.11)\cr
E_{1\over2}(-x)&={{2x}\over{\pi}}\int_0^{\infty}{{\exp(-t^2)}\over{x^2+t^2}}{\rm
d}t={\rm e}^{x^2}{\rm erfc}(x).&(7.12)\cr
E_{1\over4}(-x)&={{2x}\over{\pi}}\int_0^{\infty}{{{\rm e}^{t^2}{\rm
erfc}(-t^2)}\over{x^2+t^2}}{\rm d}t.&(7.13)\cr}
$$

\vskip.3cm \noindent {\bf Note 7.2}.\hskip.3cm Some new properties
of the Mittag-Leffler functions are recently obtained by Gupta and
Debnath (2007).

\vskip.5cm \noindent {\bf 8.\hskip.3cm The H-function and its
special cases}

\vskip.3cm The H-function is defined by means of a Mellin-Barnes
type integral in the following manner (Mathai and Saxena, 1978):

$$\eqalignno{H_{p,q}^{m,n}(z)&=H_{p,q}^{m,n}\left[z\bigg\vert_{(b_q,B_q)}^{(a_p,A_p)}\right]
=H_{p,q}^{m,n}\left[z\bigg\vert_{(b_1,B_1),...,(b_q,B_q)}^{(a_1,A_1),...,(a_p,A_p)}\right]\cr
&={{1}\over{2\pi i}}\int_{\Omega}\Theta(\zeta)z^{-\zeta}{\rm
d}\zeta&(8.1)\cr \noalign{\hbox{where $i=\sqrt{(-1)}$ and }}
\Theta(s)&={{\left\{\prod_{j=1}^m\Gamma(b_j+sB_j)\right\}\left\{\prod_{j=1}^n\Gamma(1-a_j-sA_j)\right\}}
\over{\left\{\prod_{j=m+1}^q\Gamma(1-b_j-sB_j)\right\}\left\{\prod_{j=n+1}^p\Gamma(a_j+sA_j)\right\}}},&(8.2)\cr}
$$and an empty product is interpreted as unity; $m,n,p,q\in N_0$
with $0\le n\le p, 1\le m\le q, A_i,B_j\in R_{+}, a_i,b_j\in C,
i=1,...,p; j=1,...,q$ such that

$$A_i(b_j+k)\ne B_j(a_i-\lambda-1), k,\lambda\in N_0; i=1,...,n;
j=1,...,m,\eqno(8.3)
$$where we employ the usual notations: $N_0=0,1,2,...;
R=(-\infty,\infty), R_{+}=(0,\infty),$ and $C$ being the complex
number field. The contour $\Omega$ is the infinite contour which
separates all the poles of $\Gamma(b_j+sB_j), j=1,...,m$ from all
the poles of $\Gamma(1-a_i+sA_i), i=1,...,n$. The contour $\Omega$
could be $\Omega=L_{-\infty}$ or $\Omega=L_{+\infty}$ or
$\Omega=L_{i\gamma\infty}$ where $L_{-\infty}$ is a loop starting at
$-\infty$ encircling all the poles of $\Gamma(b_j+sB_j),j=1,...,m$
and ending at $-\infty$. $L_{+\infty}$ is a loop starting at
$+\infty$, encircling all the poles of
$\Gamma(1-a_i-sA_i),i=1,...,n$ and ending at $+\infty$.
$L_{i\gamma\infty}$ is the infinite semicircle starting at
$\gamma-i\infty$ and going to $\gamma+i\infty$. A detailed and
comprehensive account of the H-function is available from the
monographs of Mathai and Saxena (1978), Prudnikov et al. (1990) and
Kilbas and Saigo (2004). The relation connecting the  Wright's
function ${_p\psi_q}(z)$ and the H-function is given for the first
time in the monograph of Mathai and Saxena (1978, p.11, eq(1.7.8))
as

$${_p\psi_q}\left[\matrix{(a_1,A_1),...,(a_p,A_p)\cr
(b_1,B_1),...,(b_q,B_q)\cr}\vert
z\right]=H_{p,q+1}^{1,p}\left[-z\bigg\vert_{(0,1),(1-b_1,B_1),...,(1-b_q,B_q)}^{(1-a_1.A_1),...,(1-a_p,A_p)}\right],\eqno(8.4)
$$where ${_p\psi_q}(z)$ is the Wright's generalized hypergeometric
function (Wright, 1935, 1940); also see Erd\'elyi, et al. (1953,
Section 4.1), defined by means of the series representation in the
form

$${_p\psi_q}(z)=\sum_{r=0}^{\infty}{{\left\{\prod_{j=1}^p\Gamma(a_j+A_jr)\right\}}\over{\left\{\prod_{j=1}^q\Gamma(b_j+B_jr)\right\}}}{{z^r}\over{r!}}\eqno(8.5)
$$where $z\in C, a_i,b_j\in C, A_i, B_j\in R_{+},
A_i\ne 0, B_j\ne 0; i=1,...,p; j=1,...,q$,

$$\sum_{j=1}^qB_j-\sum_{i=1}^pA_i=\Delta>-1.
$$The Mellin-Barnes contour integral for the generalized Wright
function is given by

$${_p\psi_q}\left[\matrix{(a_p,A_p)\cr
(b_q,B_q)\cr}\vert z\right]={{1}\over{2\pi
i}}\int_{\Omega}{{\Gamma(s)\prod_{j=1}^p\Gamma(a_j-A_js)}\over{\prod_{j=1}^q\Gamma(b_j-sB_j)}}(-z)^{-s}{\rm
d}s\eqno(8.6)
$$where the path of integration separates all the poles of
$\Gamma(s)$ at the points $s=-\nu, \nu\in N_0$ lying to the left and
all the poles of $\prod_{j=1}^p\Gamma(a_j-sA_j),j=1,...,p$ at the
points $s=(A_j+\nu_j)/A_j, \nu_j\in N_0, j=1,...,p$ lying to the
right. If $\Omega=(\gamma-i\infty,\gamma+i\infty)$ then the above
representation is valid if either of the conditions are satisfied:

$$\Delta<1, |\arg (-z)|<{{(1-\Delta)\pi}\over{2}}, z\ne
0,\leqno(i)$$
$$\Delta =1,(1+\Delta)\gamma+{1\over2}<\Re(\delta), \arg
(-z)=0,z\ne 0,
\delta=\sum_{j=1}^qb_j-\sum_{j=1}^pa_j+{{p-q)}\over{2}}.\leqno(ii)
$$This result was proved by Kilbas, Saigo and Trujillo (2002). \vskip.2cm The
generalized Wright function includes many special functions besides
the Mittag-Leffler functions defined by the equations (1.1) and
(1.2). It is interesting to observe that for $A_i=B_j=1, i=1,...,p;
j=1,...,q$, (8.5) reduces to a generalized hypergeometric function
${_pF_q}(z)$. Thus

$${_p\psi_q}\left[\matrix{(a_p,1)\cr
(b_q,1)\cr}\vert
z\right]={{\prod_{j=1}^p\Gamma(a_j)}\over{\prod_{j=1}^q\Gamma(b_j)}}{_pF_q}(a_1,...,a_p;b_1,...,b_q;z)\eqno(8.7)
$$where $a_j\ne -v,j=1,...,p,v=0,1,...$; $b_j\ne -\lambda,
j=1,...,q, \lambda=0,1,...$; $p\le q$ or $p=q+1,|z|<1$. Wright
(1933) introduced a special case of (8.5) in the form

$$\phi(a,b;z)={_0\psi_1}\left[\matrix{~~\cr
(b,a)\cr}\vert
z\right]=\sum_{r=0}^{\infty}{{1}\over{\Gamma(ar+b)}}{{z^r}\over{r!}},\eqno(8.8)
$$which widely occurs in problems of fractional diffusion. It has
been shown by Saxena, Mathai and Haubold (2004a), also see Kiryakova
(1994), that

$$E_{\alpha,\beta}(z)={_1\psi_1}\left[\matrix{(1,1)\cr
(\beta,\alpha)\cr}\vert z\right]
=H_{1,2}^{1,1}\left[-z\bigg\vert_{(0,1),(1-\beta,\alpha)}^{(0,1)}\right].\eqno(8.9)
$$If we further take $\beta=1$ in (8.9) we find that

$$E_{\alpha,1}(z)=E_{\alpha}(z)={_1\psi_1}\left[\matrix{(1,1)\cr
(1,\alpha)\cr}\vert
z\right]=H_{1,2}^{1,1}\left[-z\bigg\vert_{(0,1),(0,\alpha)}^{(0,1)}\right]\eqno(8.10)
$$where $\alpha\in C, \Re(\alpha)>0$.

\vskip.3cm \noindent{\bf Remark 8.1.}\hskip.3cm A series of papers
are devoted to the application of the Wright function in partial
differential equation of fractional order extending the classical
diffusion and wave equations. Mainardi (1997) has obtained the
result for a fractional diffusion wave equation in terms of the
fractional Green function involving the Wright function. The
scale-variant solutions of some partial differential equations of
fractional order were obtained in terms of special cases of the
generalized Wright function by Buckwar and Luchko (1998) and Luchko
and Gorenflo (1998).

\vskip.5cm \noindent {\bf 9.\hskip.3cm Mellin-Barnes integrals for
Mittag-Leffler functions}

\vskip.3cm These integrals can be obtained from the identities (8.9)
and (8.10).

\vskip.3cm \noindent {\bf Lemma 9.1.}\hskip.3cm{\it If
$\Re(\alpha)>0, \Re(\beta)>0$ and $z\in C$ the following
representations are obtained:

$$\eqalignno{E_{\alpha}(z)&={{1}\over{2\pi
i}}\int_{\gamma-i\infty}^{\gamma+i\infty}{{\Gamma(s)\Gamma(1-s)}\over{\Gamma(1-\alpha
s)}}(-z)^{-s}{\rm d}s&(9.1)\cr E_{\alpha,\beta}(z)&={{1}\over{2\pi
i}}\int_{\gamma-i\infty}^{\gamma+i\infty}{{\Gamma(s)\Gamma(1-s)}\over{\Gamma(\beta-\alpha
s)}}(-z)^{-s}{\rm d}s&(9.2)\cr}
$$where the path of integration separates all the poles of
$\Gamma(s)$ at the points $s=-\nu, \nu=0,1,...$ from those of
$\Gamma(1-s)$ at the points $s=1+v,v=0,1,...$.}

\vskip.2cm On evaluating the residues at the poles of the gamma
function $\Gamma(1-s)$ we obtain the following analytic continuation
formulas for the Mittag-Leffler functions:

$$\eqalignno{E_{\alpha}(z)&={{1}\over{2\pi
i}}\int_{\gamma-i\infty}^{\gamma+i\infty}{{\Gamma(s)\Gamma(1-s)}\over{\Gamma(1-\alpha
s)}}(-z)^{-s}{\rm
d}s=-\sum_{k=1}^{\infty}{{z^{-k}}\over{\Gamma(1-\alpha k)}}&(9.3)\cr
\noalign{\hbox{and}} E_{\alpha,\beta}(z)&={{1}\over{2\pi
i}}\int_{\gamma-i\infty}^{\gamma+i\infty}{{\Gamma(s)\Gamma(1-s)}\over{\Gamma(\beta-\alpha
s)}}(-z)^{-s}{\rm
d}s=-\sum_{k=1}^{\infty}{{z^{-k}}\over{\Gamma(\beta-\alpha
k)}}.&(9.4)\cr}
$$

\vskip.3cm \noindent {\bf 10.\hskip.3cm Relation with
Riemann-Liouville fractional calculus operators}

\vskip.3cm In this section, we present the relations of
Mittag-Leffler functions with the  left and right-sided operators of
Riemann-Liouville fractional calculus, which are defined below.

$$\eqalignno{(I_{0_{+}}^{\alpha}\phi)(x)&={{1}\over{\Gamma(\alpha)}}\int_0^x(x-t)^{\alpha-1}\phi(t){\rm
d}t,\Re(\alpha)>0&(10.1)\cr
(I_{-}^{\alpha})(x)&={{1}\over{\Gamma(\alpha)}}\int_x^{\infty}(t-x)^{\alpha-1}\phi(t){\rm
d}t,\Re(\alpha)>0&(10.2)\cr
(D_{0_{+}}^{\alpha}\phi)(x)&=\left({{{\rm d}}\over{{\rm
d}x}}\right)^{[\alpha]+1}[I_{0_{+}}^{1-\{\alpha\}}\phi](x)\cr
&={{1}\over{\Gamma(1-\{\alpha\})}}\left({{{\rm d}}\over{{\rm
d}x}}\right)^{[\alpha]+1}\int_0^x(x-t)^{\alpha-1}\phi(t){\rm d}t,
\Re(\alpha)>0&(10.3)\cr (D_{-}^{\alpha}\phi)(x)&=\left(-{{{\rm
d}}\over{{\rm
d}x}}\right)^{[\alpha]+1}[I_{-}^{1-\{\alpha\}}\phi](x)\cr
&={{1}\over{\Gamma(1-\{\alpha\})}}\left(-{{{\rm d}}\over{{\rm
d}x}}\right)^{[\alpha]+1}\int_x^{\infty}(t-x)^{-\{\alpha\}}\phi(t){\rm
d}t, \Re(\alpha)>0&(10.4)\cr}
$$where $[\alpha]$ means the maximal integer not exceeding $\alpha$
and $\{\alpha\}$ is the fractional part of $\alpha$.

\vskip.3cm \noindent {\bf Note 10.1}.\hskip.3cm The fractional
integrals (10.1) and (10.2) are connected by the relation (Kilbas,
2005, p.118)

$$[I_{-}^{\alpha}\phi(1/t)](x)=x^{\alpha-1}\left(I_{0_{+}}^{\alpha}[t^{-\alpha-1}\phi(t)]\right)\left({{1}\over{x}}\right).$$

\vskip.3cm \noindent {\bf Theorem 10.1.}\hskip.3cm {\it Let
$\Re(\alpha)
>0$ and $\Re(\beta) >0$ then there holds the formulas
$$\eqalignno{(I_{0_{+}}^{\alpha}[t^{\alpha-1}E_{\alpha,\beta}(at^{\alpha})])(x)&={{x^{\beta-1}}\over{a}}\left[E_{\alpha,\beta}(ax^{\alpha})
-{{1}\over{\Gamma(\beta)}}\right], a\ne 0&(10.5)\cr
(I_{0_{+}}^{\alpha}[E_{\alpha}(at^{\alpha})])(x)&={{1}\over{a}}\left[E_{\alpha}(ax^{\alpha})-1\right],
a\ne 0&(10.6)\cr \noalign{\hbox{which by virtue of the definitions
(1.1) and (1.2) can be written as}}
(I_{0_{+}}^{\alpha}[t^{\beta-1}E_{\alpha,\beta}(at^{\alpha})])(x)&=x^{\alpha+\beta-1}E_{\alpha,\alpha+\beta}(ax^{\alpha})&(10.7)\cr
(I_{0_{+}}^{\alpha}[E_{\alpha}(at^{\alpha})])(x)&=x^{\alpha}E_{\alpha,\alpha+1}(ax^{\alpha}).&(10.8)\cr}$$}

\vskip.3cm \noindent {\bf Theorem 10.2.}\hskip.3cm {\it Let
$\Re(\alpha)>0$ and $\Re(\beta)>0$ then there holds the formulas

$$\eqalignno{(I_{-}^{\alpha}[t^{-\alpha-\beta}E_{\alpha,\beta}(at^{-\alpha})])(x)&={{x^{\alpha-\beta}}\over{a}}\left[E_{\alpha,\beta}(ax^{-\alpha})-{{1}\over{\Gamma(\beta)}}\right],
a\ne 0&(10.9)\cr
(I_{-}^{\alpha}[t^{-\alpha-1}E_{\alpha}(at^{-\alpha})])(x)&={{1}\over{a}}x^{\alpha-1}[E_{\alpha}(ax^{-\alpha})-1],a\ne
0.&(10.10)\cr}$$}

\vskip.3cm \noindent {\bf Theorem 10.3}.\hskip.3cm{\it Let
$0<\Re(\alpha)<1$ and $\Re(\beta)>\Re(\alpha)$ then there holds the
formulas

$$\eqalignno{(D_{0_{+}}^{\alpha}[t^{\beta-1}E_{\alpha,\beta}(at^{\alpha})])(x)&={{x^{\beta-\alpha-1}}\over{\Gamma(\beta-\alpha)}}+ax^{\beta-1}E_{\alpha,\beta}(ax^{\alpha})&(10.11)\cr
(D_{0_{+}}^{\alpha}[E_{\alpha}(at^{\alpha})])(x)&={{x^{-\alpha}}\over{\Gamma(1-\alpha)}}+aE_{\alpha}(ax^{\alpha}).&(10.12)\cr}$$}

\vskip.3cm \noindent {\bf Theorem 10.4.}\hskip.3cm{\it Let
$\Re(\alpha)>0$ and $\Re(\beta)>[\Re(\alpha)]+1$ then there holds
the formula

$$(D_{-}^{\alpha}[t^{\alpha-\beta}E_{\alpha,\beta}(at^{-\alpha})])(x)={{x^{-\beta}}\over{\Gamma(\beta-\alpha)}}
+ax^{-\alpha-\beta}E_{\alpha,\beta}(ax^{-\alpha}).\eqno(10.13)$$}

\vskip.3cm \noindent {\bf 11.\hskip.3cm Generalized Mittag-Leffler
type functions}

\vskip.3cm By means of the series representation a generalization of
(1.1) and (1.2) is introduced by Prabhakar (1971) as

$$E_{\alpha,\beta}^{\gamma}(z)=\sum_{n=0}^{\infty}{{(\gamma)_n}\over{n!\Gamma(\beta+\alpha
n)}},\alpha,\beta,\gamma\in C,
\Re(\alpha)>0,\Re(\beta)>0,\eqno(11.1)
$$where
$$(\gamma)_n=\gamma(\gamma+1)...(\gamma+n-1)={{\Gamma(\gamma+n)}\over{\Gamma(\gamma)}}
$$whenever $\Gamma(\gamma)$ is defined, $(\gamma)_0=1,\gamma\ne 0$. It is an entire
function of order $\rho=[\Re(\alpha)]^{-1}$ and type
$\sigma={{1}\over{\rho}}\left[\{\Re(\alpha)\}^{\Re(\alpha)}\right]^{-\rho}$.
It is a special case of Wright's generalized hypergeometric
function, Wright (1934,1935) as well as the H-function (Mathai and
Saxena, 1978). For various properties of this function with
applications, see Prabhakar (1971). Some special cases of this
function are enumerated below.

$$\eqalignno{E_{\alpha}(z)&=E_{\alpha,1}^1(z).&(11.2)\cr
E_{\alpha,\beta}(z)&=E_{\alpha,\beta}^1(z).&(11.3)\cr \alpha\gamma
E_{\alpha,\beta}^{\gamma+1}(z)&=(1+\alpha\gamma-\beta)E_{\alpha,\beta}^{\gamma}(z)+E_{\alpha,\beta-1}^{\gamma}(z).&(11.4)\cr
\phi(\alpha,\beta;z)&=\Gamma(\beta)E_{1,\beta}^{\alpha}(z)&(11.5)\cr}
$$where $\phi(\alpha,\beta;z)$ is the Kummer's confluent
hypergeometric function. $E_{\alpha,\beta}^{\gamma}(z)$ has the
following representations in terms of the Wright's function and
H-function.

$$\eqalignno{E_{\alpha,\beta}^{\gamma}(z)&={{1}\over{\Gamma(\gamma)}}{_1\psi_1}\left[\matrix{(\gamma,1)\cr
(\beta,\alpha)\cr}; z\right]&(11.6)\cr
&={{1}\over{\Gamma(\gamma)}}H_{1,2}^{1,1}\left[-z\bigg\vert_{(0,1),(1-\beta,\alpha)}^{(1-\gamma,1)}\right]&
(11.7)\cr &={{1}\over{2\pi\omega
\Gamma(\gamma)}}\int_{c-i\infty}^{c+i\infty}{{\Gamma(s)\Gamma(\gamma-s)}\over{\Gamma(\beta-\alpha
s)}}(-z)^{-s}{\rm d}s, \Re(\gamma)>0&(11.8)\cr}
$$where ${_1\psi_1}(\cdot)$ and $H_{1,2}^{1,1}(\cdot)$ are
respectively Wright generalized hypergeometric function and the
H-function. In the Mellin-Barnes integral representation,
$\omega=\sqrt{-1}$ and  the $c$ in the contour is such that
$0<c<\Re(\gamma)$ and it is assumed that the poles of $\Gamma(s)$
and $\Gamma(\gamma-s)$ are separated by the contour. The following
two theorems are given by Kilbas, Saigo and Saxena (2004).

\vskip.3cm \noindent {\bf Theorem 11.1.}\hskip.3cm{\it If
$\alpha,\beta,\gamma,a\in C, \Re(\alpha)>0, \Re(\beta)>n,
\Re(\gamma)>0$ then for $n\in N$ the following results hold:

$$\eqalignno{{{{\rm d}^n}\over{{\rm
d}z^n}}[z^{\beta-1}E_{\alpha,\beta}^{\gamma}(az^{\alpha})]&=z^{\beta-n-1}E_{\alpha,\beta-n}^{\gamma}(az^{\alpha}).&(11.9)\cr
\noalign{\hbox{In particular,}} {{{\rm d}^n}\over{{\rm
d}z^n}}[z^{\beta-1}E_{\alpha,\beta}(az^{\alpha})]&=z^{\beta-n-1}E_{\alpha,\beta-n}(az^{\alpha})&(11.10)\cr
\noalign{\hbox{and}} {{{\rm d}^n}\over{{\rm
d}z^n}}[z^{\beta-1}\phi(\gamma,\beta;az)]&={{\Gamma(\beta)}\over{\Gamma(\beta-n)}}z^{\beta-n-1}
\phi(\gamma;\beta-n;az).&(11.11)\cr}
$$}

\vskip.3cm\noindent {\bf Theorem 11.2.}\hskip.3cm{\it If
$\alpha,\beta,\gamma,a,\nu,\sigma\in C, \Re(\alpha)>0, \Re(\beta)>0,
\Re(\gamma)>0, \Re(\nu)>0, \Re(\sigma)>0$ then

$$\int_0^x(x-t)^{\beta-1}E_{\alpha,\beta}^{\gamma}[a(x-t)^{\alpha}]t^{\nu-1}E_{\alpha,\nu}^{\sigma}(at^{\alpha}){\rm
d}t=x^{\beta+\nu-1}E_{\alpha,\beta+\nu}^{\gamma+\sigma}(ax^{\alpha}).\eqno(11.12)
$$}The proof of (11.12) can be developed with the help of the
Laplace transform formula

$$L\left[x^{\beta-1}E_{\alpha,\beta}^{\gamma}(ax^{\alpha})\right](s)=s^{-\beta}(1-as^{-\alpha})^{-\gamma}\eqno(11.13)
$$where $\alpha,\beta,\gamma,a\in C$, $\Re(\alpha)>0, \Re(\beta)>0,
\Re(\gamma)>0$, $s\in C, \Re(s)>0, |as^{-\alpha}|<1$. For
$\gamma=1$, (11.13) reduces to

$$L\left[x^{\beta-1}E_{\alpha,\beta}(ax^{\alpha})\right](s)=s^{-\beta}(1-as^{-\alpha})^{-1}.\eqno(11.14)
$$Generalization of the above two results is given by Saxena (2002).

$$L\left[t^{\rho-1}E_{\beta,\gamma}^{\delta}(at^{\alpha})\right](s)={{s^{-\rho}}\over{\Gamma(\delta)}}{_2\psi_1}\left[\matrix{(\delta,1),
(\rho,\alpha)\cr
(\gamma,\beta)\cr};{{a}\over{s^{\alpha}}}\right],\eqno(11.15)
$$where $\Re(\beta)>0, \Re(\gamma)>0, \Re(s)>0, \Re(\rho)>0,
s>|a|^{{1}\over{\Re(\alpha)}}, \Re(\delta)>0$. \vskip.2cm Relations
connecting the function defined by (11.1) and the Riemann-Liouville
fractional integrals and derivatives are given by Saxena and Saigo
(2005) in the form of nine theorems. Some of the interesting
theorems are given below.

\vskip.3cm \noindent {\bf Theorem 11.3}.\hskip.3cm{\it Let
$\alpha>0,\beta>0,\gamma>0$ and $a\in R$. Let $I_{0_{+}}^{\alpha}$
be the left-sided operator of Riemann-Liouville fractional integral.
Then there holds the formula

$$(I_{0_{+}}^{\alpha}[t^{\gamma-1}E_{\beta,\gamma}^{\delta}(at^{\beta})])(x)=x^{\alpha+\gamma-1}E_{\beta,\alpha+\gamma}^{\delta}(ax^{\beta}).\eqno(11.16)
$$}

\vskip.3cm \noindent {\bf Theorem 11.4.}\hskip.3cm {\it Let
$\alpha>0,\beta>0,\gamma>0$ and $a\in R$. Let $I_{-}^{\alpha}$ be
the right-sided operator of Riemann-Liouville fractional integral.
Then there holds the formula

$$(I_{-}^{\alpha}[t^{-\alpha-\gamma}E_{\beta,\gamma}^{\delta}(at^{-\beta})])(x)=x^{-\gamma}E_{\beta,\alpha+\gamma}^{\delta}(ax^{-\beta}).\eqno(11.17)
$$}

\vskip.3cm \noindent {\bf Theorem 11.5}.\hskip.3cm{\it Let
$\alpha>0,\beta>0,\gamma>0$ and $a\in R$. Let $D_{0_{+}}^{\alpha}$
be the left-sided operator of Riemann-Liouville fractional
derivative. Then there holds the formula

$$(D_{0_{+}}^{\alpha}[t^{\gamma-1}E_{\beta,\gamma}^{\delta}(at^{\beta})])(x)=x^{\gamma-\alpha-1}E_{\beta,\gamma-\alpha}^{\delta}(ax^{\beta}).\eqno(11.18)
$$}

\vskip.3cm \noindent {\bf Theorem 11.6.}\hskip.3cm{\it Let
$\alpha>0, \beta>0, \gamma-\alpha+\{\alpha\}>1$ and $a\in R$. Let
$D_{-}^{\alpha}$ be the right-sided operator of Riemann-Liouville
fractional derivative. Then there holds the formula

$$(D_{-}^{\alpha}[t^{\alpha-\gamma}E_{\beta,\gamma}^{\delta}(at^{-\beta})])(x)=x^{-\gamma}E_{\beta,\gamma-\alpha}^{\delta}(ax^{-\beta}).\eqno(11.19)
$$}In a series of papers by Luchko and Yakubovich (1990, 1994),
Luckho and Srivastava (1995), Al-Bassam and Luchko (1995), Hadid and
Luchko (1996), Gorenflo and Luchko (1997), Gorenflo et al. (2000),
Luchko and Gorenflo (1999), the operational method was developed to
solve in closed forms certain classes of differential equations of
fractional order and also integral equations. Solutions of the
equations and problems considered are obtained in terms of
generalized Mittag-Leffler functions. The exact solution of certain
differential equation of fractional order is given by Luchko and
Srivastava (1995) in terms of the function (11.1) by using
operational method. In other papers, the solutions are established
in terms of the following functions of Mittag-Leffler type: If
$z,\rho,\beta_j\in C, \Re(\alpha_j)>0, j=1,...,m$ and $m\in N$ then
$$E_{\rho}((\alpha_j,\beta_j)_{1,m};(z))=\sum_{k=0}^{\infty}{{(\rho)_k}\over{\prod_{j=1}^m\Gamma(a_jk+\beta_j)}}{{z^k}\over{k!}}.\eqno(11.20)
$$For $m=1$, (11.20) reduces to (11.1). The Mellin-Barnes integral
for this function is given by

$$\eqalignno{E_{\rho}((\alpha_j,\beta_j)_{1,m};(z))&={{1}\over{2\pi
 i\Gamma(\rho)}}\int_{\gamma-i\infty}^{\gamma+i\infty}{{\Gamma(s)\Gamma(\rho-s)}\over{\prod_{j=1}^m\Gamma(\beta_j-\alpha_js)}}(-z)^{-s}{\rm
d}s, i=\sqrt{-1}&(11.21)\cr
&={{1}\over{\Gamma(\rho)}}H_{1,m+1}^{1,1}\left[-z\bigg\vert_{(0,1),(1-\beta_j,\alpha_j),j=1,...,m}^{(1-\rho,1)}\right]&(11.22)\cr}
$$where $0,\gamma<\Re(\rho), \Re(\rho)>0$ and the contour separates
the poles of $\Gamma(s)$ from those of $\Gamma(\rho-s)$.
$\Re(\alpha_j)>0, j=1,...,m, \arg(-z)<\pi$. The Laplace transform of
the function defined by (11.20) is given by

$$L[E_{\rho}((\alpha_j,\beta_j)_{1,m};-t)](s)={{1}\over{s\Gamma(\rho)}}{_2\psi_m}\left[\matrix{(\rho,1),
(1,1)\cr (\alpha_j,\beta_j)_{1,m}\cr}\bigg\vert
{{1}\over{s}}\right],\eqno(11.23)
$$where $\Re(s)>0$.

\vskip.3cm \noindent {\bf Remark 11.1.}\hskip.3cm In a recent paper,
Kilbas, Saigo and Saxena (2002) obtained a closed form solution of a
fractional generalization of a free electron equation of the form:

$$D_{\tau}^{\alpha}a(\tau)=\lambda\int_0^{\tau}t^{\delta}a(\tau-t)E_{\rho,\delta+1}^b(i\nu
t^{\rho}){\rm d}t+\beta\tau^{\sigma}E_{\rho,\sigma+1}^{\gamma}(i\nu
t^{\rho}), 0\le \tau\le 1, i=\sqrt{-1}\eqno(11.24)
$$where $b,\lambda\in C, \nu,\beta\in R_{+}, \alpha>0,\rho>0;
\alpha>-1,\rho>-1,\delta>-1$ and $E_{\rho,\delta+1}^b(\cdot)$ is the
generalized Mittag-Leffler function given by (11.1), and
$\alpha(\tau)$ is the unknown function to be determined.

\vskip.3cm \noindent {\bf Remark 11.2.}\hskip.3cm The solution of
fractional differential equations by the operational methods are
also obtained in terms of certain multivariate Mittag-Leffler
functions defined below: The multivariate Mittag-Leffler function of
$n$ complex variables $z_1,...,z_n$ with complex parameters
$a_1,...,a_n,b\in C$ is defined as

$$E_{(a_1,...,a_n),b}(z_1,...,z_n)=\sum_{k=0}^{\infty}\sum_{L_1,...,L_n\ge
0}^{L_1+...+L_n=k}{{k}\choose{L_1,...,L_n}}{{\prod_{j=1}^nz_j^{L_j}}\over{\Gamma(b\pm\sum_{j=1}^na_jL_j)}}\eqno(11.25)
$$in terms of the multinomial coefficients

$${{k}\choose{L_1,...,L_n}}={{k!}\over{L_1!...L_n!}}, k,L_j,
\in N_0, j=1,...,m.\eqno(11.26)
$$
\vskip.1cm Another generalization of the Mittag-Leffler function
(1.2) was introduced by Kilbas and Saigo (1995) in terms of a
special function of the form

$$E_{\alpha,m,\beta}(z)=\sum_{k=0}^{\infty}c_kz^k,~c_0=1,c_k=\prod_{i=0}^{k-1}{{\Gamma(\alpha(im+\beta)+1)}\over{\Gamma(\alpha(im+\beta+1)+1)}},
k\in N_0=\{0,1,2,...\}\eqno(11.27)
$$where an empty product is to be interpreted as unity; $\alpha,\beta\in
C$ are complex numbers and $m\in R, \Re(\alpha)>0, m>0,
\alpha(im+\beta)\notin Z^{-}=\{0,-1,-2,...\}, i=0,1,2,...$ and for
$m=1$ the above defined function reduces to a constant multiple of
the Mittag-Leffler function, namely

$$E_{\alpha,1,\beta}(z)=\Gamma(\alpha \beta+1)E_{\alpha,\alpha
\beta+1}(z),\eqno(11.28)
$$where $\Re(\alpha)>0$ and $\alpha (i+\beta)\notin Z^{-}$. It is an
entire function of $z$ of order $[\Re(\alpha)]^{-1}$ and type
$\sigma=1/m$, see Gorenflo et al. (1998). Certain properties of this
function associated with Riemann-Liouville fractional integrals and
derivatives are obtained and exact solutions of certain integral
equations of Abel-Volterra type are derived by their applications
(Kilbas and Saigo, 1995, 1996). Its recurrence relations, connection
with hypergeometric functions and differential formulas are obtained
by Gorenflo, Kilbas and Rogosin (1998), also see, Gorenflo and
Mainardi (1996). In order to present the applications of
Mittag-Leffler functions we give definitions of Laplace transform,
Fourier transform, Riemann-Liouville fractional calculus operators,
Caputo operator and Weyl fractional operators in the next section.

\vskip.3cm \noindent {\bf 12.\hskip.3cm Laplace and Fourier
transforms, fractional calculus operators}

\vskip.3cm We will need the definitions of the well-known Laplace
and Fourier transforms of a function $N(x,t)$ and their inverses,
which are useful in deriving the solution of fractional differential
equations governing certain physical problems. The Laplace transform
of a function $N(x,t)$ with respect to $t$ is defined by

$$L\{N(x,t)\}=\tilde{N}(x,s)=\int_0^{\infty}{\rm e}^{-st}N(x,t){\rm
d}t, t>0, x\in R\eqno(12.1)
$$where $\Re(s)>0$ and its inverse transform with respect to $s$ is
given by

$$L^{-1}\{\tilde{N}(x,s)\}=L^{-1}\{\tilde{N}(x,s);t\}={{1}\over{2\pi
i}}\int_{c-i\infty}^{c+i\infty}{\rm e}^{st}\tilde{N}(x,s){\rm
d}s.\eqno(12.2)
$$The Fourier transform of a function $N(x,t)$ with respect to $x$
is defined by

$$F\{N(x,t)\}=F^{*}(k,t)=\int_{-\infty}^{\infty}{\rm
e}^{ikx}N(x,t){\rm d}x.\eqno(12.3)
$$The inverse Fourier transform with respect to $k$ is given by the
formula

$$F^{-1}\{F^{*}(k,t)\}={{1}\over{2\pi}}\int_{-\infty}^{\infty}{\rm
e}^{-ikx}F^{*}(k,t){\rm d}k.\eqno(12.4)
$$From Mathai and Saxena (1978) and Prudnikov et al.(1990, p.355,
eq(2.25.3)) it follows that the Laplace transform of the H-function
is given by

$$L\{t^{\rho-1}H_{p,q}^{m,n}\left[zt^{\sigma}\bigg\vert_{(b_q,B_q)}^{(a_p,A_p)}\right]\}
=s^{-\rho}H_{p+1,q}^{m,n+1}\left[zs^{-\sigma}\bigg\vert_{(b_q,B_q)}^{(1-\rho,\sigma),(a_p,A_p)}\right],\eqno(12.5)
$$where $\Re(s)>0, \Re(\rho+\sigma\min_{1\le j\le
m}({{b_j}\over{B_j}})>0, \sigma>0$,
$$|\arg
z|<{1\over2}\pi\Omega,~\Omega>0,
\Omega=\sum_{i=1}^nA_i-\sum_{i=n+1}^pA_i+\sum_{j=1}^mB_j-\sum_{j=m+1}^qB_j.\eqno(12.6)
$$By virtue of the cancelation law for the H-function (Mathai and
Saxena, 1978) it can be readily seen that

$$L^{-1}\{s^{-\rho}H_{p,q}^{m,n}\left[zs^{\sigma}\bigg\vert_{(b_q,B_q)}^{(a_p,A_p)}\right]\}
=t^{\rho-1}H_{p+1,q}^{m,n}\left[zt^{-\sigma}\bigg\vert_{(b_q,B_q)}^{(a_p,A_p),
(\rho,\sigma)}\right],\eqno(12.7)
$$where $\sigma>0, \Re(s)>0, \Re[\rho+\sigma\max_{1\le i\le
n}({{1-a_i}\over{A_i}})]>0$, $|\arg z|<{1\over2}\pi
\Omega_1,~\Omega_1>0, \Omega_1=\Omega-\rho$. In view of the results

$$J_{-{1\over2}}(x)=\sqrt{{{2}\over{\pi x}}}\cos x\eqno(12.8)
$$the cosine transform of the H-function (Mathai and Saxena,
1978, p.49) is given by

$$\int_0^{\infty}t^{\rho-1}\cos(kt)H_{p,q}^{m,n}[at^{\mu}\big\vert_{(b_q,B_q)}^{(a_p,A_p)}]{\rm
d}t={{\pi}\over{k^{\rho}}}H_{q+1,p+2}^{m+1,n}\left[{{k^{\mu}}\over{a}}\bigg\vert_{(\rho,\mu),(a-a_p,A_p),
({{1+\rho}\over{2}},{{\mu}\over{2}})}^{(1-b_q,B_q),
({{1+\rho}\over{2}},{{\mu}\over{2}})}\right]\eqno(12.9)
$$where $\Re[\rho+\mu\min_{1\le j\le m}({{b_j}\over{B_j}})]>0$,
$\Re[\rho+\mu\max_{1\le j\le n}({{a_j-1}\over{A_j}})<0$, $|\arg
a|<{1\over2}\pi\Omega, \Omega>0$;
$\Omega=\sum_{j=1}^mB_j-\sum_{j=m+1}^qB_j+\sum_{j=1}^nA_j-\sum_{j=n+1}^pA_j$,
 and $k>0$.
 \vskip.2cm The definitions of fractional integrals used in the
 analysis are defined below. The Riemann-Liouville fractional
 integral of order $\nu$ is defined by (Miller and Ross, 1993, p.45)

 $${_0D}_t^{-\nu}f(x,t)={{1}\over{\Gamma(\nu)}}\int_0^t(t-u)^{\nu-1}f(x,u){\rm
 d}u,\eqno(12.10)
 $$where $\Re(\nu)>0$. Following Samko, Kilbas and Marichev (1993,
 p.37) we define the Riemann-Liouville fractional derivative for
 $\alpha >0$ in the form

 $${_0D}_t^{\alpha}f(x,t)={{1}\over{\Gamma(n-\alpha)}}{{{\rm
 d}^n}\over{{\rm d}t^n}}\int_0^t(t-u)^{n-\alpha-1}f(x,u){\rm d}u,~~
 n=[\alpha]+1\eqno(12.11)
 $$where $[\alpha]$ means the integral part of the number $\alpha$.
 From Erd\'elyi, et al. (1954b, p.182) we have

 $$L\{{_0D}_t^{-\nu}f(x,t)\}=s^{-\nu}F(x,s),\eqno(12.12)
 $$where $F(x,s)$ is the Laplace transform of $f(x,t)$ with respect
 to $t$, $\Re(s)>0, \Re(\nu)>0$. The Laplace transform of the
 fractional derivative defined by (12.11) is given by Oldham and
 Spanier (1974, p.134,eq (8.1.3))

 $$L\{{_0D}_t^{\alpha}f(x,t)\}=s^{\alpha}F(x,s)-\sum_{k=1}^ns^{k-1}{_0D}_t^{\alpha-k}f(x,t)|_{t=0},~
 n-1<\alpha \le n,\eqno(12.13)
 $$In certain boundary-value problems arising in the theory of
 visco-elasticity and in the hereditary solid mechanics the
 following fractional derivative of order $\alpha>0$ is introduced
 by Caputo (1969) in the form

 $$\eqalignno{D_t^{\alpha}f(x,t)&={{1}\over{\Gamma(m-\alpha)}}\int_0^t(t-\tau)^{m-\alpha-1}f^{(m)}(x,t){\rm
 d}t, m-1<\alpha\le m, \Re(\alpha)>0,m\in N&(12.14)\cr
 &={{\partial^m}\over{\partial t^m}}f(x,t),\hbox{  if
 }\alpha=m,&(12.15)\cr}
 $$where ${{\partial^m}\over{\partial t^m}}f$ is the $m$-th partial
 derivative of the function $f(x,t)$ with respect to $t$. The
 Laplace transform of this derivative is given by Podlubny (1999) in
 the form

 $$L\{D_t^{\alpha}f(t);s\}=s^{\alpha}F(s)-\sum_{k=0}^{m-1}s^{\alpha-k-1}f^{(k)}(0_{+}),~m-1<\alpha\le
 m.\eqno(12.16)
 $$The above formula is very useful in deriving the solution of
 differintegral equations of fractional order governing certain
 physical problems of reaction and diffusion. Making use of the
 definitions (12.10) and (12.11) it readily follows that for
 $f(t)=t^{\rho}$ we obtain

 $${_0D}_t^{-\nu}t^{\rho}={{\Gamma(\rho+1)}\over{\Gamma(\rho+\nu+1)}}t^{\rho+\nu},
 \Re(\nu)>0, \Re(\rho)>-1; t>0\eqno(12.17)
 $$and

 $${_0D}_t^{\nu}t^{\rho}={{\Gamma(\rho+1)}\over{\Gamma(\rho-\nu+1)}}t^{\rho-\nu},
 \Re(\nu)<0, \Re(\rho)>-1; t>0.\eqno(12.18)
 $$On taking $\rho=0$ in (12.18) we find that

 $${_0D}_t^{\nu}[1]={{1}\over{\Gamma(1-\nu)}}t^{-\nu},
 t>0, \Re(\nu)<1.\eqno(12.19)
 $$From the above result, we infer that the Riemann-Liouville
 derivative of unity is not zero. We also need the Weyl fractional
 operator defined by

 $${_{-\infty}D}_x^{\mu}={{1}\over{\Gamma(n-\mu)}}{{{\rm
 d}^n}\over{{\rm d}t^n}}\int_{-\infty}^t (t-u)^{n-\mu-1}f(u){\rm
 d}u,\eqno(12.20)
 $$where $n=[\mu]+1$ is the integer part of $\mu>0$. Its Fourier
 transform is (Metzler and Klafter, 2000, p.59, A.11)

 $$F\{{_{-\infty}D}_x^{\mu}f(x)\}=(ik)^{\mu}f^{*}(k),\eqno(12.21)
 $$where we define the Fourier transform as

 $$h^{*}(q)=\int_{-\infty}^{\infty}h(x)\exp(iqx){\rm
 d}x.\eqno(12.22)
 $$Following the convention initiated by Compte (1996) we suppress
 the imaginary unit in Fourier space by adopting a slightly
 modified form of the above result in our investigations (Metzler
 and Klafter, 2000,p.59, A.12).

 $$F\{{_{-\infty}D}_x^{\mu}f(x)\}=-|k|^{\mu}f^{*}(k),\eqno(12.23)
 $$instead of (12.21). \vskip.2cm We now proceed to discuss the
 various applications of Mittag-Leffler functions in applied
 sciences. In order to discuss the application of Mittag-Leffler
 function in kinetic equations, we derive the solution of two
 kinetic equations in the next section.

 \vskip.3cm
 \noindent
 {\bf 13.\hskip.3cm Application in kinetic equations}

 \vskip.3cm
 \noindent
 {\bf Theorem 13.1.}\hskip.3cm{\it If $\Re(\nu)>0$ then the solution
 of the integral equation

 $$N(t)-N_0=-c^{\nu}{_0D}_t^{-\nu}N(t)\eqno(13.1)
 $$is given by

 $$N(t)=N_0E_{\nu}(-c^{\nu}t^{\nu})\eqno(13.2)
 $$where $E_{\nu}(t)$ is the Mittag-Leffler function defined in
 (1.1).}

 \vskip.3cm
 \noindent
 {\bf Proof.}\hskip.3cm Applying Laplace transform to both sides of
 (13.1) and using (12.12) it gives

 $$\tilde{N}(s)=L\{N(t);s\}=N_0s^{-1}[1+(s/c)^{-\nu}]^{-1}.\eqno(13.3)
 $$By virtue of the relation

 $$\eqalignno{L^{-1}\{s^{-\rho}\}&={{t^{\rho-1}}\over{\Gamma(\rho)}},
 \Re(\rho)>0, \Re(s)>0, s\in C,&(13.4)\cr
 \noalign{\hbox{it is seen that}}
 L^{-1}[N_0s^{-1}[1+(s/c)^{-\nu}]^{-1}]&=N_0\sum_{k=0}^{\infty}(-1)^kc^{\nu
 k}L^{-1}\{s^{-\nu k-1}\}&(13.5)\cr
 &=N_0\sum_{k=0}^{\infty}(-1)^kc^{\nu k}{{t^{\nu
 k}}\over{\Gamma(1+\nu k)}}=N_0E_{\nu}(-c^{\nu}t^{\nu}).&(13.6)\cr}
 $$This completes the proof of Theorem 13.1.

 \vskip.3cm
 \noindent
 {\bf Remark 13.1}.\hskip.3cm If we apply the operator
 ${_0D}_t^{\nu}$ from the left to (13.1) and make use of the formula

 $${_0D}_t^{\nu}[1]={{1}\over{\Gamma(1-\nu)}}t^{-\nu},
 t>0, \Re(\nu)<1\eqno(13.7)
 $$we obtain the fractional diffusion equation

 $${_0D}_t^{\nu}N(t)-N_0{{t^{-\nu}}\over{\Gamma(1-\nu)}}=-c^{\nu}N(t),
 t>0, \Re(\nu)<1\eqno(13.8)
 $$whose solution is also given by (13.6).

 \vskip.3cm
 \noindent
 {\bf Remark 13.2}.\hskip.3cm We note that Haubold and Mathai (2000)
 have given the solution of (13.1) in terms of the series given by
 (13.5). The solution in terms of the Mittag-Leffler function is given in
 Saxena, Mathai and Haubold (2002).

 \vskip.3cm
 \noindent
 {\bf Alternate procedure}.
 \hskip.3cm We now present an alternate method similar to that
 followed by Al-Saqabi and Tuan (2006) for solving some differintegral
 equations, also, see Saxena and Kalla (2008) for details.
 \vskip.2cm Applying the operator $(-c^{\nu})^m{_0D}_t^{-m\nu}$ to
 both sides of (13.1) we find that
 $$(-c^{\nu})^m{_0D}_t^{-m\nu}N(t)-(-c^{\nu})^{m+1}{_0D}_t^{-\nu(m+1)}N(t)=N_0~~{_0D}_t^{-m\nu}[1],
 m=0,1,2,...\eqno(13.9)
 $$Summing up the above expression with respect to $m$ from $0$ to
 $\infty$, it gives

 $$\sum_{m=0}^{\infty}(-c^{\nu})^m{_0D}_t^{-m\nu}N(t)-\sum_{m=0}^{\infty}(-c^{\nu})^{m+1}{_0D}_t^{-\nu(m+1)}N(t)
 =N_0\sum_{m=0}^{\infty}(-c^{\nu})^m{_0D}_t^{-m\nu}[1]
 $$which can be written as

 $$\sum_{m=0}^{\infty}(-c^{\nu})^m{_0D}_t^{-m\nu}N(t)-\sum_{m=1}^{\infty}(-c^{\nu})^m{_0D}_t^{-\nu
 m}N(t)=N_0\sum_{m=0}^{\infty}(-c^{\nu})^m{_0D}_t^{-m\nu}[1].
 $$Simplifying the above equation by using the result

 $${_0D}_t^{-\nu}t^{\mu-1}={{\Gamma(\mu)}\over{\Gamma(\mu+\nu)}}t^{\mu+\nu-1}\eqno(13.10)
 $$where $\min\{\Re(\nu),\Re(\mu)\}>0$, we obtain

 $$N(t)=N_0\sum_{m=0}^{\infty}(-c^{\nu})^{m}{{t^{m\nu}}\over{\Gamma(1+m\nu)}}\eqno(13.11)
 $$for $m=0,1,2,...$ Rewriting the series on the right in terms
of the Mittag-Leffler
 function, it yields the desired result (13.6). The next theorem can
 be proved in a similar manner.

 \vskip.3cm
 \noindent
 {\bf Theorem 13.2}.\hskip.3cm{\it If $\min\{\Re(\nu),\Re(\mu)\}>0$
 then the solution of the integral equation

 $$\eqalignno{N(t)-N_0t^{\mu-1}&=-c^{\nu}{_0D}_t^{-\nu}N(t)&(13.12)\cr
 \noalign{\hbox{is given by}}
 N(t)&=N_0\Gamma(\mu)t^{\mu-1}E_{\nu,\mu}(-c^{\nu}t^{\nu}),&(13.13)\cr}
 $$where $E_{\nu,\mu}(t)$ is the generalized Mittag-Leffler function
 defined in (1.2).}

 \vskip.3cm
 \noindent
 {\bf Proof.}\hskip.3cm Applying Laplace transform to both sides of
 (13.12) and using (1.11), it gives

 $$\eqalignno{\tilde{N}(s)&=\{N(t);s\}=N_0\Gamma(\mu)s^{-\mu}[1+(s/c)^{-\nu}]^{-1}.&(13.14)\cr
 \noalign{\hbox{Using the relation (13.4), it is seen that}}
 L^{-1}\{N_0\Gamma(\mu)s^{-\mu}[1+(s/c)^{-\nu}]^{-1}\}&=
 N_0\sum_{k=0}^{\infty}(-1)^kc^{\nu k}L^{-1}\{s^{-\mu-\nu k}\}\cr
 &=N_0\Gamma(\mu)t^{\mu-1}\sum_{k=0}^{\infty}(-1)^kc^{\nu k}{{t^{\nu
 k}}\over{\Gamma(\mu+\nu k)}}\cr
 &=N_0\Gamma(\mu)t^{\mu-1}E_{\nu,\mu}(-c^{\nu}t^{\nu}).&(13.15)\cr}
 $$This completes the proof of Theorem 13.2.

 \vskip.3cm
 \noindent
 {\bf Alternate procedure}.\hskip.3cm We now give an alternate
 method similar to that followed by Al-Saqabi and Tuan (2006) for
 solving the differintegral equations. Applying the operator
 $(-c^{\nu})^m{_0D}_t^{-m\nu}$ to both sides of (13.12), we find
 that

 $$(-c^{\nu})^m{_0D}_t^{-m\nu}N(t)-(-c^{\nu})^{m+1}{_0D}_t^{-\nu(m+1)}N(t)=N_0(-c^{\nu})^m{_0D}_t^{-m\nu}t^{\mu-1}\eqno(13.16)
 $$for $m=0,1,2,...$. Summing up the above expression with respect
 to $m$ from $0$ to $\infty$, it gives
 $$\eqalignno{\sum_{m=0}^{\infty}(-c^{\nu})^m{_0D}_t^{-m\nu}N(t)&-\sum_{m=0}^{\infty}(-c^{\nu})^{m+1}{_0D}_t^{-\nu(m+1)}N(t)\cr
 &=N_0\sum_{m=0}^{\infty}(-c^{\nu})^m{_0D}_t^{-m\nu}t^{\mu-1}\cr}
 $$which can be written as

 $$\sum_{m=0}^{\infty}(-c^{\nu})^m{_0D}_t^{-m\nu}N(t)-\sum_{m=1}^{\infty}(-c^{\nu})^m{_0D}_t^{-m\nu}N(t)=N_0\sum_{m=0}^{\infty}(-c^{\nu})^m{_0D}_t^{-m\nu}t^{\mu-1}.
 $$Simplifying by using the result (13.10) we obtain

 $$N(t)=N_0\Gamma(\mu)\sum_{m=0}^{\infty}(-c^{\nu})^m{{t^{\mu+m\nu-1}}\over{\Gamma(m\nu+\mu)}};
 m=0,1,2,...\eqno(13.17)
 $$Rewriting the series on the right of (13.17) in terms of the
 generalized Mittag-Leffler function, it yields the desired result
 (13.5). Next we present a general theorem given by Saxena et
 al.(2004).

 \vskip.3cm
 \noindent
 {\bf Theorem 13.3.}\hskip.3cm{\it If $c>0, \Re(\nu)>0$ then for the
 solution of the integral equation

 $$N(t)-N_0f(t)=-c^{\nu}{_0D}_t^{-\nu}N(t)\eqno(13.18)
 $$where $f(t)$ is any integrable function on the finite interval
 $[0,b]$, there exists the formula

 $$N(t)=c~N_0\int_0^tH_{1,2}^{1,1}\left[c^{\nu}(t-\tau)^{\nu}\bigg\vert_{(-1/\nu,1),(-1,\nu)}^{(-1/\nu,1)}\right]
 f(\tau){\rm
 d}\tau,\eqno(13.19)
 $$where $H_{1,2}^{1,1}(\cdot)$ is the H-function defined by (8.1).}
 \vskip.2cm
 The proof can be developed by identifying the Laplace transform of
 $N()+c^{\nu}{_0D}_t^{-\nu}N(t)$ as an H-function and then using the
 convolution property for the Laplace transform.  In what follows, $E_{\beta,\gamma}^{\delta}(\cdot)$ will
 be employed to denote the generalized Mittag-Leffler function,
 defined by (11.1).

 \vskip.3cm
 \noindent
 {\bf Note 13.1.}\hskip.3cm For an alternate derivation of  this
 theorem see Saxena and Kalla (2008).

 \vskip.2cm Next we will discuss time-fractional diffusion.

 \vskip.3cm
 \noindent
 {\bf 14.\hskip.3cm Application to time-fractional diffusion}

 \vskip.3cm
 \noindent
 {\bf Theorem 14.1}.\hskip.3cm{\it Consider the following
 time-fractional diffusion equation

 $${{\partial^{\alpha}}\over{\partial
 t^{\alpha}}}N(x,t)=D{{\partial^2}\over{\partial x^2}}N(x,t),
 0<\alpha<1,\eqno(14.1)
 $$where $D$ is the diffusion constant and $N(x,t=0)=\delta(x)$ is
 the Dirac delta function and
 $\lim_{x\rightarrow\pm\infty}N(x,t)=0$. Then its fundamental
 solution is given by

 $$N(x,t)={{1}\over{|x|}}H_{1,1}^{1,0}\left[{{|x|^2}\over{D
 t^{\alpha}}}\bigg\vert_{(1,2)}^{(1,\alpha)}\right].\eqno(14.2)
 $$}

 \vskip.3cm
 \noindent
 {\bf Proof}.\hskip.3cm In order to find a closed form
 representation of the solution in terms of the H-function, we use
 the method of joint Laplace-Fourier transform, defined by

 $${\tilde{N}}^{*}(k,s)=\int_0^{\infty}\int_{-\infty}^{\infty}{\rm
 e}^{-st+ikx}N(x,t){\rm d}x~{\rm d}t\eqno(14.3)
 $$where, according to the convention followed, $\sim $ will
 denote the Laplace transform and * the Fourier transform.
 Applying the Laplace transform with respect to time variable $t$,
 Fourier transform with respect to space variable $x$ and using the
 given condition, we find that

 $$\eqalignno{s^{\alpha}{\tilde{N}}^{*}(k,s)-s^{\alpha-1}&=-Dk^2{\tilde{N}}^{*}(k,s)\cr
 \noalign{\hbox{and then}}
 {\tilde{N}}^{*}(k,s)&={{s^{\alpha-1}}\over{s^{\alpha}+Dk^2}}.\cr}$$Inverting the Laplace transform, it yields

 $$N^{*}(k,t)=L^{-1}\{{{s^{\alpha-1}}\over{s^{\alpha}+Dk^2}}\}=E_{\alpha}(-Dk^2t^{\alpha}),\eqno(14.4)
 $$where $E_{\alpha}(\cdot)$ is the Mittag-Leffler function defined
 by (1.1). In order to invert the Fourier transform, we will make
 use of the integral

 $$\int_0^{\infty}\cos(kt)E_{\alpha,\beta}(-at^2){\rm
 d}t={{\pi}\over{k}}H_{1,1}^{1,0}\left[{{k^2}\over{a}}\bigg\vert_{(1,2)}^{(\beta,\alpha)}\right],\eqno(14.5)
 $$where $\Re(\alpha)>0, \Re(\beta)>0, k>0,a>0$; and the formula

 $${{1}\over{2\pi}}\int_{-\infty}^{\infty}{\rm e}^{-ikx}f(k){\rm
 d}k={{1}\over{\pi}}\int_0^{\infty}f(k)\cos(kx){\rm d}k.\eqno(14.6)
 $$Then it yields the required solution as

 $$N(x,t)={{1}\over{|x|}}H_{3,3}^{2,1}\left[{{|x|^2}\over{Dt^{\alpha}}}\big\vert_{(1,2),(1,1),(1,1)}^{(1,1),(1,\alpha),(1,1)}\right]
 ={{1}\over{|x|}}H_{1,1}^{1,0}\left[{{|x|^2}\over{Dt^{\alpha}}}\big\vert_{(1,2)}^{(1,\alpha)}\right].\eqno(14.7)
 $$

 \vskip.3cm
 \noindent
 {\bf Note 14.1.}\hskip.3cm For $\alpha=1$, (14.7) reduces to the
 Gaussian density

 $$N(x,t)={{1}\over{2(\pi
 Dt)^{1\over2}}}\exp(-{{|x|^2}\over{4Dt}}).\eqno(14.8)
 $$Fractional space-diffusion will be discussed in the next section.

 \vskip.3cm
 \noindent
 {\bf 15.\hskip.3cm Application to fractional-space diffusion}

 \vskip.3cm
 \noindent
 {\bf Theorem 15.1}.\hskip.3cm{\it Consider the following
 fractional space-diffusion equation

 $${{\partial}\over{\partial
 t}}N(x,t)=D{{\partial^{\alpha}}\over{\partial x^{\alpha}}}N(x,t),
 0<\alpha <1,\eqno(15.1)
 $$where $D$ is the diffusion constant,
 ${{\partial^{\alpha}}\over{\partial x^{\alpha}}}$ is the operator
 defined by (12.20) and $N(x,t=0)=\delta(x)$ is the Dirac delta
 function and $\lim_{x\rightarrow\pm\infty}N(x,t)=0$. Then its
 fundamental solution is given by

 $$N(x,t)={{1}\over{a|x|}}H_{2,2}^{1,1}\left[{{|x|}\over{(Dt)^{{1}\over{\alpha}}}}
 \bigg\vert_{(1,1),(1,{1\over2})}^{(1,{{1}\over{\alpha}}),(1,{1\over2})}\right].\eqno(15.2)
 $$}The proof can be developed on similar lines to that of the
 theorem of the preceding section.

 \vskip.3cm
 \noindent
 {\bf 16.\hskip.3cm Application to fractional reaction-diffusion
 model}

 \vskip.3cm In the same way, we can establish the following theorem,
 which gives the fundamental solution of the reaction-diffusion
 model given below.

 \vskip.3cm
 \noindent
 {\bf Theorem 16.1}.\hskip.3cm {\it Consider the following
 reaction-diffusion model

 $${{\partial^{\beta}}\over{\partial
 t^{\beta}}}N(x,t)=\eta~~{_{-\infty}D}_x^{\alpha}N(x,t), 0<\beta\le
 1\eqno(16.1)
 $$with the initial condition $N(x,t=0)=\delta(x),
 \lim_{x\rightarrow\pm\infty}N(x,t)=0$, where $\eta$ is a diffusion
 constant and $\delta(x)$ is the Dirac delta function. Then for the
 solution of (16.1) there holds the formula

 $$N(x,t)={{1}\over{a|x|}}H_{3,3}^{2,1}\left[{{|x|}\over{(\eta
 t^{\beta})^{{1}\over{\alpha}}}}\bigg\vert_{(1,1),(1,{{1}\over{\alpha}}),
 (1,{1\over2})}^{(1,{{1}\over{\alpha}}),
 (1,{{\beta}\over{\alpha}}),(1,{1\over2})}\right].\eqno(16.2)
 $$}For details of the proof, the reader is referred to the original
 paper by Saxena, Mathai and Haubold (2006).

 \vskip.3cm
 \noindent
 {\bf Corollary 16.1}.\hskip.3cm{\it For the solution of the fraction
 reaction-diffusion equation

 $${{\partial}\over{\partial t}}N(x,t)=\eta~~
 {_{-\infty}D}_x^{\alpha}N(x,t),\eqno(16.3)
 $$with initial condition $N(x,t=0)=\delta(x)$ there holds the
 formula

 $$N(x,t)={{1}\over{\alpha|x|}}H_{2,2}^{1,1}\left[{{|x|}\over{(\eta
 t)^{{1}\over{\alpha}}}}\bigg\vert_{(1,1), (1,{1\over2})}^{(1,{{1}\over{\alpha}}),(1,{1\over2})}\right]\eqno(16.4)
 $$where $\alpha>0$.}

 \vskip.3cm
 \noindent
 {\bf Note 16.1.}\hskip.3cm It may be noted that (16.4) is a closed
 form representation of a L\'evy stable law. It is interesting to
 note that as $\alpha\rightarrow 2$ the classical Gaussian solution
 is recovered since

 $$\eqalignno{N(x,t)&={{1}\over{2|x|}}H_{2,2}^{1,1}\left[{{|x|}\over{(\eta
 t)^{{1}\over{\alpha}}}}\bigg\vert_{(1,1), (1,{1\over2})}^{
 (1,{1\over2}), (1,{1\over2})}\right]\cr
 &={{1}\over{2\pi^{1\over2}|x|}}\sum_{k=0}^{\infty}{{(-1)^k}\over{k!}}\left[{{|x|}\over{2(\eta
 t)^{{1}\over{\alpha}}}}\right]^{2k+1}&(16.5)\cr
 &=(4\pi (\eta
 t)^{{2}\over{\alpha}})^{-{1\over2}}\exp\left[-{{|x|^2}\over{4(\eta
 t)^{{2}\over{\alpha}}}}\right].&(16.6)\cr}
 $$

 \vskip.3cm
 \noindent
 {\bf 17.\hskip.3cm Application to nonlinear waves}

 \vskip.3cm It will be shown in this section that by the application
 of the inverse Laplace transforms of certain algebraic functions
 derived in Saxena, Mathai and Haubold (2006), we can establish the
 following theorem for nonlinear waves.

 \vskip.3cm
 \noindent
 {\bf Theorem 17.1.}\hskip.3cm{\it Consider the fractional
 reaction-diffusion equation

 $$\eqalignno{{_0D}_t^{\alpha}N(x,t)&+a{_0D}_t^{\beta}N(x,t)\cr
 &=\nu^2{_{-\infty}D}_x^{\gamma}N(x,t)+\zeta^2N(x,t)+\phi(x,t)&(17.1)\cr}
 $$for $x\in R, t>0, 0\le \alpha\le 1, 0\le \beta\le 1$ with initial
 conditions $N(x,0)=f(x)$, $\lim_{x\rightarrow\pm\infty}N(x,t)=0$
 for $x\in R$, where $\nu^2$ is a diffusion constant, $\zeta$ is a
 constant which describes the nonlinearity in the system, and
 $\phi(x,t)$ is nonlinear function which belongs to the area of
 reaction-diffusion, then there holds the following formula for the
 solution of (17.1).

$$\eqalignno{N(x,t)&=\sum_{r=0}^{\infty}{{(-a)^r}\over{2\pi}}\int_{-\infty}^{\infty}t^{(\alpha-\beta)r}f^{*}(k)\exp(-kx)\cr
&\times
[E_{\alpha,(\alpha-\beta)r+1}^{r+1}(-bt^{\alpha})+t^{\alpha-\beta}E_{\alpha,(\alpha-\beta)(r+1)+1}^{r+1}(-bt^{\alpha})]{\rm
d}k\cr
&+\sum_{r=0}^{\infty}{{(-a)^r}\over{2\pi}}\int_0^t\zeta^{\alpha+(\alpha-\beta)r-1}\int_{-\infty}^{\infty}\phi(k,t-\zeta)\exp(-ikx)\cr
&\times
E_{\alpha,\alpha+(\alpha-\beta)r}^{r+1}(-b\zeta^{\alpha}){\rm
d}k~{\rm d}\zeta&(17.2)\cr}
$$where $\alpha>\beta$ and $E_{\beta,\gamma}^{\delta}(\cdot)$ is the
generalized Mittag-Leffler function, defined by (11.1), and
$b=\nu^2|k|^{\gamma}-\zeta^2.$}

\vskip.3cm \noindent {\bf Proof.}\hskip.3cm Applying the Laplace
transform with respect to the time variable $t$ and using the
boundary conditions, we find that

$$\eqalignno{s^{\alpha}\tilde{N}(x,s)&-s^{\alpha-1}f(x)+as^{\beta}\tilde{N}(x,s)-as^{\beta-1}f(x)\cr
&=\nu^2{_{-\infty}D}_x^{\gamma}\tilde{N}(x,s)+\zeta^2\tilde{N}(x,s)+\tilde{f}(x,s).&(17.3)\cr}
$$If we apply the Fourier transform with respect to the space
variable $x$ to (17.3) it yields

$$\eqalignno{s^{\alpha}{\tilde{N}}^{*}(k,s)&-s^{\alpha-1}f^{*}(k)+as^{\beta}{\tilde{N}}^{*}(k,s)-as^{\beta-1}f^{*}(k)\cr
&=-\nu^2|k|^{\gamma}{\tilde{N}}^{*}(k,s)+\zeta^2{\tilde{N}}^{*}(k,s)+{\tilde{f}}^{*}(k,s).&(17.4)\cr}
$$Solving for ${\tilde{N}}^{*}$ it gives

$${\tilde{N}}^{*}(k,s)={{(s^{\alpha-1}+as^{\beta-1})f^{*}(k)+{\tilde{f}}^{*}(k,s)}\over{s^{\alpha}+as^{\beta}+b}}\eqno(17.5)
$$where $b=\nu^2|k|^{\gamma}-\zeta^2$. For inverting (17.5) it is
convenient to first invert the Laplace transform and then the
Fourier transform. Inverting the Laplace transform with the help of
the result Saxena et al.(2004a, eq(28))

$$L^{-1}\left\{{{s^{\rho-1}}\over{s^{\alpha}+as^{\beta}+b}};t\right\}
=t^{\alpha-\rho}\sum_{r=0}^{\infty}(-a)^rt^{(\alpha-\beta)r}E_{\alpha,\alpha+(\alpha-\beta)r-\rho+1}^{r+1}(-bt^{\alpha}),\eqno(17.6)
$$where $\Re(\alpha)>0, \Re(\beta)>0, \Re(\rho)>0$,
$|{{as^{\beta}}\over{s^{\alpha}+b}}|<1$ and provided that the series
in (17.6) is convergent, it yields

$$\eqalignno{N^{*}(k,t)&=\sum_{r=0}^{\infty}(-a)^rt^{(\alpha-\beta)r}f^{*}(k)\cr
&\times[E_{\alpha,(\alpha-\beta)r+1}^{r+1}(-bt^{\alpha})+t^{\alpha-\beta}E_{\alpha,(\alpha-\beta)(r+1)+1}^{r+1}(-bt^{\alpha})]\cr
&+\sum_{r=0}^{\infty}(-a)^r\int_0^t\phi^{*}(k,t-\zeta)\zeta^{\alpha+(\alpha-\beta)r-1}E_{\alpha,(\alpha-\beta)r+\alpha}^{r+1}(-b\zeta^{\alpha}){\rm
d}\zeta.&(17.7)\cr}
$$Finally, the inverse Fourier transform gives the desired solution
(17.3).

\vskip.3cm \noindent {\bf 18.\hskip.3cm Generalized Mittag-Leffler
type functions}

\vskip.3cm The multiindex ($m$-tuple) Mittag-Leffler function is
defined in Kiryakova (2000) by means of the power series

$$E_{({{1}\over{\rho_i}}),(\mu_i)}(z)=\sum_{k=0}^{\infty}\phi_kz^k=\sum_{k=0}^{\infty}{{z^k}\over{\prod_{j=1}^m\Gamma(\mu_j+{{k}\over{\rho_j}})}}.\eqno(18.1)
$$Here $m>1$ is an integer, $\rho_1,...,\rho_m$ and
$\mu_1,...,\mu_m$ are arbitrary real parameters. \vskip.2cm The
following theorem is proved by Kiryakova (2000, p.244) which shows
that the multiindex Mittag-Leffler function is an entire function
and also gives its asymptotic estimate, order and type.

\vskip.3cm \noindent {\bf Theorem 18.1.}\hskip.3cm {\it For
arbitrary sets of indices $\rho_i>0, -\infty<\mu_i<\infty,
i=1,...,m$ the multiindex Mittag-Leffler function defined by (18.1)
is an entire function of order

$$\eqalignno{\rho&=[\sum_{i=1}^m{{1}\over{\rho_i}}]^{-1},
\hbox{  that is,~~
}{{1}\over{\rho}}={{1}\over{\rho_1}}+...+{{1}\over{\rho_m}}&(18.2)\cr
\noalign{\hbox{and type}}
\sigma&=({{\rho_1}\over{\rho}})^{{\rho}\over{\rho_1}}...({{\rho_m}\over{\rho}})^{{\rho}\over{\rho_m}}.&(18.3)\cr}
$$}Furthermore, for every positive $\epsilon$, the asymptotic
estimate

$$|E_{(1/{\rho_i}),(\mu_i)}(z)|<\exp((\sigma+\epsilon)|z|^{\rho}),\eqno(18.4)
$$holds for $|z|\ge r_0(\epsilon), r_0(\epsilon)$ sufficiently
large. \vskip.2cm It is interesting to note that for $m=2$, (18.2)
reduces to the generalized Mittag-Leffler function considered by
Dzherbashyan (1960) denoted by $\phi_{\rho_1,\rho_2}(z;\mu_1,\mu_2)$
and defined in the following form (Kiryakova, 1994, Appendix)

$$\eqalignno{E_{(1/{\rho_1},1/{\rho_2});(\mu_1,\mu_2)}(z)&=\phi_{\rho_1,\rho_2}(z;\mu_1,\mu_2)\cr
&=\sum_{k=0}^{\infty}{{z^k}\over{\Gamma(\mu_1+{{k}\over{\rho_1}})\Gamma(\mu_2+{{k}\over{\rho_2}})}},&(18.5)\cr}
$$and shown to be an entire function of order

$$\rho={{\rho_1\rho_2}\over{\rho_1+\rho_2}}\hbox{  and type
}\sigma=({{\rho_1}\over{\rho_2}})^{{\rho_2}\over{\rho_1}}({{\rho_2}\over{\rho_1}})^{{\rho_1}\over{\rho_2}}.\eqno(18.6)
$$Relations between multiindex Mittag-Leffler function with
H-function, generalized Wright function and other special functions
are given by Kiryakova; for details, see the original papers
Kiryakova (1999, 2000). Saxena, Kalla and Kiryakova (2003)
investigated the relations between the multiindex Mittag-Leffler
function and the Riemann-Liouville fractional integrals and
derivatives. The results derived are of general nature and give rise
to a number of known as well as unknown results in the theory of
generalized Mittag-Leffler functions, which serve as a backbone for
the fractional calculus. Two interesting theorems established by
Saxena et al. (2003) are described below.

\vskip.3cm \noindent {\bf Theorem 18.2.}\hskip.3cm {\it Let
$\alpha>0,\rho_i>0,\mu_i>0, i=1,...,m$ and further, let
$I_{0_{+}}^{\alpha}$ be the left-sided Riemann-Liouville fractional
integral. Then there holds the relation

$$\eqalignno{(I_{0_{+}}^{\alpha}[t^{\rho_i-1}&E_{(1/{\rho_i}),(\mu_i)}(at^{1/{\rho_i}})])(x)\cr
&=x^{\alpha+\rho_i-1}E_{(1/{\rho_i}),(\mu_1+\alpha,\mu_2,...,\mu_m)}(at^{{1}\over{\rho_i}}).&(18.7)\cr}
$$}

\vskip.3cm \noindent {\bf Theorem 18.3}.\hskip.3cm{\it Let
$\alpha>0,\rho_i>0,\mu_i>0, i=1,...,m$ and further, let
$I_{-}^{\alpha}$ be the right-sided Riemann-Liouville fractional
integral. Then there holds the relation

$$\eqalignno{(I_{-}^{\alpha}[t^{-\rho_i-\alpha}&E_{(1/{\rho_i}),(\mu_i)}(at^{-{{1}\over{\rho_i}}})])(x)\cr
&=x^{-\mu_i}E_{(1/{\rho_i}),(\mu_1+\alpha,\mu_2,...,\mu_m)}(at^{-{{1}\over{\rho_i}}}).&(18.8)\cr}
$$}

\vskip.2cm Another generalization of the Mittag-Leffler function is
recently given by Sharma (2008) in terms of the M-series defined by

$$\eqalignno{{_pM}_q^{\alpha}&={_pM}_q^{\alpha}(a_1,...,a_p;b_1,...,b_q
;\alpha;z)\cr
&=\sum_{r=0}^{\infty}{{(a_1)_r...(a_p)_r}\over{(b_1)_r...(b_q)_r\Gamma(\alpha
r+1)}}z^r,~~p\le q+1.&(18.9)\cr}
$$

\vskip.3cm \noindent {\bf Remark 18.1.}\hskip.3cm According to
Saxena (2009), the M-series discussed by Sharma (2008) is not a new
special function. It is, in disguise, a special case of the
generalized Wright function ${_p\psi_q}(z)$, which was introduced by
Wright (1935), as shown below.

$$\eqalignno{\kappa~~~{_{p+1}\psi}_{q+1}\left[\matrix{(a_1,1),...,(a_p,1),(1,1)\cr
(b_1,1),...,(b_q,1),(1,\alpha)\cr};z\right]&=\sum_{r=0}^{\infty}{{(a_1)_r...(a_p)_r(1)_r}\over{(b_1)_r...(b_q)_r\Gamma(\alpha
r+1)}}{{z^r}\over{r!}}\cr
&=\sum_{r=0}^{\infty}{{(a_1)_r...(a_p)_r}\over{(b_1)_r...(b_q)_r\Gamma(\alpha
r+1)}}z^r\cr
&={_pM}_q^{\alpha}(a_1,...,a_p;b_1,...,b_q;\alpha;z)&(18.10)\cr}
$$where

$$\kappa={{\prod_{j=1}^q\Gamma(b_j)}\over{\prod_{j=1}^p\Gamma(a_j)}}.
$$Fractional integration and fractional differentiation of the
M-series are discussed by Sharma (2008). The two results proved in
Sharma (2008) for the function defined by (18.9) are reproduced
below. For $\Re(\nu)>0$

$$\eqalignno{(I_{0_{+}}^{\nu}[{_pM}_q^{\alpha}(a_1,...,a_p;&b_1,...,b_q;\alpha;z)])(x)\cr
&={{z^{\nu}}\over{\Gamma(1+\nu)}}{_{p+1}M}_{q+1}^{\alpha}(a_1,...,a_p,1;b_1,...,b_q,1+\nu;\alpha;z)\cr
\noalign{\hbox{and for $\Re(\nu)<0$}}
(D_{0_{+}}^{\nu}[{_pM}_q^{\alpha}(a_1,...,a_p;&b_1,...,b_q;\alpha;z)])(x)\cr
&={{z^{\nu}}\over{\Gamma(1-\nu)}}{_{p+1}M}_{q+1}^{\alpha}(a_1,...,a_p,1;b_1,...,b_q,1-\nu;\alpha;z).\cr}
$$

\vskip.5cm \noindent{\bf 19.\hskip.3cm Mittag-Leffler Distributions
and Processes}

\vskip.3cm\noindent{\bf 19.1.\hskip.3cm Mittag-Leffler Statistical
Distribution and Its Properties}

\vskip.3cm A statistical distribution in terms of the Mittag-Leffler
function $E_{\alpha}(y)$ was defined by Pillai (1990) in terms of
the distribution function or cumulative density function as follows:

$$G_y(y)=1-E_{\alpha}(-y^{\alpha})=\sum_{k=1}^{\infty}{{(-1)^{k+1}y^{\alpha
k}}\over{\Gamma(1+\alpha k)}},0<\alpha\le 1, y>0\eqno(19.1.1)
$$and $G_y(y)=0$ for $y\le 0$. Differentiating on both sides with
respect to $y$ we obtain the density function $f(y)$ as follows:

$$\eqalignno{f(y)&={{{\rm d}}\over{{\rm d}y}}G_y(y)\cr
&={{{\rm d}}\over{{\rm
d}y}}\left[\sum_{k=1}^{\infty}{{(-1)^{k+1}y^{\alpha
k}}\over{\Gamma(1+\alpha
k)}}\right]=\sum_{k=1}^{\infty}{{(-1)^{k+1}\alpha k y^{\alpha
k-1}}\over{\Gamma(1+\alpha k)}}\cr
&=\sum_{k=1}^{\infty}{{(-1)^{k+1}y^{\alpha k}}\over{\Gamma(1+\alpha
k)}}=\sum_{k=0}^{\infty}{{(-1)^ky^{\alpha+\alpha
k-1}}\over{\Gamma(\alpha+\alpha k)}}\cr \noalign{\hbox{by replacing
$k$ by $k+1$}}
&=y^{\alpha-1}E_{\alpha,\alpha}(-y^{\alpha}),~0<\alpha\le 1,
y>0&(19.1.2)\cr}
$$where $E_{\alpha,\beta}(x)$ is the generalized Mittag-Leffler
function.
 \vskip.2cm It is straightforward to observe that for the
density in (19.1.2) the distribution function is that in (19.1.1).
The Laplace transform of the density in (19.1.2) is the following:
$$\eqalignno{L_f(t)=\int_0^{\infty}{\rm e}^{-tx}f(x){\rm
d}x&=\int_0^{\infty}{\rm
e}^{-tx}x^{\alpha-1}E_{\alpha,\alpha}(-x^{\alpha}){\rm d}x\cr
&=(1+t^{\alpha})^{-1}, |t^{\alpha}|<1.&(19.1.3)\cr}
$$Note that (19.1.3) is a special case of the general class of Laplace
transforms discussed in Section 2.3.7 [Haubold and Mathai (2008)].
From (19.1.3) one can also note that there is a structural
representation in terms of positive L\'evy distribution. A positive
L\'evy random variable $u>0$, with parameter $\alpha$ is such that
the Laplace transform of the density of $u>0$ is given by ${\rm
e}^{-t^{\alpha}}$. That is,

$$E[{\rm e}^{-tu}]={\rm e}^{-t^{\alpha}}\eqno(19.1.4)
$$ where $E(\cdot)$ denotes the  expected value of
$(\cdot)$ or the statistical expectation of $(\cdot)$. When
$\alpha=1$ the random variable is degenerate with the density
function

$$f_1(x)=\cases{1, \hbox{ for  } x=1\cr
0,\hbox{  elsewhere.}\cr}
$$Consider an exponential random variable with density
function

$$f_1(x)=\cases{{\rm e}^{-x},\hbox{  for  }0\le x<\infty\cr
0,\hbox{  elsewhere;}\cr}~~~~~~L_{f_1}(t)=(1+t)^{-1}\eqno(19.1.5)
$$and with the Laplace transform $L_{f_1}(t)$.

\vskip.3cm \noindent{\bf Theorem 19.1.1}\hskip.3cm{\it Let $y>0$ be
a L\'evy random variable with Laplace transform as in (19.1.4) and
let $x$ and $y$ be independently distributed. Then
$u=yx^{{1}\over{\alpha}}$ is distributed as a Mittag-Leffler random
variable with Laplace transform as in (19.1.3).}

\vskip.3cm \noindent{\bf Proof}.\hskip.3cm For establishing this
result we will make use of a basic result on conditional
expectations, which will be stated as a lemma.

\vskip.3cm \noindent{\bf Lemma 19.1.1}\hskip.3cm{\it For two random
variables $x$ and $y$ having a joint distribution,

$$E(x)=E[E(x|y)]\eqno(19.1.6)
$$whenever all the expected values exist, where the inside
expectation is taken in the conditional space of $x$ given $y$ and
the outside expectation is taken in the marginal space of $y$.}

\vskip.3cm Now by applying (19.1.6) we have the following: Let the
density of $u$ be denoted by $g(u)$. Then the Laplace transform of
$g$ is given by

$$E[{\rm e}^{-(tx^{{1}\over{\alpha}})y}|x]={\rm
e}^{-t^{\alpha}x.}\eqno(19.1.8)
$$But the right side of (19.1.8) is in the form of a Laplace
transform of the density of $x$ with parameter $t^{\alpha}$. Hence
the expected value of the right side is

$$L_g(t)=(1+t^{\alpha})^{-1}\eqno(19.1.9)
$$which establishes the result. From (19.1.8) one property is
obvious. Suppose that we consider an arbitrary random variable $y$
with the Laplace transform of the form

$$L_g(t)={\rm e}^{-[\phi(t)]}\eqno(19.1.10)
$$whenever the expected value exists, where $\phi(t)$ be such that
$$\phi(tx^{{1}\over{\alpha}})=x\phi(t), \lim_{t\rightarrow
0}\phi(t)=0.
$$Then from (19.1.8) we have

$$E[{\rm e}^{-(tx^{{1}\over{\alpha}})y}|x]={\rm
e}^{-x[\phi(t)]}.\eqno(19.1.11) $$Now, let $x$ be an arbitrary
positive random variable having Laplace transform, denoted by
$L_x(t)$ where $L_x(t)=\psi(t)$. Then from (19.1.9) we have

$$L_g(t)=\psi[\phi(t)].\eqno(19.1.12)
$$For example, if $y$ is a random variable whose density has the Laplace
transform, denoted by $L_y(t)=\phi(t)$, with
$\phi(tx^{{1}\over{\alpha}})=x\phi(t)$, and if $x$ is a real random
variable having the gamma density,

$$f_x(x)={{x^{\beta-1}{\rm
e}^{-{{x}\over{\delta}}}}\over{\delta^{\beta}\Gamma(\beta)}}, x\ge
0, \beta>0,\delta>0\eqno(19.1.13)
$$and $f_x(x)=0$ elsewhere, and if $x$ and $y$ are statistically
independently distributed and if $u=yx^{{1}\over{\alpha}}$ then the
Laplace transform of the density of $u$, denoted by $L_u(t)$ is
given by
$$L_u(t)=[1+\delta\{\phi(t)\}]^{-\beta}.\eqno(19.1.14)
$$

\vskip.3cm \noindent {\bf Note 19.1.1.}\hskip.3cm Since we did not
put any restriction on the nature of the random variables, except
that the expected values exist, the result in (19.1.12) holds
whether the variables are continuous, discrete or mixed.

\vskip.3cm \noindent{\bf Note 19.1.2}.\hskip.3cm Observe that for
the result in (19.1.12) to hold we need only the conditional Laplace
transform of $y$ given $x$ be of the form in (19.1.11) and the
marginal Laplace transform of $x$ be $\psi(t)$. Then the result in
(19.1.12) will hold. Thus statistical independence of $x$ and $y$ is
not a basic requirement for the result in (19.1.12) to hold

\vskip.2cm Thus from (19.1.13) we may write a particular case as

$$z=yx^{{1}\over{\alpha}}\eqno(19.1.15)
$$where $x$ is distributed as in (19.1.5) and $y$ as in (19.1.4)
then $z$ will be distributed as in (19.1.3) or (19.1.9) when $x$ and
$y$ are assumed to be independently distributed.

\vskip.3cm \noindent{\bf Note 19.1.3}\hskip.3cm The representation
of the Mittag-Leffler variable as well as the properties described
on page 1432 of Jayakumar (2003) and on page 53 of Jayakumar and
Suresh (2003) are to be rewritten and corrected because the
exponential variable and L\'evy variable seem to be interchanged
there.

\vskip.2cm By taking the natural logarithms on both sides of
(19.1.15) we have

$${{1}\over{\alpha}}\ln x+\ln y=\ln z.\eqno(19.1.16)
$$Then the first moment of $\ln z$ is available from (19.1.16) by
computing $E[\ln x]$ and $E[\ln y]$. But $E[\ln x]$ is available
from the following procedure:

$$\eqalignno{E[{\rm e}^{-t\ln x}]&=E[{\rm e}^{\ln
x^{-t}}]=E[x^{-t}]=\int_0^{\infty}x^{-t}{\rm e}^{-x}{\rm d}x\cr
&=\Gamma(1-t)\hbox{  for  }\Re(1-t)>0&(19.1.17)\cr}
$$which will be $\Gamma(\beta-t)/\Gamma(\beta)$ for the density in
(19.1.13). Hence
$$\eqalignno{E[\ln x]&=-{{{\rm d}}\over{{\rm d}t}}E[{\rm e}^{-t\ln
x}]|_{t=0}=-{{{\rm d}}\over{{\rm d}t}}\Gamma(1-t)|_{t=0}.\cr
\noalign{\hbox{But}} {{{\rm d}}\over{{\rm
d}t}}\Gamma(1-t)&=-\Gamma(1-t){{{\rm d}}\over{{\rm d}t}}\ln
\Gamma(1-t)=-\Gamma(1-t)\psi(1-t)&(19.1.18)\cr}
$$where $\psi(\cdot)$ is the psi function of $(\cdot)$, see Mathai
(1993) for details. Hence by taking the limits $t\rightarrow 0$

$$E[\ln
x]=-\Gamma(1-t)\psi(1-t)|_{t=0}=-\psi(1)=\gamma\eqno(19.1.19)
$$where $\gamma$ is Euler's constant, see Mathai (1993) for details.

\vskip.5cm \noindent{\bf 19.2.\hskip.3cm Mellin-Barnes
representation of the Mittag-Leffler density}

\vskip.3cm Consider the density function in (19.1.2). After writing
in series form and then looking at the corresponding Mellin-Barnes
representation we have the following:

$$\eqalignno{g(x)&=x^{\alpha-1}E_{\alpha,\alpha}(-x^{\alpha})=\sum_{k=0}^{\infty}\Gamma(1+k){{(-1)^{k}}\over{k!}}{{x^{\alpha-1+\alpha
k}}\over{\Gamma(\alpha+\alpha k)}}&(19.2.1)\cr &={{1}\over{2\pi
i}}\int_{c-i\infty}^{c+i\infty}{{1}\over{\alpha}}{{\Gamma({{1}\over{\alpha}}-{{s}\over{\alpha}})\Gamma(1-{{1}\over{\alpha}}+{{s}\over{\alpha}})}\over{\Gamma(1-s)}}x^{-s}{\rm
d}s,1-\alpha<c<1&(19.2.2)\cr \noalign{\hbox{[by expanding as the sum
of residues at the poles of
$\Gamma(1-{{1}\over{\alpha}}+{{s}\over{\alpha}})$]}}
&={{1}\over{2\pi
i}}\int_{c_1-i\infty}^{c_1+i\infty}{{\Gamma(s)\Gamma(1-s)}\over{\Gamma(\alpha
s)}}x^{\alpha s-1}{\rm d}s=g(x), 0<c_1<1, 0<\alpha\le 1&(19.2.3)\cr}
$$by putting ${{1}\over{\alpha}}-{{s}\over{\alpha}}=s_1$. Here the
point $s=0$ is removable. By taking the Laplace transform of $g(x)$
from (19.2.1) we have

$$\eqalignno{L_g(t)&=\sum_{k=0}^{\infty}{{(-1)^k}\over{\Gamma(\alpha+\alpha
k)}}\int_0^{\infty}x^{\alpha+\alpha k-1}{\rm e}^{-tx}{\rm d}x\cr
&=\sum_{k=0}^{\infty}(-1)^kt^{-\alpha-\alpha
k}=t^{-\alpha}(1+t^{-\alpha})^{-1}=(1+t^{\alpha})^{-1},|t^{\alpha}|<1.&(19.2.4)\cr}
$$

\vskip.3cm \noindent{\bf 19.2.1.\hskip.3cm Generalized
Mittag-Leffler density}

\vskip.3cm Consider the generalized Mittag-Leffler function

$$\eqalignno{g_1(x)&={{1}\over{\Gamma(\gamma)}}\sum_{k=0}^{\infty}{{(-1)^k\Gamma(\gamma+k)}\over{k!\Gamma(\alpha
k+\alpha \gamma)}}x^{\alpha\gamma-1+\alpha k}\cr
&=x^{\alpha\gamma-1}E_{\alpha\gamma,\alpha}^{\gamma}(-x^{\alpha}),\alpha>0,\gamma>0.&(19.2.6)\cr}
$$Laplace transform of $g_1(x)$ is the following:

$$\eqalignno{L_{g_1}(t)&=\sum_{k=0}^{\infty}{{(-1)^k}\over{k!}}{{(\gamma)_k}\over{\Gamma(\alpha\gamma+\alpha
k)}}\int_0^{\infty}x^{\alpha\gamma+\alpha k-1}{\rm e}^{-tx}{\rm
d}x\cr &=\sum_{k=0}^{\infty}(-1)^k{{(\gamma)_k}\over{k!}}t^{-\alpha
\gamma-\alpha
k}=(1+t^{\alpha})^{-\gamma},|t^{\alpha}|<1.&(19.2.7)\cr}
$$In fact, this is a special case of the general class of Laplace
transforms connected with Mittag-Leffler function considered in
Mathai, et al. (2006).

\vskip.3cm \noindent{\bf 19.3.\hskip.3cm Mittag-Leffler Density as
an H-function}

\vskip.3cm $g_1(x)$ of (19.2.6) can be written as a Mellin-Barnes
integral and then as an H-function.

$$\eqalignno{g_1(x)&={{1}\over{\Gamma(\eta)}}{{1}\over{2\pi
i}}\int_{c-i\infty}^{c+i\infty}{{\Gamma(s)\Gamma(\eta-s)}\over{\Gamma(\alpha\eta-\alpha
s)}}x^{\alpha\eta-1}(x^{\alpha})^{-s}{\rm
d}s,\Re(\eta)>0,0<c<\Re(\eta)\cr
&={{x^{\alpha\eta-1}}\over{\Gamma(\eta)}}H_{1,2}^{1,1}\left[x^{\alpha}
\bigg\vert_{(0,1),(1-\alpha\eta,\alpha)}^{(1-\eta,1)}\right]&(19.3.1)\cr
&={{1}\over{\alpha\Gamma(\eta)}}{{1}\over{2\pi
i}}\int_{c_1-i\infty}^{c_1+i\infty}
{{\Gamma(\eta-{{1}\over{\alpha}}+{{s}\over{\alpha}})\Gamma({{1}\over{\alpha}}-{{s}\over{\alpha}})}
\over{\Gamma(1-s)}}x^{-s}{\rm d}s&(19.3.2)\cr \noalign{\hbox{[by
taking $\alpha\eta-1-\alpha s=-s_1$]}}
&={{1}\over{\alpha\Gamma(\eta)}}H_{1,2}^{1,1}\left[x\bigg\vert_{(\eta-{{1}\over{\alpha}},
{{1}\over{\alpha}}),(0,1)}^{(1-{{1}\over{\alpha}},{{1}\over{\alpha}})}\right].&(19.3.3)\cr}
$$Since $g$ and $g_1$ are represented as inverse Mellin transforms,
in the Mellin-Barnes representation, one can obtain the $(s-1)$-th
moments of $g$ and $g_1$ from (19.3.2). That is,

$$\eqalignno{M_{g_1}(s)&=E(x^{s-1})\hbox{  in  }g_1\cr
&={{1}\over{\Gamma(\gamma)}}{{\Gamma(\eta-
{{1}\over{\alpha}}+{{s}\over{\alpha}})\Gamma({{1}\over{\alpha}}-{{s}\over{\alpha}})}
\over{\alpha\Gamma(1-s)}},&(19.3.4)\cr \noalign{\hbox{for
$1-\alpha<\Re(s)<1,0<\alpha\le 1, \eta>0.$}}
M_g(t)&=E(x^{s-1})\hbox{ in  }g\cr
&={{\Gamma(1-{{1}\over{\alpha}}+{{s}\over{\alpha}})\Gamma({{1}\over{\alpha}}-{{s}\over{\alpha}})}
\over{\alpha\Gamma(1-s)}}&(19.3.5)\cr}
$$for $1-\alpha<\Re(s)<1, 0<\alpha\le 1$, obtained by putting
$\eta=1$ in (19.3.4) also. Since

$$\lim_{\alpha\rightarrow
1}{{\Gamma({{1}\over{\alpha}}-{{s}\over{\alpha}})}\over{\Gamma(1-s)}}=1
$$for $\alpha\rightarrow 1$, (19.3.4) reduces to

$$M_{g_1}(t)={{1}\over{\Gamma(\eta)}}\Gamma(\eta-1+s)\hbox{  for
}\alpha\rightarrow 1.\eqno(19.3.6)
$$Its inverse Mellin transform is then

$$g_1={{1}\over{\Gamma(\eta)}}{{1}\over{2\pi
i}}\int_{c-i\infty}^{c+i\infty}\Gamma(\eta-1+s)x^{-s}{\rm
d}s={{1}\over{\Gamma(\eta)}}x^{\eta-1}{\rm e}^{-x},x\ge
0,\eta>0\eqno(19.3.7)
$$which is the one-parameter gamma density and for $\eta=1$ it
reduces to the exponential density. Hence the generalized
Mittag-Leffler density $g_1$ can be taken as an extension of a gamma
density such as the one in (19.3.7) and the Mittag-Leffler density
$g$ as an extension of the exponential density for $\eta=1$. Is
there a structural representation for the random variable giving
rise to the Laplace transform in (19.2.4) corresponding to
(19.1.10)? The answer is in the affirmative and it is illustrated in
(19.1.14).

\vskip.3cm\noindent{\bf Note 19.3.1.}\hskip.3cm Pillai (1990,
Theorem 2.2), Lin (1998, Lemma 3) and others list the $\rho$-th
moment of the Mittag-Leffler density $g$ in (19.1.2) as follows:

$$\eqalignno{E(x^{\rho})&={{\Gamma(1-{{\rho}\over{\alpha}})\Gamma(1+{{\rho}\over{\alpha}})}
\over{\Gamma(1-\rho)}}, -\alpha<\Re(\rho)<\alpha<1.\cr
\noalign{\hbox{Therefore}}
E(x^{s-1})&={{\Gamma(1+{{1}\over{\alpha}}-{{s}\over{\alpha}})\Gamma(1-{{1}\over{\alpha}}+{{s}\over{\alpha}})}
\over{\Gamma(2-s)}}\cr
&={{({{1}\over{\alpha}}-{{s}\over{\alpha}})\Gamma({{1}\over{\alpha}}-{{s}\over{\alpha}})
\Gamma(1-{{1}\over{\alpha}}+{{s}\over{\alpha}})}\over{(1-s)\Gamma(1-s)}}
={{1}\over{\alpha}}{{\Gamma({{1}\over{\alpha}}-{{s}\over{\alpha}})
\Gamma(1-{{1}\over{\alpha}}+{{s}\over{\alpha}})}\over{\Gamma(1-s)}}\cr}
$$which is the expression in (19.3.5). Hence the two expressions are
one and the same.

\vskip.3cm\noindent{\bf Note 19.3.2.}\hskip.3cm If $y=ax,a>0$ and if
$x$ has a Mittag-Leffler distribution then the density of $y$ can
also be represented as a Mittag-Leffler function with the Laplace
transform

$$L_x(t)=(1+t^{\alpha})^{-1}\Rightarrow
L_y(t)=(1+(at)^{\alpha})^{-1},a>0,|(at)^{\alpha}|<1.
$$
\vskip.3cm\noindent{\bf Note 19.3.3.}\hskip.3cm From the
representation that

$$E(x^h)={{\Gamma(1-{{h}\over{\alpha}})\Gamma(1+{{h}\over{\alpha}})}\over{\Gamma(1-h)}},
-\alpha<\Re(h)<\alpha< 1
$$we have
$$E(x^0)={{\Gamma(1)\Gamma(1)}\over{\Gamma(1)}}=1.
$$Further, $g(x)$ in (19.2.1) is a non-negative function for all
$x$, with
$$E(x^h)=\int_0^{\infty}x^hg(x){\rm d}x=1\hbox{  for  }h=0.
$$Hence $g(x)$ is a density function for a positive random variable
$x$. Note that from the series form for the Mittag-Leffler function
it is not possible to show that $\int_0^{\infty}g(x){\rm d}x=1$
directly.

\vskip.3cm \noindent{\bf 19.4.\hskip.3cm Structural Representation
of the Generalized Mittag-Leffler Variable}

\vskip.3cm Let $u$ be the random variable corresponding to the
Laplace transform (19.2.7) with $t^{\alpha}$ replaced by $\delta
t^{\alpha}$ and $\gamma$ by $\eta$. Let $u$ be a positive L\'evy
variable with the Laplace transform ${\rm
e}^{-t^{\alpha}},0<\alpha\le 1$ and let $v$ be a gamma random
variable with parameters $\eta$ and $\delta$ or with he Laplace
transform $(1+\delta t)^{-\eta},\eta>0,\delta>0$. Let $u$ and $v$ be
statistically independently distributed.

\vskip.3cm \noindent{\bf Lemma 19.4.1.}\hskip.3cm{\it Let $u,v$ as
defined above. Then

$$w\sim uv^{{1}\over{\alpha}}\eqno(19.4.1)
$$where $w$ is a generalized Mittag-Leffler variable with Laplace transform
$(1+\delta t^{\alpha})^{-\beta},|\delta t^{\alpha}|<1$, where $\sim$
means `distributed as' or both sides have the same distribution.}

\vskip.3cm \noindent{\bf Proof.}\hskip.3cm Denoting the Laplace
transform of the density of $w$ by $L_w(t)$ and treating it as an
expected value

$$\eqalignno{L_w(t)&=E[{\rm e}^{-tv^{{1}\over{\alpha}}u}]=E[E[{\rm
e}^{-(tv^{{1}\over{\alpha}}u)}|v]]\cr &=E[{\rm
e}^{-(tv^{{1}\over{\alpha}})^{\alpha}}]=E[{\rm
e}^{-t^{\alpha}v}]=(1+\delta t^{\alpha})^{-\eta}&(19.4.2)\cr}
$$from the Laplace transform of a gamma variable. This establishes
the result.

\vskip.2cm From the structural representation in (19.4.1), taking
the Mellin transforms and writing as expected values, we have

$$E(w^{s-1})=E(u^{s-1})E(v^{{1}\over{\alpha}})^{s-1}\eqno(19.4.2)
$$due to statistical independence of $u$ and $v$. The left side is
available from (19.3.4) as

$$E(w^{s-1})={{\Gamma(\eta-{{1}\over{\alpha}}+{{s}\over{\alpha}})
\Gamma({{1}\over{\alpha}}-{{s}\over{\alpha}})\delta^{{s-1}\over{\alpha}}}
\over{\alpha\Gamma(\eta)\Gamma(1-s)}}.\eqno(19.4.3)
$$Let us compute $E[v^{{1}\over{\alpha}}]^{s-1}$ from the gamma
density. That is,

$$E[v^{{1}\over{\alpha}}]^{s-1}={{1}\over{\delta^{\eta}\Gamma(\eta)}}\int_0^{\infty}(v^{{1}\over{\alpha}})^{s-1}v^{\eta-1}{\rm
e}^{-{{v}\over{\delta}}}{\rm
d}v={{\Gamma(\eta-{{1}\over{\alpha}}+{{s}\over{\alpha}})\delta^{{s-1}\over{\alpha}}}\over{\Gamma(\eta)}}\eqno(19.4.4)
$$for $\Re(s)>1-\alpha\eta, 0<\alpha\le 1, \eta>0$. Comparing
(19.4.3) and (19.4.4) we have the $(s-1)$-th moment of a L\'evy
variable

$$E[u^{s-1}]={{\Gamma({{1}\over{\alpha}}-{{s}\over{\alpha}})}\over{\alpha\Gamma(1-s)}}
={{\Gamma(1+{{1}\over{\alpha}}-{{s}\over{\alpha}})}\over{\Gamma(2-s)}},
\Re(s)<1,0<\alpha\le 1.\eqno(19.4.5)
$$

\vskip.3cm \noindent{\bf Note 19.4.1.}\hskip.3cm Lin (1998) gives
the $\rho$-th moment of a L\'evy variable with parameter $\alpha$ as

$$\eqalignno{E[u^{\rho}]&={{\Gamma(1-{{\rho}\over{\alpha}})}\over{\Gamma(1-\rho)}}.&(19.4.6)\cr
\noalign{\hbox{Hence for $\rho=s-1$ we have}}
E[u^{s-1}]&={{\Gamma(1+{{1}\over{\alpha}}-{{s}\over{\alpha}})}\over{\Gamma(2-s)}}
={{({{1}\over{\alpha}}-{{s}\over{\alpha}})\Gamma({{1}\over{\alpha}}-{{s}\over{\alpha}})}
\over{(1-s)\Gamma(1-s)}}\cr
&={{\Gamma({{1}\over{\alpha}}-{{s}\over{\alpha}})}\over{\alpha\Gamma(1-s)}}\hbox{
for  }\Re(s)<1.\cr}
$$This is (19.4.5) and hence both the representations are one and
the same.

\vskip.2cm Hence the L\'evy density, denoted by $g_2(u)$, can be
written as

$$g_2(u)={{1}\over{2\pi
i}}\int_{c-i\infty}^{c+i\infty}{{\Gamma({{1}\over{\alpha}}-{{s}\over{\alpha}})}
\over{\alpha\Gamma(1-s)}}u^{-s}{\rm
d}s.\eqno(19.4.7)
$$Its Laplace transform is then

$$\eqalignno{L_{g_2}(t)&={{1}\over{2\pi
i}}\int_{c-i\infty}^{c+i\infty}{{\Gamma({{1}\over{\alpha}}-{{s}\over{\alpha}})}
\over{\alpha\Gamma(1-s)}}[\int_0^{\infty}u^{1-s-1}{\rm e}^{-tu}{\rm
d}u]{\rm d}s\cr &={{1}\over{2\pi
i}}\int_{c-i\infty}^{c+i\infty}{{1}\over{\alpha}}\Gamma({{1}\over{\alpha}}-{{s}\over{\alpha}})t^{-1+s}{\rm
d}s\cr &={{1}\over{2\pi
i}}\int_{c_1-i\infty}^{c_1+i\infty}{{1}\over{\alpha}}\Gamma({{s}\over{\alpha}})t^{-s}{\rm
d}s&(19.4.7)\cr}
$$by making the substitution $-1+s=-s_1$. Then evaluating as the sum
of the residues at ${{s}\over{\alpha}}=-\nu, \nu=0,1,2,...$ we have

$$L_{g_2}(t)=\sum_{\nu=0}^{\infty}{{(-1)^{\nu}}\over{\nu!}}t^{\alpha\nu}={\rm
e}^{-t^{\alpha}}.\eqno(19.4.8)
$$This verifies the result about the Laplace transform of the
positive L\'evy variable with parameter $\alpha$. Note that when
$\alpha=1$, (19.4.8) gives the Laplace transform of a degenerate
random variable taking the value $1$ with probability $1$.

\vskip.2cm Mellin convolution of certain L\'evy variables can be
seen to be again a L\'evy variable.

\vskip.3cm \noindent {\bf Lemma 19.4.2}.\hskip.3cm{\it Let $x_j$ be
a positive L\'evy variable with parameter $\alpha_j, 0<\alpha_j\le
1$ and let $x_1,...,x_p$ be statistically independently distributed.
Then

$$u=x_1x_2^{{1}\over{\alpha_1}}...x_p^{{1}\over{\alpha_1\alpha_2...\alpha_{p-1}}}
$$is distributed as a L\'evy variable with parameter
$\alpha_1\alpha_2...\alpha_p$.}

\vskip.3cm \noindent{\bf Proof.}\hskip.3cm From (19.4.5)

$$\eqalignno{E[{\rm e}^{-tx_j}]&={\rm
e}^{-t^{\alpha_j}},0<\alpha_j\le 1, j=1,...,p\cr E[{\rm
e}^{-tu}]&=E[{\rm
e}^{-tx_1x_2^{{1}\over{\alpha_1}}...x_p^{{1}\over{\alpha_1...\alpha_{p-1}}}}]=E[E[{\rm
e}^{-tx_1...x_p^{{1}\over{\alpha_1...\alpha_{p-1}}}}|_{x_2,...,x_p}]]\cr
&=E[{\rm
e}^{-t^{\alpha_1}x_2x_3^{{1}\over{\alpha_2}}...x_p^{{1}\over{\alpha_2...\alpha_{p-1}}}}].\cr}
$$Repeated application of the conditional argument gives the final
result as

$$E[{\rm e}^{-tu}]={\rm e}^{-t^{\alpha}},
\alpha=\alpha_1\alpha_2...\alpha_p, 0<\alpha_1...\alpha_p\le 1
$$which means that $u$ is distributed as a L\'evy with parameter
$\alpha_1...\alpha_p$.

\vskip.2cm From the representation in (19.4.1) we can compute the
moments of the natural logarithms of Mittag-Leffler, L\'evy and
gamma variables.

$$w=uv^{{1}\over{\alpha}}\Rightarrow \ln w=\ln
u+{{1}\over{\alpha}}\ln v.\eqno(19.4.9)
$$But from (19.4.3) and (19.4.6) we have the $h$-th moments of $u$
and $v$ given by

$$\eqalignno{E[u^h]&={{\Gamma(1-{{h}\over{\alpha}})}\over{\Gamma(1-h)}},
\Re(h)<\alpha\le 1&(19.4.10)\cr \noalign{\hbox{and}}
E[v^{{1}\over{\alpha}}]^h&={{\Gamma(\eta+{{h}\over{\alpha}})\delta^{{h}\over{\alpha}}}\over{\Gamma(\eta)}},
\Re(\eta+{{h}\over{\alpha}})>0.&(19.4.11)\cr \noalign{\hbox{From
(19.4.3)}}
E[w^h]&={{\Gamma(\eta+{{h}\over{\alpha}})\Gamma(1-{{h}\over{\alpha}})\delta^{{h}\over{\alpha}}}\over{\Gamma(\eta)\Gamma(1-h)}}.&(19.4.12)\cr}
$$But for a positive random variable $z$

$$\eqalignno{E[z^h]&=E[{\rm e}^{\ln z^h}].\cr
\noalign{\hbox{Hence}} {{{\rm d}}\over{{\rm
d}h}}E[z^h]|_{h=0}&=E[\ln z {\rm e}^{h\ln z}]_{h=0}=E[\ln z].\cr}
$$Therefore from (19.4.9) to (19.4.12) we have the following:

$$\eqalignno{E[\ln w]&={{{\rm d}}\over{{\rm
d}h}}\left\{{{\delta^{{h}\over{\alpha}}\Gamma(\eta+{{h}\over{\alpha}})\Gamma(1-{{h}\over{\alpha}})}
\over{\Gamma(\eta)\Gamma(1-h)}}\right\}|_{h=0}\cr
&={{1}\over{\alpha}}\psi((\eta)-{{1}\over{\alpha}}\psi(1)+\psi(1)+{{1}\over{\alpha}}\ln
\delta &(19.4.13)\cr}
$$by taking the logarithmic derivative, where $\psi$ is a psi
function, see, for example Mathai (1993).

$$\eqalignno{E[\ln v^{{1}\over{\alpha}}]&={{{\rm d}}\over{{\rm
d}h}}\left\{{{\delta^{{h}\over{\alpha}}\Gamma(\eta+{{h}\over{\alpha}})}
\over{\Gamma(\eta)}}\right\}|_{h=0}={{1}\over{\alpha}}\psi(\eta)+{{1}\over{\alpha}}\ln
\delta&(19.4.14)\cr \noalign{\hbox{or}} E[\ln v]&=\psi(\eta)+\ln
\delta\cr \noalign{\hbox{and}} E[\ln u]&={{{\rm d}}\over{{\rm
d}h}}\left\{{{\Gamma(1-{{h}\over{\alpha}})}\over{\Gamma(1-h)}}\right\}|_{h=0}
=-{{1}\over{\alpha}}\psi(1)+\psi(1)&(14.4.15)\cr}
$$where $\psi(1)=-\gamma$ where $\gamma$ is the Euler's constant.

\vskip.3cm \noindent{\bf Note 19.4.2.}\hskip.3cm The relations on
the expected values of the logarithms of Mittag-Leffler variable,
positive L\'evy variable and exponential variable, given on page
1432 of Jayakumar (2003), where $\eta=1$, are not correct.

\vskip.3cm \noindent{\bf 19.5.\hskip.3cm A Pathway from
Mittag-Leffler Distribution to Positive L\'evy Distribution}

\vskip.3cm Consider the function

$$\eqalignno{f(x)&=\sum_{k=0}^{\infty}{{(-1)^k(\eta)_k}\over{k!\Gamma(\alpha\eta+\alpha
k)}}{{x^{\alpha\eta-1+\alpha
k}}\over{(a^{{1}\over{\alpha}})^{\alpha\eta+\alpha
k}}},\eta>0,a>0,0<\alpha\le 1\cr
&={{x^{\alpha\eta-1}}\over{a^{\eta}}}E_{\alpha,\alpha\eta}^{\eta}(-{{x^{\alpha}}\over{a}}).\cr}
$$Thus $x=a^{{1}\over{\alpha}}y$ where $y$ is a generalized
Mittag-Leffler variable. The Laplace transform of $f$ is given by
the following:

$$\eqalignno{L_f(t)&=\sum_{k=0}^{\infty}{{(-1)^k(\eta)_k}\over{k!a^{\eta+k}}}\int_0^{\infty}{{x^{\alpha\eta+\alpha
k-1}{\rm e}^{-tx}}\over{\Gamma(\alpha\eta+\alpha k)}}{\rm d}x\cr
&=\sum_{k=0}^{\infty}{{(-1)^k(\eta)_k}\over{k!a^{\eta+k}}}t^{-\alpha\eta-\alpha
k}=[1+at^{\alpha}]^{-\eta}, |at^{\alpha}|<1.&(19.5.1)\cr}
$$If $\eta$ is replaced by ${{\eta}\over{q-1}}$ and $a$ by $a(q-1)$
with $q>1$ then we have a Laplace transform

$$\eqalignno{L_f(t)&=[1+a(q-1)t^{\alpha}]^{-{{\eta}\over{q-1}}},
q>1&(19.5.2)\cr \noalign{\hbox{If $q\rightarrow 1_{+}$ then}}
L_f(t)&\rightarrow {\rm e}^{-a\eta
t^{\alpha}}=L_{f_1}(t)&(19.5.3)\cr}
$$which is the Laplace transform of a constant multiple of a positive
L\'evy variable with parameter $\alpha$. Thus $q$ here creates a
pathway of going from the general Mittag-Leffler density $f$ to a
positive L\'evy density $f_1$ with parameter $\alpha$, the
multiplying constant being $(a\eta)^{{1}\over{\alpha}}$. For a
discussion of a general rectangular matrix-variate pathway model see
Mathai (2005). The result in (19.5.3) can be put in a more general
setting. Consider an arbitrary real random variable $y$ with the
Laplace transform, denoted by $L_y(t)$, and given by

$$L_y(t)={\rm e}^{-\phi(t)}\eqno(19.5.4)
$$where $\phi(t)$ is a function such that
$\phi(tx^{\gamma})=x\phi(t), \phi(t)\ge 0, \lim_{t\rightarrow
0}\phi(t)=0$ for some real positive $\gamma$. Let

$$u=yx^{\gamma}\eqno(19.5.5)
$$where $x$ and $y$ are independently distributed with $y$ having
the Laplace transform in (19.5.4) and $x$ having a two-parameter
gamma density with shape parameter $\beta$ and scale parameter
$\delta$ or with the Laplace transform

$$L_x(t)=(1+\delta t)^{-\beta}.\eqno(19.5.6)
$$Now consider the Laplace transform of $u$ in (19.5.5), denoted by
$L_u(t)$. Then

$$\eqalignno{L_u(t)&=E[{\rm e}^{-tu}]=E[{\rm
e}^{-tyx^{\gamma}}]=E[E[{\rm e}^{tyx^{\gamma}}|x]]\cr &=E[{\rm
e}^{-\phi(tx^{\gamma})}]=E[{\rm e}^{-x\phi(t)}]\cr
\noalign{\hbox{from the assumed property of $\phi(t)$}}
&=[1+\delta\phi(t)]^{-\beta}.&(19.5.7)\cr}
$$If $\delta$ is replaced by $\delta(q-1)$ and $\beta$ by
${{\beta}\over{q-1}}$ with $q>1$ then we get a path through $q$.
That is, when $q\rightarrow 1_{+}$,

$$L_u(t)=[1+\delta(q-1)\phi(t)]^{-{{\beta}\over{q-1}}}\rightarrow
{\rm e}^{-\delta\beta\phi(t)}={\rm
e}^{-\phi(\delta^{\gamma}\beta^{\gamma}t)}.\eqno(19.5.8)
$$If $\phi(t)=t^{\alpha},0<\alpha\le 1$ then

$$L_u(t)={\rm e}^{-(\delta\beta)^{\gamma\alpha}t^{\alpha}}
$$which means that $u$ goes to a constant multiple of a positive
L\'evy variable with parameter $\alpha$, the constant being
$(\delta\beta)^{\gamma}$.

\vskip.3cm \noindent{\bf 19.6.\hskip.3cm Linnik or $\alpha$-Laplace
Distribution}

\vskip.3cm A Linnik random variable is defined as that real scalar
random variable whose characteristic function is given by

$$\phi(t)={{1}\over{1+|t|^{\alpha}}},0<\alpha\le
2,-\infty<t<\infty.\eqno(19.6.1)
$$For $\alpha=2$, (19.6.1) corresponds to the characteristic function
of a Laplace random variable and hence Pillai (1995) called the
distribution corresponding to (19.6.1) as the $\alpha$-Laplace
distribution. For positive variable, (19.6.1) reduces to the
characteristic function of a Mittag-Leffler variable. Infinite
divisibility, characterizations, other properties and related
materials may be seen from the review paper Jayakumar and Suresh
(2003) and the many references therein, Pakes (1998) and Mainardi
and Pagnini (2008). Multivariate generalization of Mittag-Leffler
and Linnik distributions may be seen from Lim and Teo (2009). Since
the steps for deriving results on Linnik distribution are parallel
to those of the Mittag-Leffler variable, further discussion of
Linnik distribution is omitted.

\vskip.3cm \noindent{\bf 19.7.\hskip.3cm Multivariable
Generalization of Mittag-Leffler, Linnik and L\'evy Distributions}

\vskip.3cm A multivariate Linnik distribution can be defined in
terms of a multivariate L\'evy vector. Let $T'=(t_1,...,t_p),
X'=(x_1,...,x_p)$, prime denoting the transpose. A vector variable
having positive L\'evy distribution is given by the characteristic
function

$$E[{\rm e}^{iT'X}]={\rm e}^{-(T'\Sigma T)^{{\alpha}\over2}},
0<\alpha\le 2,\eqno(19.7.1)
$$where $\Sigma=\Sigma'>0$ is a real positive definite $p\times p$
matrix. Consider the representation

$$u=y^{{1}\over{\alpha}}X\eqno(19.7.2)
$$where the $p\times 1$ vector $X$, having a multivariable L\'evy
distribution with parameter $\alpha$, and $y$ a real scalar gamma
random variable with the parameters $\delta$ and $\beta$, are
independently distributed. Then the characteristic function of the
random vector variable $u$ is given by the following:

$$\eqalignno{E[{\rm e}^{iy^{{1}\over{\alpha}}T'X}]&=E[E[{\rm
e}^{iy^{{1}\over{\alpha}}T'X}]|_{y}]\cr &=E[{\rm e}^{-y[|T'\Sigma
T|]^{{\alpha}\over2}}]=[1+\delta |T'\Sigma
T|^{{\alpha}\over2}]^{-\beta}.&(19.7.3)\cr}
$$Then the distribution of $u$, with the characteristic function in
(19.7.3) is called a vector-variable Linnik distribution. Some
properties of this distribution are given in Lim and Teo (2009).

\vskip.5cm \noindent{\bf 20.\hskip.3cm Mittag-Leffler Stochastic
Processes}

\vskip.3cm The stochastic process $\{x(t), t>0\}$ having stationary
independent increment with $x(0)=0$ and $x(1)$ having the Laplace
transform
$$L_{x(1)}(\lambda)=(1+\lambda^{\alpha})^{-1}, 0<\alpha\le 1,
\lambda>0,
$$which is the Laplace transform of a Mittag-Leffler random
variable, is called the Mittag-Leffler stochastic process. Then the
Laplace transform of $x(t)$, denoted by $L_{x(t)}(\lambda)$, is
given by

$$L_{x(t)}(\lambda)=[(1+\lambda^{\alpha})^{-1}]^t=[1+\lambda^{\alpha}]^{-t}.\eqno(20.1.1)
$$The density corresponding to the Laplace transform (20.1.1) or the
density of $x(t)$ is then available as the following:

$$\eqalignno{f_{x(t)}(x)&=\sum_{k=0}^{\infty}(-1)^k{{(t)_k}\over{k!}}{{x^{\alpha
k+\alpha t-1}}\over{\Gamma(\alpha k+\alpha t)}}\cr &=x^{\alpha
t-1}E_{\alpha,\alpha t}^{t}(-x^{\alpha}),0<\alpha\le 1, x\ge 0,
t>0.&(20.1.2)\cr}
$$The distribution function of $x(t)$ is given by

$$\eqalignno{F_{x(t)}(x)&=\int_0^xf_{x(t)}(y){\rm
d}y=\sum_{k=0}^{\infty}(-1)^k{{(t)_k}\over{k!}}{{x^{\alpha k+\alpha
t}}\over{(\alpha k+\alpha t)\Gamma(\alpha k+\alpha t)}}\cr
&=\sum_{k=0}^{\infty}{{(-1)^k}\over{k!}}{{\Gamma(t+k)}\over{\Gamma(t)}}{{x^{\alpha
k+\alpha t}}\over{\Gamma(1+\alpha k+\alpha t)}}, 0<\alpha\le 1,
t>0.&(20.1.3)\cr}
$$This form is given by Pillai (1990) and by his students.

\vskip.3cm\noindent{\bf 20.2.\hskip.3cm Linear First Order
Autoregressive processes}

\vskip.3cm Consider the stochastic process

$$x_n=\cases{e_n,\hbox{  with probability  }p, 0\le p\le 1\cr
e_n+ax_{n-1}\hbox{  with probability  }1-p, 0<a\le
1.\cr}\eqno(20.2.1)
$$Let the sequence $\{e_n\}$ be independently and identically
distributed with Laplace transform $L_{e}(\lambda)$ and let
$\{x_n\}$ be identically distributed with Laplace transform
$L_{x}(\lambda)$. From the representation in (20.2.1)

$$\eqalignno{L_{x_n}(\lambda)&=pL_{e}(\lambda)+(1-p)L_{e}(\lambda)L_{x_{n-1}}(a\lambda)\cr
\noalign{\hbox{Therefore}}
L_{e}(\lambda)&={{L_{x_n}(\lambda)}\over{p+(1-p)L_{x_{n-1}}(a\lambda)}}
={{L_{x}(\lambda)}\over{p+(1-p)L_{x}(a\lambda)}}&(20.2.1)\cr}
$$assuming stationarity. When $p=0$,

$$L_{e}(\lambda)={{L_{x}(\lambda)}\over{L_{x}(a\lambda)}}, 0<a\le
1\eqno(20.2.3)
$$which defines class $L$ distributions, for all $a, 0<a<1$. When
$p=0$, (20.2.3) implies that the innovation sequence $\{e_n\}$
belongs to class $L$ distributions. Then (20.2.3) can lead to two
autoregressive situations, the first order exponential
autoregressive process $EAR(1)$ and the first order Mittag-Leffler
autoregressive process $MLAR(1)$.

\vskip.3cm \noindent {Concluding remarks} \vskip.3cm The various
Mittag-Leffler functions discussed in this paper will be useful for
investigators in various disciplines of applied sciences and
engineering. The importance of Mittag-Leffler function in physics is
steadily increasing. It is simply said that deviations of physical
phenomena from exponential behavior could be governed by physical
laws through Mittag-Leffler functions (power-law). Currently more
and more such phenomena are discovered and studied. \vskip.2cm It is
particularly important for the disciplines of stochastic systems,
dynamical systems theory, and disordered systems. Eventually, it is
believed that all these new research results will lead to the
discovery of truly nonequilibrium statistical mechanics. This is
statistical mechanics beyond Boltzmann and Gibbs. This
nonequilibrium statistical mechanics will focus on entropy
production, reaction, diffusion, reaction-diffusion, etc and may be
governed by fractional calculus. \vskip.2cm Right now, fractional
calculus and H-function (Mittag-Leffler function) are very important
in research in physics.

\vskip.3cm \noindent \centerline{\bf Acknowledgment}

\vskip.3cm \noindent The authors would like to thank the Department
of Science and Technology, Government of India, New Delhi, for the
financial assistance under Project No. SR/S4/MS:287/05 which made
this research collaboration possible.\vskip.2cm The authors are grateful 
to Constantino Tsallis, Enrico Scalas, Francesco Mainardi, Henkel Malte, 
Zivorad Tomovski, and the Wolfram Demonstration Team for suggestions 
to improve the paper.

\vskip.5cm \noindent \centerline{\bf References}

\vskip.5cm \noindent [1] M. Abramowitz and I.A. Stegun, {\it
Handbook of Mathematical Functions}, Dover, New York, 1965.

\vskip.2cm \noindent [2] M.A. Al-Bassam and Yu F. Luchko, On
generalized fractional calculus and its application to the solution
of integro-differential equations, Journal of Fractional Calculus 7(1995), 69-88.

\vskip.2cm \noindent [3] S.N. Agal and C.L. Koul, Weyl fractional
calculus and Laplace transform, Proceedings of the Indian Academy of Sciences (Mathematical
Sciences), 92(1983), 167-170.

\vskip.2cm \noindent [4] R.P. Agarwai, A propos d'une note de M.
Pierre Humbert, C.R. Acad. Sci. Paris, 236(1953), 2031-2032.

\vskip.2cm \noindent [5] W.A. Al-Salam, Some fractional q-integrals
and q-derivatives, Proceedings of the Edinburgh Mathematical Society, 15(1967), 135-140.

\vskip.2cm \noindent [6] B. N. Al-Saqabi, Solution of a class of
differintegral equations by means of Riemann-Liouville operator, Journal of
Fractional Calculus 8(1995), 95-102.

\vskip.2cm \noindent [7] B.N. Al-Sqabi, S.L. Kalla and H.M.
Srivastava, A certain family of infinite series associated with
digamma functions, Journal of Mathametical Analysis and its Applications, 159(1991), 361-372.

\vskip.2cm \noindent [8] B.N. Al-Saqabi, and V.K. Tuan, Solution of
a fractional differintegral equation. Integral Transforms and Special
Functions, 4(2006), 321-326.

\vskip.2cm \noindent [9] V.V. Anh and N.N. Leonenko, Spectral
analysis of fractional kinetic equations with random data, Journal of
Statisticl Physics, 104(2001), 1349-1387.

\vskip.2cm \noindent [10] Yu. I. Babenko, {\it Heat and Mass Transfer},
Leningrad, Khimiya, 1986 (in Russian).

\vskip.2cm\noindent [11] R.L. Bagley, {\it Applications of
Generalized Derivatives to Viscoelasticity}, Ph.D. thesis, Air
Force Institute of Technology, 1979.

\vskip.2cm \noindent [12] R.L. Bagley, On the equivalence of the
Riemann-Liouville and the Caputo fractional order derivatives in
modeling of linear viscoelastic materials, Fractional Calculus and Applied Analysis,
10(2007), 123-126.

\vskip.2cm \noindent [13] R.L. Bagley and P.J. Torvik, On the
appearance of the fractional derivative in the behaviour of real
materials, Journal of Applied Mechanics, 51(1984), 294-298.

\vskip.2cm \noindent [14] J.H. Barret, Differential equations of
non-integer order, Canadian Journal of Mathematics, 6(1954), 529-541.

\vskip.2cm \noindent [15] L. Beghin and E. Orsingher, The
distribution of the local time for ``pseudoprocesses'' and its
connection with fractional diffusion equations, Stochastic Processes
and their Applications, 115(2005), 1017-1040.

\vskip.2cm \noindent [16] L. Beghin and E. Orsingher, The telegraph
process stopped at stable distributed times and its connection with
the fractional telegraph equation, Fractional Calculus and Applied Analysis,
6(2003), 187-204.

\vskip.2cm \noindent [17] L. Beghin and E. Orsingher, Iterated
elastic Brownian motions and fractional diffusion equations,
Stochastic Processes and Their Applications, 119(2009), 1975-2003.

\vskip.2cm \noindent [18] M.N. Berberan-Santos, Relation between
the inverse Laplace transforms $I(t^{\beta})$ and $I(t)$:
Application to the Mittag-Leffler and asymptotic inverse power law
relaxation functions, Journal of Mathematical Chemistry, 38(2005),
265-270.

\vskip.2cm \noindent [19] M.N. Berberan-Santos, Analytic inversion
of the Laplace transform without contour integration: Application to
luminescence decay laws and other relaxation functions, Journal of
Mathematical Chemistry, 38(2)(2005a), 165-173.

\vskip.2cm \noindent [20] M.N. Berberan-Santos, Properties of the
Mittag-Leffler relaxation function, Journal of Mathematical Chemistry,
38(4)(2005b), 629-635.

\vskip.2cm \noindent [21] G.W.S. Blair, Psychorheology: Links
between the past and the present, Journal of Texture Studies,
5(1974), 3-12.

\vskip.2cm \noindent [22] G.W.S. Blair, The role of psychophysics
in rheology, Journal of Colloid Sciences, (1947), 21-32.

\vskip.2cm \noindent [23] L. Bondesson, {\it Generalized Gamma
Convolutions and Related Classes of Distributions and Densities},
Lecture Notes in Statistics, 76 (1992), Springer-Verlag, New York.

\vskip.2cm \noindent [24] B.L.J. Braaksma, Asymptotic expansions
and analytic continuations for a class of Barnes-integrals,
Compositio Mathematica, 15(1964), 239-341.

\vskip.2cm \noindent [25] E. Buckwar and Yu. Luchko, Invariance of
a partial differential equation of fractional order under the Lie
group of scaling transformations, Journal of Mathematical Analysis and its Applications, 227(1998),
81-97.

\vskip.2cm \noindent [26] A. Buhl, Series analytiques
sommabilit\'e, Mem. des Sci. Acad. Sci., Paris Fasc. 7,
Gauthier-Villars, Paris, 1925, Ch.3.

\vskip.2cm \noindent [27] R.G. Buschman and H.M. Srivastava, The
$\bar{H}$-function associated with a certain Feynman integrals, Journal of
Physics A: Mathematical and General 23(1990), 4707-4710.

\vskip.2cm\noindent [28] M. Caputo, {\it Elasticit\'a e Dissipazione},
Zanichelli, Bologna, 1969 (in Italian).

\vskip.2cm \noindent [29] M. Caputo and F. Mainardi, Linear models
of dissipation in anelastic solids, Rivista del Nuovo Cimento, Ser.
II, 1(1971), 161-198.

\vskip.2cm \noindent [30] B.C. Carlson, Appell functions and
multiple averages, SIAM Journal of Mathematical Analysis, 2(1971), 420-430.

\vskip.2cm \noindent [31] B.C. Carlson, {\it Special Functions of
Applied Mathematics}, Academic Press, New York, 1977.

\vskip.2cm \noindent [32] H. Chamati and N.S. Tonchev, Generalized
Mittag-Leffler functions in the theory of finite-size scaling for
systems with strong anisotropy and/or long-range interaction, Journal of
Physics A: Mathematical and General, 39(2006), 469-478.

\vskip.2cm \noindent [33] A.V. Chechkin, R. Gorenflo, I.M. Sokolov
and V. Yu. Gochar, Distributed order time fractional diffusion
equation, Fractional Calculus and Applied Analysis, 6(3)(2003), 259-279.

\vskip.2cm\noindent [34] A. Compte, Stochastic foundation of
fractional dynamics, Physical Review E, 53(1996), 4191-4193.

\vskip.2cm \noindent [35] H.T. Davis, The applications of
fractional operators to functional equations, American Journal of Mathematics,
49(1)(1927), 123-142.

\vskip.2cm \noindent [36] H.T. Davis, {\it The Theory of Linear
Operators}, Principia Press, Bloomington, Indiana, 1936.

\vskip.2cm \noindent [37] L. Debnath, {\it Integral Transforms and Their
Applications} (2nd edition), Boca Raton, CRC Press, Chapman \& Hall,
2005.

\vskip.2cm \noindent [38] L. Debnath, Applications of fractional
integral and differential equations in fluid mechanics, Fractional Calculus
and Applied Analysis, 6(2003), 119-155.

\vskip.2cm \noindent [39] L. Debnath and D.D. Bhatta, Solution to
few linear fractional inhomogeneous partial differential equations
in fluid mechanics, Fractional Calculus and Applied Analysis, 7(2004), 21-36.

\vskip.2cm \noindent [40] Y. Deora, P.K. Banerji and M. Saigo,
Fractional integral and Dirichlet averages, Journal of Fractional Calculus,
6(1994), 55-60.

\vskip.2cm \noindent [41] L. Devroye, A note on Linnik
distribution, Statistics and Probability Letters, 5(1990), 305-306.

\vskip.2cm\noindent [42] M.M. Dzherbashyan, On the integral
transformations  generated by the generalzied Mittag-Leffler
fuction, (in Russian), Izv. Akad. Nauk Arm. SSR, 13(3)(1960), 21-63.

\vskip.2cm \noindent [43] M.M. Dzherbashyan, About integral
representation of functions continuous on some rays (generalization
of Fourier integral), Izv. Akad. Nauk Arm. SSR, Ser. Mat.,
18(1964), 427-448 (in Russian).

\vskip.2cm \noindent [44] M.M. Dzherbashyan, {\it Integral Transforms
and Representations of Functions in the Complex Plane}, Nauka,
Moscow, 1966 (in Russian).

\vskip.2cm \noindent [45]  M.M. Dzherbashyan, On the integral
representation and uniqueness of some classes of entire functions
(in Russian), Dokl. Akad. Nauk SSR, 85(1)(1952), 29-32.

\vskip.2cm \noindent [46] A.M. A. El-Sayed, Fractional order
diffusion-wave equation, Journal of Theoretical Physics, 35(2)(1996),
371-382.

\vskip.2cm \noindent [47] A. Erd\'elyi, W. Magnus, F. Oberhettinger
and F.G. Tricomi, {\it Higher Transcendental Functions}, Vol. 1, New York,
McGraw-Hill, 1953.

\vskip.2cm \noindent [48] A. Erd\'elyi, W. Magnus, F. Oberhettinger
and F.G. Tricomi, {\it Higher Transcendental Functions}, Vol. 3, New York,
McGraw-Hill, 1955.

\vskip.2cm \noindent [49] A. Erd\'elyi, W. Magnus, F. Oberhettinger
and F.G. Tricomi, {\it Tables of Integral Transforms}, Vol. 1, New York,
McGraw-Hill, 1954a.

\vskip.2cm \noindent [50] A. Erd\'elyi, W. Magnus, F. Oberhettinger
and F.G. Tricomi, {\it Tables of Integral Transforms}, Vol. 2, New York,
McGraw-Hill, 1954b.

\vskip.2cm\noindent [51] S. Ferraro, M. Manzim, A. Masoero and E.
Scalas, A random telegraph signal of Mittag-Leffler type, Physica A,
388(2009), 3991-3999.

\vskip.2cm \noindent [52] C. Fox, The G and H-functions as
symmetrical Fourier kernels, Transactions of the American Mathematical Society, 98(1961),
395-429.

\vskip.2cm \noindent [53] Y. Fujita, A characterization of the
results of Pillai, Annals of the Institute of Statistical Mathematics, 45(2)(1993), 361-365.

\vskip.2cm \noindent [54] Y. Fujita, Integro-differential equation
which interpolates the beat equation and the wave equation (I),
Osaka Journal of Mathematics, 27(2)(1990), 309-321.

\vskip.2cm \noindent [55] Y. Fujita, Integro-differential equations
which interpolates the beat equation and wave equation (II), Osaka
Journal of Mathematics, 27(4)(1990a), 797-804.

\vskip.2cm \noindent [56] Y. Fujita, Cauchy problems of fractional
order and stable processes, Japan Journal of Applied Mathematics, 7(3)(1990b),
459-470.

\vskip.2cm \noindent [57] Y. Fujita, Energy inequalities for
integro-partial differential equations with Riemann-Liouville
integrals, SIAM Journal of Mathematical Analysis, 23(5)(1992),  1182-1188.

\vskip.2cm \noindent [58] M. Fujiwara, On the integration and
differentiation of an arbitrary order, Tohoku Mathematical Journal, 37(1993),
110-121.

\vskip.2cm \noindent [59] M. Fujiwara, On uniqueness of the
solution of the initial value problems for ordinary fractional
differential equations, International Journal of Applied Mathematics, 2(10)(2002), 177-189.

\vskip.2cm \noindent [60] M. Fujiwara, On uniqueness of the
solution of initial value problems for ordinary fractional
differential equations II, International Journal of Applied Mathematics, 3(10)(2002a),
269-283.

\vskip.2cm \noindent [61] D. Fulger, E. Scalas and G. Germano,
Monte-Carlo  simulation of uncoupled continuous-time random walks
yielding a stochastic solution of the space-time fractional
diffusion equation, Physical Review E 77(2008), 021122.

\vskip.2cm \noindent [62] L. Galu\'e, Composition of hypergeometric
fractional operators, Kuwait Journal of Science and Engineering, 27(2000),
1-14.

\vskip.2cm \noindent [63] L. Galu\'e, S.L. Kalla and Y.K. Tuan,
Composition of Erd\'elyi-Kober fractional operators, Integral
Transforms and Special Functions, 9(2000), 185-196.

\vskip.2cm \noindent [64] L. Galu\'e, S.L. Kalla and H.M.
Srivastava, Further results on an H-function generalized fractional
calculus, Journal of Fractional Calculus, 4(1993), 89-102.

\vskip.2cm \noindent [65] L. Galu\'e, V.S. Kiryakova and S.L.
Kalla, Solution of dual integral equations by fractional calculus,
Mathematica Balkanica, 7(1993), 53-72.

\vskip.2cm \noindent [66] M. Garg, A. Rao and S.L.Kalla, On
Mittag-Leffler type functions and fractional calculus operators,
Mathematica Balkanica, New Series, 21(2007), Fasc. 3-4, 349-360.

\vskip.2cm \noindent [67] R.S. Garg, On multidimensional Mellin
convolutions and H-function transformations, Indian Journal of Pure and Applied
Mathematics, 13(1982), 30-38.

\vskip.2cm \noindent [68] G. Gasper and M. Rahman, {\it Basic
Hypergeometric Series: Encyclopedia of Mathematics and Its
Applications}, Vol. 35, Cambridge University Press, Cambridge, 1990.

\vskip.2cm\noindent [69] I.M. Gel'fand and G.E. Shilov, {\it
Generalized Functions, Vol. 1}, Academic Press, New York, 1964.

\vskip.2cm\noindent [70] G. Germano, M. Politi, E. Scalas and R.L.
Schilling, Stochastic calculus for uncoupled continuous-time random
walks, Physical Reviews E, 79(2009), 066102.

\vskip.2cm \noindent [71] W.G. Gl\"ockle and T.F. Nonnenmacher,
Fractional integral operators and Fox functions in the theory of
viscoelasticity, Macromolecules, American Chemical Society,
24(24)(1991), 6426-6434.

\vskip.2cm \noindent [72] W.G. Gl\"ockle and T.F. Nonnenmacher, Fox
function representation of non-Debye relaxation processes, Journal
of Statistical Physics, 71(1993), 741-747.

\vskip.2cm \noindent [73] R. Gorenflo, A. Iskenderov and Y. Luchko,
Mapping between solutions of fractional diffusion wave equations,
Fractional Calculus and Applied Analysis, 3(2000), 75-86.

\vskip.2cm \noindent [74] R. Gorenflo, A.A. Kilbas and S.V.
Rogosin, On the generalized Mittag-Leffler type function, Integral
Transforms and Special Functions, 7(3-4)(1998), 215-224.

\vskip.2cm \noindent [75] R. Gorenflo, J. Loutschko and Yuri
Luchko, Computation of the Mittag-Leffler function and its
derivatives, Fractional Calculus and Applied Analysis, 5(4)(2002), 491-518.

\vskip.2cm \noindent [76] R. Gorenflo and Yu. F. Luchko,
Operational methods for solving generalized Abel equations of second
kind, Integral Transforms and Special Functions, 5(1997), 47-58.

\vskip.2cm \noindent [77] R. Gorenflo, Yu. Luchko and F. Mainardi,
Analytical properties and applications of Wright function, Fractional
Calculus and Applied Analysis, 2(1999), 383-414.

\vskip.2cm \noindent [78] R. Gorenflo, Yu. Luchko and F. Mainardi,
Wright function as scale-invariant solutions of the diffusion wave
equation, Journal of Computational and Applied Mathematics,
118(2000), 175-191.

\vskip.2cm \noindent [79] R. Gorenflo, Yu. F. Luchko and S.V.
Rogosin, Mittag-Leffler type functions, notes on growth properties
and distribution of zeros, Preprint No. A04-97, Freie Universit\"at
Berlin, Serie A Mathematik, Berlin, 1997.

\vskip.2cm \noindent [80] R. Gorenflo, Yu. F. Luchko and H.M.
Srivastava, Operational method for solving Gauss' hypergeometric
function as a kernel, International Journal of Mathematical and Statistical Sciences, 6(1997), 179-200.

\vskip.2cm \noindent [81] R. Gorenflo and F. Mainardi, Fractional
oscillations and Mittag-Leffler functions, Preprint No. 1-14/96,
Free University of Berlin, Berlin, 1994.

\vskip.2cm \noindent [82] R. Gorenflo and F. Mainardi, The
Mittag-Leffler function in the Riemann-Liouville fractional
calculus, In: A.A. Kilbas (ed) {\it Boundary Value Problems, Special
Functions and Fractional Calculus}, Minsk, 1996, pp.215-225.

\vskip.2cm \noindent [83] R. Gorenflo and F. Mainardi, Fractional
calculus: integral and differential equations of fractional order,
In: {\it Fractals and Fractional Calculus in Continuum Mechanics} (eds. A.
Carpinteri and F. Mainardi), Springer-Verlag, Wien, 1997,
pp.223-276.

\vskip.2cm \noindent [84] R. Gorenflo and R. Rutman, On ultraslow
and intermediate processes, In: P. Rusev, I. Dimovski, V. Kiryakova
(eds) Transform Methods and Special Functions, Sofia, 1994, 61-81,
Science Culture Technology Publ., Singapore, 1995, pp.171-183..

\vskip.2cm \noindent [85] R. Gorenflo and S. Vessella, {\it Abel Integral
Equations: Analysis and Applications}, Lecture Notes in Mathematics
1461, Springer-Verlag, Berlin, 1991.

\vskip.2cm \noindent [86] K.C. Gupta, New relationships of the
H-function with functions of practical utility in fractional
calculus, Ganita Sandesh, 15(2001), 63-66.

\vskip.2cm \noindent [87] K.C. Gupta and R.C. Soni, A unified
inverse Laplace transform formula, functions of practical importance
and H-functions, Journal of the Rajasthan Academy of Physical Sciences, 1(2002),7-16.

\vskip.2cm \noindent [88] I.S. Gupta and L. Debnath, Some
properties of the Mittag-Leffler functions, Integral Transforms and
Special Functions, 18(5)(2007), 329-336.

\vskip.2cm \noindent [89] S.B. Hadid and Yu. F. Luchko, An
operational method for solving fractional differential equations of
an arbitrary real order, Panamerican Mathematical Journal, 6(1996), 57-73.

\vskip.2cm \noindent [90] H.J. Haubold and A.M. Mathai, The
fractional kinetic equation and thermonuclear functions,
Astrophysics and Space Science, 273(2000), 53-63.

\vskip.2cm \noindent [91] H.J. Haubold, A.M. Mathai and R.K.
Saxena, Solution of fractional reaction-diffusion equations in terms
of the H-function, Bulletin of the Astronomical Society, India, 35(2007), 681-689.

\vskip.2cm \noindent [92] R. Hilfer, Fractional diffusion based on
Riemann-Liouville fractional derivative, Journal of Physical Chemistry B,
104(3)(2000), 914-924.

\vskip.2cm \noindent [93] R. Hilfer, On fractional diffusion and
continuous time random walks, Physica A, 329(1-2)(2003), 35-40.

\vskip.2cm \noindent [94] R. Hilfer (ed.), {\it Applications of
Fractional Calculus in Physics}, World Scientific, Singapore, 2000.

\vskip.2cm \noindent [95] R. Hilfer, Y. Luchko and Z. Tomovski,
Fractional method for the solution of fractional differential
equations with generalized Riemann-Liouville fractional derivatives,
Fractional Calculus and Applied Analysis, 12(3)(2009), 299-318.

\vskip.2cm \noindent [96] R. Hilfer and H.J. Seybold, Computation
of generalized Mittag-Leffer function and its inverse in the complex
plane, Integral Transforms and Special Functions, 17(2006), 637-652.

\vskip.2cm \noindent [97] E. Hille and J.D. Tamarkin, On the theory
of linear integral equations, Annals of Mathematics, 31(1930), 479-528.

\vskip.2cm \noindent [98] P. Humbert, Quelques resultants retifs a
la fonction de Mittag-Leffler, C.R. Acad. Sci. Paris, 236(1953),
1467-1468.

\vskip.2cm \noindent [99] P. Humbert and R.P.Agarwal, Sur la
fonction de Mittag-Leffler et quelques unes de ses generalizations,
Bull. Sci. Math.,(Ser.II), 77(1953), 180-185.

\vskip.2cm \noindent [100] M.I. Imanaliev and V.K. Weber, On a
generalization of functions of Mittag-Leffler type and its
application (in Russian), Issled. po Integrao-Diff. Uravneniam v
Kirgizii, 13(1980), 49-59.

\vskip.2cm \noindent [101] A.A. Inayat-Hussain, New properties of
hypergeometric series derivable from Feynman integrals: II, A
generalization of the H-function, Journal of Physics A: Mathematical and General, 20(1987),
4119-4128.

\vskip.2cm \noindent [102] C. Jacques; B. Remillard and R.
Theodorescu, Estimation of Linnik law parameters, Statistics and
Decision, 17(1999), 213-235.

\vskip.2cm \noindent [103] K. Jayakumar, On Mittag-Leffler process,
Mathematical and Computer Modelling, 37(2003), 1427-1434.

\vskip.2cm \noindent [104] K. Jayakumar and R.N. Pillai, The first
order autoregressive Mittag-Leffler process, Journal of Applied
Probability, 30(1993), 462-466.

\vskip.2cm \noindent [105] K. Jayakumar and R.P. Suresh,
Mittag-Leffler distribution, Journal of the Indian Society of
Probability and Statistics, 7(2003), 51-71.

\vskip.2cm \noindent [106] K. Jayakumar and Thomas Mathew,
Autoregressive Mittag-Leffler process, Far East Journal of
Theoretical Statistics, 8(2006), 159-173.

\vskip.2cm \noindent [107] K.K. Jose and R.N. Pillai, Generalized
autoregressive time series models in Mittag-Leffler variables,
Recent Advances in Statistics (ed. P. Yageen Thomas), 1986, pp.
96-103.

\vskip.2cm \noindent [108] R.N. Kalia (ed.), {\it Recent Advances on
Fractional Calculus}, Global Publishing Company, Sauk Rapids,
Minnesota, 1993.

\vskip.2cm \noindent [109] S.L. Kalla, Fractional relations by
means of Riemann-Liouville operator, Serdica, 13(1987), 170-173.

\vskip.2cm \noindent [110] S.L. Kalla and B.N. Al-Saqabi, Summation
of certain infinite series with digamma functions, C.R. Acad. Bulgare
Sci., 41(1988), 15-17.

\vskip.2cm \noindent [111] S.L. Kalla, R.K. Yadav and S.D. Purohit,
On the Riemann-Liouville fractional q-integral operator involving a
basic analogue of Fox's H-function, Fractional Calculus and Applied Analysis,
8(3)(2005), 313-322.

\vskip.2cm \noindent [112] S. Kant and C.L. Koul, On fractional
integral operators, Journal of the Indian Mathematical Society, 56(1991), 97-107.

\vskip.2cm \noindent [113] Anitha Kattuveettil, On Dirichlet
averages, STARS, 2(2008), 78-88.

\vskip.2cm \noindent [114] A.A. Kilbas, Fractional calculus of the
generalized Wright function, Fractional Calculus and Applied Analysis, 8(2)(2005),
113-126.

\vskip.2cm\noindent [115] A.A. Kilbas, Luis Rodrigues-Germa, M.
Saigo, R.K. Saxena and J.J. Trujillo, The Kr\"atzel function and
evaluation of integrals, Computers and Mathematics with
Applications, 58(2009), (to appear).

\vskip.2cm \noindent [116] A.A. Kilbas and Anitha Kattuveettil,
Representations of Dirichlet averages of generalized Mittag-Leffler
functions via fractional integrals and special functions, Fractional
Calculus and Applied Analysis, 11(4)(2008), 471-492.

\vskip.2cm \noindent [117] A.A. Kilbas and Natalia V. Zhukovaskaya,
Euler-type non-homogeneous differential equations with three
Liouville fractional derivatives, Fractional Calculus and Applied Analysis,
12(2)(2009), 204-234.

\vskip.2cm \noindent [118] A.A. Kilbas and M. Saigo, On solutions
of integral equations of Abel-Volterra type, Differential and
Integral Equations, 8(1995), 993-1011.

\vskip.2cm \noindent [119] A.A. Kilbas and M. Saigo, Fractional
integrals and derivatives of Mittag-Leffler type function (Russian),
Dokl. Akad. Nauk Belarusi 39(4)(1995), 22-26.

\vskip.2cm \noindent [120] A.A. Kilbas and M. Saigo, On
Mittag-Leffler type functions, fractional calculus operators and
solution of integral equations, Integral Transforms and Special
Functions, 4(1996), 335-370.

\vskip.2cm \noindent [121] A.A. Kilbas and M. Saigo, Solution in
closed form of a class of linear differential equations of
fractional order (Russian), Differentsialnye Uravnenija 33(1997),
185-204; Translation in Differential Equations 33(1997), 194-204.

\vskip.2cm \noindent [122] A.A. Kilbas and M. Saigo, On the solution
of integral equations of Abel-Volterra type, Differential and
Integral Equations, 8(1995), 993-1011.

\vskip.2cm \noindent [123] A.A. Kilbas and M. Saigo, {\it H-transforms:
Theory and Applications, Analytic Methods and Special Functions},
Vol. 9, Chapman \& Hall, CRC Press, Boca Raton, 2004.

\vskip.2cm \noindent [124] A.A. Kilbas, M. Saigo and R.K. Saxena,
Solution of Volterra integro-differential equations with generalized
Mittag-Leffler function in the kernels, Journal of Integral
Equations and Applications, 14(4)(2002), 377-386.

\vskip.2cm \noindent [125] A.A. Kilbas, M. Saigo and R.K. Saxena,
Generalized Mittag-Leffler function and generalized fractional
calculus operators, Integral Transforms and Special Functions, 15(2004),
31-49.

\vskip.2cm \noindent [126] A.A. Kilbas, M. Saigo and J.J. Trujillo,
On the generalized Wright function, Fractional Calculus and Applied Analysis,
5(4)(2002), 437-460.

\vskip.2cm \noindent [127] A.A. Kilbas, R.K. Saxena, M.Saigo and
J.J. Trujillo, Generalized Wright function as the H-function, in
{\it Analytic Methods of Analysis and Differential Equations}, AMADE 2003,
Cambridge Scientific Publishers, (2006), 117-134.

\vskip.2cm \noindent [128] A.A. Kilbas, R.K. Saxena and J.J.
Trujillo, Kr\"atzel function as a function of hypergeometric type,
Fractional Calculus and Applied Analysis, 9(2006), 109-131.

\vskip.2cm \noindent [129] A.A. Kilbas and N. Sebastian,
Generalized fractional differentiation of the Bessel function of the
first kind, Mathematica Balkanica, New Series, 22(3-4)(2008),

\vskip.2cm \noindent [130] A.A. Kilbas, H.M. Srivastava and J.J.
Trujillo, {\it Theory and Applications of Fractional Differential
Equations}, Elsevier, Amsterdam, 2006.

\vskip.2cm \noindent [131] A.A. Kilbas, J.J. Trujillo and A.A. 
Voroshilov, Cauchy-type problem for diffusion-wave equation with
the Riemann-Liouville partial derivative, Fractional Calculus and Applied Analysis,
8(2)(2005), 403-430.

\vskip.2cm \noindent [132] V. Kiryakova, All the special functions
are fractional differintegrals of elementary functions, Journal of Physics A:
Mathematical and General, 30(1997), 5085-5103.

\vskip.2cm \noindent [133] V. Kiryakova, On two Saigo's fractional
integral operators, Fractional Calculus and Applied Analysis, 9(2)(2006), 159-178.

\vskip.2cm \noindent [134] V. Kiryakova, Some special functions
related to fractional calculus and fractional non-integer order
control systems and equations, Facta Universitatis Ser. Automatic
Control and Robotics, Univ. Nis (2008),

\vskip.2cm \noindent [135] V. Kiryakova, Multiple (multiindex)
Mittag-Leffler functions and relations to generalized fractional
calculus, Journal of Computational and Applied Mathematics, 118(2000), 241-259.

\vskip.2cm\noindent [136] V. Kiryakova, The special functions of
fractional calculus as generalized fractional calculus operators of
some basic functions, Computers and Mathematics with Applications,
58(2009) (to appear).

\vskip.2cm \noindent [137] V.S. Kiryakova, Multiindex Mittag-Leffler
functions related to Gelfond-Leontiev operators and Laplace type
integral transforms. Fractional Calculus and Applied Analysis, 2(1999), 4445-462.

\vskip.2cm \noindent [138] V.S. Kiryakova, {\it Generalized Fractional
Calculus and Applications}, Longman, Harlow, [Pitman Research Notes
in Mathematics, Vol. 301], Wiley New York, 1994.

\vskip.2cm \noindent [139] V.S. Kiryakova, Special functions of
fractional calculus: recent list, results, applications, 3rd IFC
Workshop, FDA 08: Fractional Differentiation and its Applications,
Cankaya University, Ankara, Turkey, 5-7 November 2008b, pp.1-23.

\vskip.2cm \noindent [140] A.N. Kochubei, Fractional order
diffusion, Differential Equations, 26(1990), 485-492 (English
translation from Russian journal, Differentsial'nye Uravneniya).

\vskip.2cm \noindent [141] J.D.E. Konhauser, Biorthogonal
polynomials suggested by the Laguerre polynomials, Pacific Journal of Mathematics,
21(1967), 303-314.

\vskip.2cm \noindent [142] S. Kotz and I.V. Ostrovskii, A mixture
representation of the Linnik distribution, Statistics and
Probability Letters, 26(1996), 61-64.

\vskip.2cm \noindent [143] T.J. Kozubowski, Fractional moment
estimate for Linnik and Mittag-Leffler parameters, Mathematical and
Computer Modeling, special issue: Stable Non-Gaussian Models in
Finance and Econometrics, 34(2001), 1023-1035.

\vskip.2cm \noindent [144] K.R. Lang, {\it Astrophysical Formulae, Vol.
1: Radiation, Gas Processes and  High-energy Astrophysics}, 3rd
edition, Revised edition, Springer-Verlag, New York, 1999a.

\vskip.2cm \noindent [145] K.R. Lang, {\it Astrophysical Formulae, Vol.
2: Space, Time, Matter and Cosmology}, Springer-Verlag, New York,
1999b.

\vskip.2cm \noindent [146] T.A.M. Langlands, Solution of a modified
fractional diffusion equation, Physica A, 367 (2006), 136-144.

\vskip.2cm \noindent [147] J.H. Lavoie, T.J. Osler and R. Tremblay,
Fractional derivatives and special functions, SIAM Review, 18(1976),
240-268.

\vskip.2cm \noindent [148] Gwo Dong Lin, A note on the
characterization of positive Linnik laws, Australian and
New Zealand Journal of Statistics, 40(2001), 17-20.

\vskip.2cm \noindent [149] Gwo Dong Lin, On the Mittag-Leffler
distribution, Journal of Statistical Planning and Inference,
74(1998), 1-9.

\vskip.2cm \noindent [150] S.C. Lim and L.P. Teo, Analytic and
asymptotic properties of multivariate generalized Linnik's
probability densities, arXiv:0903.5344 [math.PR]30 March 2009.

\vskip.2cm \noindent [151] C.F. Lorenzo and T.T. Hartley,
Generalized functions for the factional calculus,
NASA/TP-1999-209424, (1999), 1-17.

\vskip.2cm \noindent [152] C.F. Lorenzo and T.T. Hartley,
Initialization, Conceptualization and Applications in the
Generalized Fractional Calculus, NASA/TP-1998-208415, (1998), 1-107.

\vskip.2cm \noindent [153] C.F. Lorenzo and T.T. Hartley,
Initialized fractional calculus, International Journal of Applied
Mathematics, 3(2000), 249-265.

\vskip.2cm \noindent [154] Yu. F. Luchko, Operational method in
fractional calculus, Fractional Calculus and Applied Analysis, 2(1999), 463-488.

\vskip.2cm \noindent [155] Yu.F.  Luchko, Asymptotics of zeros of
the Wright function, Journal for Analysis and Applications,
19(2)(2000), 583-595.

\vskip.2cm \noindent [156] Yu.F. Luchko, On the distribution of the
zeros of the Wright function. Integral Transforms and Special Functions,
11(2)(2001), 195-200.

\vskip.2cm \noindent [157] Yu. F. Luchko and R. Gorenflo,
Scale-invariant solutins of a partial differential equation of
fractional order, Fractional Calculus and Applied Analysis, 1(1998), 63-78.

\vskip.2cm \noindent [158] Yu. F. Luchko and R. Gorenflo, An
operational method for solving fractional differential equations
with a Caputo derivative, Acta Mathematica Vietnam, 24(1999), 207-234.

\vskip.2cm \noindent [159] Yu. F. Luchko and H.M. Srivastava, The
exact solution of certain differential equations of fractional order
by using fractional calculus, Computational Mathematics and Applications, 29(1995), 73-85.

\vskip.2cm \noindent [160] Yu. F. Luchko and S.B. Yaskubovich,
Operational calculus for the generalized fractional differential
operator and applications, Mathematica Balkanica, New Series, 4(2)(1990),
119-130

\vskip.2cm \noindent [161] Yu. F. Luchko and S.B. Yakubovich, An
operational method for solving some classes of integro-differential
equations (Russian), Differ. Uravn., 30(1994), 269-280; translated
in Differential Equations, 30(1994), 247-256.

\vskip.2cm \noindent [162] F. Mainardi, Fractional calculus: some
basic problems in continuum and statistical mechanics, In: {\it
Fractals and Fractional Calculus in Continuum Mechanics} (eds: A.
Carpinteri and F. Mainardi), Springer-Verlag, Wien, 1997,
pp.291-348.

\vskip.2cm\noindent [163] F. Mainardi, {\it Fractional Calculus and Waves
in Linear Viscoelasticity}, Imperial College Press, London (to
appear in 2010).

\vskip.2cm \noindent [164] F. Mainardi and R. Gorenflo, On
Mittag-Leffler type functions in fractional evolution processes,
Computational and Applied Mathematics, 118(2000), 283-299.

\vskip.2cm \noindent [165] F. Mainardi, R. Gorenflo and A. Vivoli,
Renewal processes of Mittag-Leffler and Wright type, Fractional Calculus and
Applied Analysis, 8(1)(2005), 7-38.

\vskip.2cm\noindent [166] F. Mainardi and R. Gorenflo,
Time-fractional derivatives in relaxation processes a tutorial
survey?, Fractional Calculus and Applied Analysis, 10(3)(2007), 269-308.

\vskip.2cm \noindent [167] F. Mainardi, Yu. Luchko and G. Pagnini,
The fundamental solution of the space-time fractional diffusion
equation, Fractional Calculus and Applied Analysis, 4(2001), 153-192.

\vskip.2cm \noindent [168] F. Mainardi and G. Pagnini,
Mellin-Barnes integrals for stable distributions and their
convolutions, Fractional Calculus and Applied Analysis, 11(4)(2008),
443-456.

\vskip.2cm \noindent [169] F. Mainardi and G. Pagnini, The Wright
functions as solutions of the time-fractional diffusion equation,
Applied Mathematics and Computation, 141(2003), 51-62.

\vskip.2cm \noindent [170] F. Mainardi, G. Pagnini and R.K. Saxena,
Fox H-functions in fractional diffusion, Journal of Computational and Applied Mathematics,
178(2005), 321-331.

\vskip.2cm \noindent [171] O.I. Marichev, {\it Handbook of Integral
Transforms of Higher Transcendental Functions: Theory and
Algorithmic Tables}, Ellis Horwood, Chicherster, 1983.

\vskip.2cm\noindent [172] A.M. Mathai, Some properties of
Mittag-Leffler functions and matrix-variate analogues: A statistical
perspective, Fractional Calculus and Applied Analysis, 13(2010)(to
appear)

\vskip.2cm \noindent [173] A.M. Mathai, A pathway to matrix-variate
gamma and normal densities, Linear Algebra and Its Applications,
396(2005), 317-328.

\vskip.2cm \noindent [174] A.M. Mathai, {\it A Handbook of Generalized
Special Functions for Statistical and Physical Sciences}, Clarendon
Press, Oxford, 1993.

\vskip.2cm \noindent [175] A.M. Mathai, Fox's H-function with matrix
argument, Jornal de Matematica e Estatistica, 1(1979), 91-106.

\vskip.2cm \noindent [176] A.M. Mathai and H.J. Haubold, {\it Special
Functions for Applied Scientists}, Springer, New York, 2008.

\vskip.2cm \noindent [177] A.M. Mathai and R.K. Saxena, {\it The
H-function with Applications in Statistics and Other Disciplines},
Wiley Halsted New York, 1978.

\vskip.2cm \noindent [178] A.M. Mathai, R.K. Saxena and H.J.
Haubold, A certain class of Laplace transforms with application in
reaction and reaction-diffusion equations, Astrophysics and Space
Science, 305(2006), 283-288.

\vskip.2cm \noindent [179] R. Metzler and J. Klafter, The random
walk's guide to anomalous diffusion: a fractional dynamics approach,
Physics Reports, 339(2000), 1-77.

\vskip.2cm \noindent [180] J. Mikusinski, On the function whose
Laplace transform is $\exp(-s^{\alpha}),0<\alpha<1$, Studia Mathematica,
18(1959), 191-198.

\vskip.2cm \noindent [181] K.S. Miller, The Mittag-Leffler function
and related functions, Integral Transforms and Special Functions,
1(1993), 41-49.

\vskip.2cm\noindent [182] K.S. Miller and S.G. Samko, A note on the
complete monotonicity of the generalized Mittag-Leffler function.
Real Analysis Exchange, 23(1997), 753-755.

\vskip.2cm\noindent [183] K.S. Miller and S.G. Sanko, Completely
monotonic functions, Integral Transforms and Special Functions,
12(2001), 389-402.

\vskip.2cm \noindent [184] K.S. Miller and B. Ross, {\it An Introduction
to the Fractional Calculus and Fractional Differential Equations},
Wiley, New York, 1993.

\vskip.2cm \noindent [185] G.M. Mittag-Leffler, Sur l'integrale de
Laplace-Abel, C.R. Acad. Sci. Paris (Ser. II), 136(1902), 937-939.

\vskip.2cm \noindent [186] G.M. Mittag-Leffler, Une generalisation
de l'integrale de Laplace-Abel, C.R. Acad. Sci. Paris (Ser. II),
137(1903), 537-539.

\vskip.2cm \noindent [187] G.M. Mittag-Leffler, Sur la nouvelle
fonction $E_{\alpha}(x)$, C.R. Acad. Sci. Paris, 137(1903), 554-558.

\vskip.2cm \noindent [188] G.M. Mittag-Leffler, Sopra la funzione
$E_{\alpha}(x)$, Rendiconti della Reale Accademia dei Lincei (Ser. v),
13(1904), 3-5.

\vskip.2cm \noindent [189] G.M. Mittag-Leffler, Sur la
representation analytiqie d'une fonction monogene (cinquieme note),
Acta Mathematica, 29(1905), 101-181.

\vskip.2cm \noindent [190] S.S. Nair, Fractional calculus on a
H-function and its special cases, STARS, 2(1)(2008), 50-64.

\vskip.2cm\noindent [191] Seema S. Nair, Pathway fractional
integration operator, Fractional Calculus and Applied Analysis, 12(3)(2009),
237-252.

\vskip.2cm \noindent [192] R. Nigmatulin, On the theory of
relaxation with remnant memory, Phys. Statist. Soc., B, 124(1984),
389-393, translated from Russian.

\vskip.2cm \noindent [193] K. Nishimoto, {\it Fractional Calculus} I
(1984), II( 1987), III (1989), IV (1991), Descartes Press, Koriyama,
Japan.

\vskip.2cm \noindent [194] K. Nishimoto and R.K. Saxena, An
application of Riemann-Liouville operator in the unification of
certain functional relations, Journal of College of Engineering, Nihon
University, Series B, 32(1991), 133-139.

\vskip.2cm \noindent [195] T.F. Nonnenmacher, Fractional integral
and differential equations for a class of L\'evy-type probability
densities, Journal of Physics A: Mathematical and General, 23(1990), L697-L700.

\vskip.2cm \noindent [196] T.F. Nonnenmacher and R. Metzler, On the
Riemann-Liouville fractional calculus and some recent applications.
Fractals, 3(3)(1995), 557-566.

\vskip.2cm \noindent [197] T.F. Nonnenmacher and R. Metzler,
Applications of fractional calculus techniques to problems of
biophysics; In: R. Hilfer (ed), {\it Applications of Fractional Calculus
in Physics}, World Scientific, Singapore, 2001, pp.377-427.

\vskip.2cm \noindent [198] K.B. Oldham and J. Spanier, {\it The
Fractional Calculus}, Academic Press, New York, 1974.

\vskip.2cm \noindent [199] E. Orsingher and L. Beghin,
Time-fractional telegraph equations and telegraph processes with
Brownian time, Probability Theory and Related Fields, 128(2004), 141-160.

\vskip.2cm \noindent [200] E. Orsingher and X. Zhao, The
space-fractional telegraph equation and the related fractional
telegraph process, Chinese Annals of Mathematics, 24B(2003), 45-56.

\vskip.2cm \noindent [201] E. Orsingher and L. Beghin, Fractional
diffusion equations and processes with randomly varying time, The
Annals of Probability, 37(2009), 206-249.

\vskip.2cm \noindent [202] A.G. Pakes, Mixture representation for
symmetric generalized Linnik laws, Statistics and Probability
Letters, 37(1998), 213-221.

\vskip.2cm \noindent [203] R.B. Paris and D. Kaminski, {\it Asymptotics
and Mellin-Barnes Integrals}, Cambridge University Press, Cambridge,
2001.

\vskip.2cm \noindent [204] R.N. Pillai, On Mittag-Leffler functions
and related distributions, Annals of the Institute of Statistical Mathematics, 42(1990),
157-161.

\vskip.2cm \noindent [205] R.N. Pillai and K. Jayakumar, Discrete
Mittag-Leffler distributions, Statistics and Probability Letters,
23(1995), 271-274.

\vskip.2cm \noindent [206] E. Pitcher and W.E. Sewell, Existence
theorems for solutions of differential equations of non-integer
order, Bulletin of the American Mathemaical Society, 44(1938), 100-107.

\vskip.2cm \noindent [207] I. Podlubny, {\it Fractional Differential
Equations}, Academic Press, San Diego, 1999.

\vskip.2cm\noindent [208] T.K. Pog\'any, Integral expressions for
Mathieu-type series whose terms contain Fox's H-function, Applied
Mathematics Letters, 30(2007), 764-769.

\vskip.2cm\noindent [209] T.K. Pog\'any and Z. Tomovski, On
Mathieu-type series whose terms contain a generalized hypergeometric
function ${_pF_q}$ and Meijer's G-function, Mathematical and
Computer Modelling, 47(2008), 952-969.

\vskip.2cm \noindent [210] H. Pollard, The completely monotonic
character of the Mittag-Leffler function $E_{\alpha}(-x)$, Bulletin of the
American Mathematical Society, 54(1948), 1115-1116.

\vskip.2cm\noindent [211] A.Y. Popov, The spectral values of a
boundary values problem and the zeros of Mittag-Leffler functions,
Differential Equations, 38(2002), 642-653.

\vskip.2cm\noindent [212] A.Y. Popov and A.M. Sedletskii,
Distribution of zeros of Mittag-Leffler function, Dokl. Akad. Nauk
390(2003), 165-168.

\vskip.2cm\noindent [213] A.Y. Popov and A.M. Sedletskii,
Distribution of the zeros of the Mittag-Leffler functions, Dokl.
Akad. Nauk, 390(2003), 165-168.

\vskip.2cm \noindent [214] T.R. Prabhakar, A singular integral
equation with a generalized Mittag-Leffler function in the kernel,
Yokohama Mathematical Journal, 19(1971), 7-15.

\vskip.2cm \noindent [215] A.P. Prudnikov, Yu. Brychkov and O.I.
Marichev, {\it Integrals and Series, Vol. 3: More Special Functions},
Gordon and Breach, New York, 1990.

\vskip.2cm \noindent [216] R.K. Raina and M. Bolia, On distortion
theorems involving generalized fractional calculus operators,
Tamkang Journal of Mathematics, 27(2)(1986), 233-241.

\vskip.2cm \noindent [217] R.K. Raina and M. Bolia, The
decomposition structure of a generalized hypergeometric
transformation of convolution type, Computers and Mathematics with Applications,
34(9)(1997), 87-93.

\vskip.2cm \noindent [218] R.K. Raina and C.L. Koul, On Weyl
fractional calculus and H-function transform, Kyungpook Mathematical Journal,
21(2)(1981), 271-279.

\vskip.2cm \noindent [219] R.K. Raina and R.N. Kalia, On
convolution structures for H-function transformations, Analysis
Mathematica, 24(1998), 221-239.

\vskip.2cm\noindent [220] M. Rivero, L. Rodrigues-Germa, J.J.
Trujillo and M. Palar Velasco, Fractional operators and some special
functions, Computers and Mathematics with Applications, 58(2009),
(to appear).

\vskip.2cm \noindent [221] B. Ross (ed), {\it Fractional Calculus and
Its Applications}, Lecture Notes in Mathematics, Vol. 457,
Springer-Verlag, Berlin, 1975 [Proc. Int. Conf. held at University
of New Haven, U.S.A., 1974]

\vskip.2cm \noindent [222] B. Ross and B.K. Sachdeva, The solution
of certain integral equations by means of operators, Amer. Math.
Monthly, 97(1990), 498-503.

\vskip.2cm \noindent [223] A. I. Saichev and G.M. Zaslavsky,
Fractional kinetic equations, solutions and applications, Chaos,
7(4)(1997), 763-764.

\vskip.2cm \noindent [224] M. Saigo, A remark on integral operators
involving the Gauss hypergeometric functions, Math. Rep. Kyushu
Univ., 11(1978), 135-143.

\vskip.2cm \noindent [225] M. Saigo and N. Maeda, More
Generalizations of Fractional Calculus, Transform Methods and
Special Functions, Varna, 1996, IMI-BAS, Sofia (1998), 386-400.

\vskip.2cm \noindent [226] M. Saigo, R.K. Saxena and J. Ram,
Certain properties of operators of fractional integration associated
with Mellin and Laplace transformations, Current Topics in Analytic
Function Theory, pp.291-301, World Scientific Publishing, River
Edge, 1992.

\vskip.2cm \noindent [227] M. Saigo and R.K. Saxena, Applications
of generalized fractional calculus operators in the solution of an
integral equation, Journal of Fractional Calculus, 14(1998), 53-63.

\vskip.2cm \noindent [228] S.G. Samko, A.A. Kilbas and O.I.
Marichev, {\it Fractional Integrals and Derivatives: Theory and
Applications}, Gordon and Breach, New York, 1993, [Translation from
Russian edition, Nauka i Tekhnika, Minsk (1997)].

\vskip.2cm \noindent [229] R.K. Saxena, On fractional integration
operators, Mathematische Zeitschrift, 96(1967), 288-291.

\vskip.2cm \noindent [230] R.K. Saxena, A remark on a paper on
M-series, Fractional Calculus and Applied Analysis, 12(1)(2009), 109-1110.

\vskip.2cm \noindent [231] R.K. Saxena, Certain properties of
generalized Mittag-Leffler function, Conference of the Society for
Special Functions and Their Applications, [Proceedings of the Third
Annual Conference], pp.77-81, Chennai, India, 2002.

\vskip.2cm \noindent [232] R.K. Saxena, Alternative derivation of
the solution of certain integro-differential equations of
Volterra-type, Ganita Sandesh, 17(2003), 51-56.

\vskip.2cm\noindent [233] R.K. Saxena, On a unified fractional
generalization of free electron laser equation, Vijnana Parishad
Anusandhan Patrika, 47(1)(2004), 17-27.

\vskip.2cm \noindent [234] R.K. Saxena, O.P. Gupta and R.K.
Kumbhat, On two dimensional Weyl fractional calculus, Comtes Rendus
de l'Academie Bulgare des Sciences, Tome, 42(1989), 11-14.

\vskip.2cm\noindent [235] R.K. Saxena and S.L. Kalla, On a
fractional generalization of the free electron laser equation, Applied
Mathematics and Computation, 143(2003), 89-97.

\vskip.2cm\noindent [236] R.K. Saxena and S.L. Kalla, Solutions of
Volterra-type integro-differential equations with a generalized
Lauricella confluent hypergeometric function in the kernel, Journal of Mathematical
Sciences, 8(2005), 1155-1170.

\vskip.2cm \noindent [237] R.K. Saxena and S.L. Kalla, On the
solution of certain kinetic equations, Applied Mathematics and Computation,
199(2008), 504-511.

\vskip.2cm \noindent [238] R.K. Saxena and S.L. Kalla, On a
generalized Kr\"atzel function and associated inverse Gaussian
probability distribution, Algebras, Groups and Geometries,
24(3)(2007), 303-324.

\vskip.2cm \noindent [239] R.K. Saxena, S.L. Kalla and V.S.
Kiryakova, Relations connecting multiindex Mittag-Leffler functions
and Riemann-Liouville fractional calculus, Algebras, Groups and
Geometries, 20(2003), 363-385.

\vskip.2cm \noindent [240] R.K. Saxena, A.M. Mathai and H.J.
Haubold, On fractional kinetic equations, Astrophysics and Space
Science, 282(2002), 281-287.

\vskip.2cm \noindent [241] R.K. Saxena, A. M. Mathai and H.J.
Haubold, On generalized fractional kinetic equations, Physica A,
344(2004), 657-664.

\vskip.2cm \noindent [242] R.K. Saxena, A.M. Mathai and H.J.
Haubold, Unified fractional kinetic equations and a fractional
diffusion equation, Astrophysics and Space Science, 290(2004a),
241-245.

\vskip.2cm \noindent [243] R.K. Saxena, A.M. Mathai and H.J.
Haubold, Astrophysical thermonuclear functions for Boltzmann-Gibbs
statistics and Tsallis statistics, Physica A, 344(2004b), 649-656.

\vskip.2cm \noindent [244] R.K. Saxena, A.M. Mathai and H.J.
Haubold, Fractional reaction-diffusion equations, Astrophysics and
Space Science, 305(2006), 289-296.

\vskip.2cm \noindent [245] R.K. Saxena, A.M. Mathai and H.J.
Haubold, Reaction-diffusion systems and nonlinear waves,
Astrophysics and Space Science, 305(2006a), 297-303.

\vskip.2cm \noindent [246] R.K. Saxena, A.M. Mathai and H.J.
Haubold, Solution of generalized fractional reaction-diffusion
equations, Astrophysics and Space Science, 305(2006b), 305-313.

\vskip.2cm \noindent [247] R.K. Saxena, A.M. Mathai and H.J.
Haubold, Solution of fractional reaction-diffusion equation in terms
of Mittag-Leffler functions, International Journal of Scientific Research, 15(2006c), 1-17.

\vskip.2cm \noindent [248] R.K. Saxena, A.M. Mathai and H.J.
Haubold, Solutions of certain fractional kinetic equations and a
fractional diffusion equation, International Journal of Scientific Research 17(2008), 1-8.

\vskip.2cm \noindent [249] R.K. Saxena and K. Nishimoto, Fractional
integral formula for the H-function, Journal of Fractional Calculus, 6(1994),
65-75.

\vskip.2cm \noindent [250] R.K. Saxena and T.F. Nonnenmacher,
Application of H-function to Markovian and non-Markovian chain
models, Fractional Calculus and Applied Analysis, 7(2004), 135-148.

\vskip.2cm \noindent [251] R.K. Saxena and M.A. Pathan, Asymptotic
formulas for unified elliptic-type integrals, Demonstratio
Mathematica, 36(2003), 579-589.

\vskip.2cm \noindent [252] R.K. Saxena and J. Ram, Fractional
integral operators associated with a general class of polynomials,
Indian Acad. Math., 16(1994), 163-171.

\vskip.2cm \noindent [253] R.K. Saxena, J. Ram, S. Chandak and S.L.
Kalla, Unified fractional integral formulae for the Fox-Wright
generalized hypergeometric function, Kuwait Journal of Science and Engineering,
35(1A)(2008), 1-20.

\vskip.2cm \noindent [254] R.K. Saxena, J. Ram and A. R. Chauhan,
Fractional integration of the product of I-function and Appell
function $F_3$, Vijnana Parishad Anusandhan Patrika,
45(4)(2002),345-371.

\vskip.2cm \noindent [255] R.K. Saxena, J. Ram and A.R. Chauhan,
Fractional integration of Appell function $F_3$ associated with a
general class of multivariable polynomials, Proc. Third Annual
Conference of the Society for Special Functions and Their
Applications, Chennai, India, March 4-6, 2002a.

\vskip.2cm \noindent [256] R.K. Saxena, C. Ram and S.L. Kalla, 
Applications of generalized H-function in bivariate distributions, 
Rev. Acad. Canar., 14(1-2)(2002), 111-120.

\vskip.2cm \noindent [257] R.K. Saxena, J. Ram and D.L. Suthar,
Integral formulas for the H-function generalized fractional calculus
-II, South East Asian J. Math.\& Math. Sci., 5(2007), 23-31.

\vskip.2cm \noindent [258] R.K. Saxena, J. Ram and D.L. Suthar,
Fractional calculus of generalized Mittag-Leffler functions, J. Nat.
Acad. Math. (2009) (in press).

\vskip.2cm \noindent [259] R.K. Saxena and M. Saigo, Generalized
fractional Calculus of the H-function associated with Appell
function $F_3$, Journal of Fractional Calculus, 19(2001), 89-104.

\vskip.2cm \noindent [260] R.K. Saxena and M. Saigo, Certain
properties of fractional calculus operators associated with
generalized Wright function, Fractional Calculus and Applied Analysis, 6(2005),
141-154.

\vskip.2cm \noindent [261] R.K. Saxena, R.K. Yadav, S.D. Purohit
and S.L.Kalla Kober fractional q-integral operator of the basic
analogue of the H-function, Rev. Tec. Ing. Univ., 28(2)(2005),
154-158.

\vskip.2cm \noindent [262] M. Saigo and R.K. Saxena, Unified
fractional integral formulas for the multivariable H-function -II,
Journal of Fractional Calculus, 16(1999), 99-110.

\vskip.2cm \noindent [263] M. Saigo and R.K. Saxena, Unified
fractional integral formulas for the multivariable H-function -III,
Journal of Fractional Calculus, 20(1999a), 45-68.

\vskip.2cm \noindent [264] M. Saigo, R.K. Saxena and J. Ram,
Fractional integration of the product of Appell function $F_3$ and
multlvariable H-function, Journal of Fractional Calculus, 27(2005), 31-42.

\vskip.2cm\noindent [265] M.A. Sayed, S.H. Bahiry and W.E. Rasian,
Adomian decomposition method for solving an intermediate fractional
advection-dispersion equation, Computers and Mathematics with
Applications, 58(2009), (to appear).

\vskip.2cm\noindent [266] N. Sazuka, Jun-Ichi Inoue and E. Scalas,
The distribution of first-passage times and durations in FOREX and
future markets, Physica A, 388(2009), 2839-2853.

\vskip.2cm \noindent [267] W.R. Schneider, Stable distributions,
Fox function representation and generalization, In: S. Albeverio, G.
Casati and D. Merilini (eds) {\it Stochastic Processes in Classical and
Quantum Systems}, Lecture Notes in Physics, Vol.262, Springer,
Berlin, 1985, pp.497-511.

\vskip.2cm \noindent [268] W.R. Schneider, Complete monotone
generalized Mittag-Leffler functions, Expositiones Mathematicae,
14(1996), 3-16.

\vskip.2cm \noindent [269] W.R. Schneider and W. Wyss, Fractional
diffusion and wave equations, Journal of Mathematical Physics, 30(1989), 134-144.

\vskip.2cm \noindent [270] N. Sebastian, Certain fractional
integral and differential operators on modified Bessel function of
the first kine, STARS, 2(1)(2008), 50-64.

\vskip.2cm\noindent [271] A.M. Sedletskii, Asymptotic formulas for
the zeros of functions of Mittag-Leffler type, Anal. 20(1994),
117-132.

\vskip.2cm\noindent [272] A.M. Sedletskii, On the zeros of a
function of Mittag-Leffler type, Math. Notes 68(2000), 602-613.

\vskip.2cm \noindent [273] A.M. Sedletskii, Some nonasymptotic
properties of roots of a function of Mittag-Leffler type of order
1/2, Math. Montisnigri, 13(2001), 75-81.

\vskip.2cm \noindent [274] H.J. Seybold and R. Hilfer, Numerical
results for the generalized Mittag-Leffler functions, Fractional Calculus
and Applied Analysis 8(2)(2005), 127-140.

\vskip.2cm \noindent [275] S. Sharma, Fractional differentiation
and fractional integration of the M-series, Fractional Calculus and Applied
Analysis, 11(2008), 187-191.

\vskip.2cm \noindent [276] A.K. Shukla and J.C. Prajapati, Some
remarks on generalized Mittag-Leffler function, Proyecciones
(Chile), 28(1) (2009), 27-34.

\vskip.2cm \noindent [277] A.K. Shukla and J.C. Prajapati, On a
generalized Mittag-Leffler type function and generalized integral
operators, Math. Sci. Res. J., 12(12) (2008), 283-290.

\vskip.2cm \noindent [278] A.K. Shukla and J.C. Prajapati, A
general class of polynomials associated with generalized
Mittag-Leffler function, Integral Transforms and Special Functions,
19(1)(2005), 23-34.

\vskip.2cm \noindent [279] A.K. Shukla and J.C. Prajapati, On a
generalization of Mittag-Leffler function and its properties, Journal
of Mathematical Analysis and its Applications, 336(2007), 797-811.

\vskip.2cm \noindent [280] J.M. Sixdeniers, K.A. Penson and A.I.
Solomon, Mittag-Leffler coherent states, Journal of Physics A: Math. Gen.,
32(1999), 7543-7563.

\vskip.2cm \noindent [281] G.L. Slonimski, About a law of
deformation of polymers, Dokl. Akad. Nauk. SSSR (Physics),
140(1961), 826-829.

\vskip.2cm \noindent [282] I.M. Sokolov, A.V. Chechkin and J.
Klafer, Distributed-order fractional kinetics, Acta Physica Polonica B,
35(2004), 1323-1341.

\vskip.2cm \noindent [283] I.M. Sokolov and J.K. Klafter, From
diffusion in anomalous diffusion: a century after Einstein's
Brownian motion, Chaos, 15(2005), 026103.

\vskip.2cm \noindent [284] H.M. Srivastava, On an extension of
Mittag-Leffler function, Yokohama Math. J., 16(1968), 77-88.

\vskip.2cm \noindent [285] H.M. Srivastava, A certain family of
sub-exponential series, Int. J. Math. Educ. Technol., 25(2)(1994),
211-216.

\vskip.2cm \noindent [286] H.M. Srivastava, Fractional calculus and
its applications, Cubo Matematica Educacional, 5(1)(2003), 33-48.

\vskip.2cm \noindent [287] H.M. Srivastava and P.W. Karlsson,
{\it Multiple Gaussian Hypergeometric Function}, Ellis Horwood, Chichester
and Wiley Halsted, New York.

\vskip.2cm \noindent [288] H.M. Srivastava and S. Owa (eds.)
{\it Univalent Functions, Fractional Calculus and Their Applications},
Ellis Horwood, Chichester and Wiley, New York, 1989.

\vskip.2cm \noindent [289] H.M. Srivastava and R.K. Saxena,
Operators of fractional integration and their applications, Applied
Mathematics and Computation, 118(2001), 1-52.

\vskip.2cm \noindent [290] H.M. Srivastava and R.K. Saxena, Some
Volterra-type fractional integro-differential equations with a
multivariable confluent hypergeometric function in the kernel, Journal of
Integral Equations and Applications, 17(2)(2005), 199-217.

\vskip.2cm \noindent [291] H.M. Srivastava, R.K. Saxena and C. Ram,
A unified presentation of the gamma-type functions occurring in
diffraction theory and associated probability distributions, Applied
Mathematics and Computations, 162(2005), 931-947.

\vskip.2cm\noindent [292] H.M. Srivastava and Z. Tomovski,
Fractional calculus with an integral operator containing a
generalized Mittag-Leffler function in the kernel, Applied Mathematics and
Computations, 211(2009), 198-210.

\vskip.2cm \noindent [293] H.M. Srivastava, S.B. Yakubovich and
Yu. F. Luchko, The convolution method for the development of new
Leibnitz rules involving fractional derivatives and of their
integral analogues, Integral Transforms and Special Functions,
1(2)(1993), 119-134.

\vskip.2cm \noindent [294] A.A. Stanislavsky, Probability
interpretation of the integral of fractional order, Theoretical and
Mathematical Physics, 138(3)(2004), 418-431.

\vskip.2cm \noindent [295] B. Stankovic, On the function of E.M.
Wright, Publ. de l'Institut Mathematique, Nouvelle Ser.,
10(24)(1970), 115-124.

\vskip.2cm \noindent [296] Z. Tomovski, Integral representations of
generalized Mathieu series via Mittag-Leffler type functions, Fractional
Calculus and Applied Analysis, 10(2)(2007), 127-138.

\vskip.2cm \noindent [297] S. Westerlund and L. Ekstam, Capacitor
theory, IEEE Transactions on Dielectrics and Electrical Insulation,
1(1994), 826-839.

\vskip.2cm\noindent [298] T. Mathew and K. Jayakumar,
Generalized Linnik distribution and process, Stochastic Modelling
and Applications, 6(2003), 27-37.

\vskip.2cm \noindent [299] N.S. Tonchev, Finite-size scaling in
systems with strong anisotropy: An analytical example, Communication
of the Joint Institute for Nuclear Research, Dubna, 2005, 1-10.

\vskip.2cm \noindent [300] D.N. Vyas, P.K. Banerji and M. Saigo, On
Dirichlet average and fractional integral of a general class of
polynomials, Journal of Fractional Calculus, 6(1994), 61-64.

\vskip.2cm \noindent [301] E.W. Weisstein, Mittag-Leffler function
from Math world - A Wolfram Web Resource, 10, March 2003.

\vskip.2cm \noindent [302] K. Weron and M. Kotulski, On the
Cole-Cole relaxation function and related Mittag-Leffler
distributions, Physica A, 232(1996), 180-188.

\vskip.2cm \noindent [303] A. Wiman, \"Uber den Fundamental satz in
der Theorie der Funcktionen, $E_{\alpha}(x)$, Acta Mathematica, 29(1905),
191-201.

\vskip.2cm \noindent [304] A. Wiman, \"Uber die Nullstellun der
Funktionen $E_{\alpha}(x)$, Acta Mathematica, 29(1905a), 217-234.

\vskip.2cm \noindent [305] E.M. Wright, On the coefficients of
power series having exponential singularities, J. London Math. Soc.,
8(1933), 71-79.

\vskip.2cm \noindent [306] E.M. Wright, The asymptotic expansion of
the generalized Bessel function, Proc. London Math. Soc.,(2)
38(1934), 257-270.

\vskip.2cm \noindent [307] E.M. Wright, The asymptotic expansion of
the generalized hypergeometric function, J. London Math. Soc.,
10(1935), 286-293.

\vskip.2cm \noindent [308] E.M. Wright, The asymptotic expansion of
the integral functions defined by Taylor series, Philos. Trans. Roy.
Soc. London Ser. A, 238(1940), 423-451.

\vskip.2cm \noindent [309] E.M. Wright, The asymptotic expansion of
the generalized hypergeometric function, Proc. London Math. Soc.,
46(1940a), 389-408.

\vskip.2cm \noindent [310]  E.M. Wright, The generalized Bessel
function of order greater than one, Quarterly Journal of Mathematics Oxford, 
11(1940b), 36-48.

\vskip.2cm \noindent [311] W. Wyss, Fractional diffusion equations,
Journal of Mathematical Physics, 27(1986), 2782-2785.

\vskip.2cm \noindent [312] W. Wyss, The fractional Black-Scholes
equation, Fractional Calculus and  Applied Analysis, 3(1)(2000), 51-61.

\vskip.2cm \noindent [313] R. Yu and H. Zhang, New function of
Mittag-Leffler type and its application in the fractional
diffusion-wave equation, Chaos, Solitons and Fractals, 30(2006),
946-955.

 \bye